\documentclass[preprint]{elsarticle}

\usepackage{jmm3d}

\begin{document}

\begin{frontmatter}

\title{Numerical geometric acoustics: an eikonal-based approach for modeling sound propagation in 3D environments}

\author[1]{Samuel F. Potter\corref{cor1}}\ead{sfp@cims.nyu.edu}
\author[2]{Maria K. Cameron}\ead{mariakc@umd.edu}
\author[3]{Ramani Duraiswami}\ead{ramanid@umd.edu}

\cortext[cor1]{Corresponding author}

\affiliation[1]{
  organization={Courant Institute of Mathematical Sciences, New York University},
  addressline={251 Mercer St.},
  postcode={10012},
  city={New York},
  country={USA}
}

\affiliation[2]{
  organization={University of Maryland Department of Mathematics},
  city={College Park},
  state={MD},
  country={USA}
}

\affiliation[3]{
  organization={University of Maryland Department of Computer Science},
  city={College Park},
  state={MD},
  country={USA}
}

\begin{abstract}
  We present algorithms for solving high-frequency acoustic scattering problems in complex domains. The eikonal and transport partial differential equations from the WKB/geometric optic approximation of the Helmholtz equation are solved recursively to generate boundary conditions for a tree of eikonal/transport equation pairs, describing the phase and amplitude of a geometric optic wave propagating in a complicated domain, including reflection and diffraction. Edge diffraction is modeled using the uniform theory of diffraction. For simplicity, we limit our attention to domains with piecewise linear boundaries and a constant speed of sound. The domain is discretized into a conforming tetrahedron mesh. For the eikonal equation, we extend the jet marching method to tetrahedron meshes. Hermite interpolation enables second order accuracy for the eikonal and its gradient and first order accuracy for its Hessian, computed using cell averaging. To march the eikonal on an unstructured mesh, we introduce a new method of rejecting unphysical updates by considering Lagrange multipliers and local visibility. To handle accuracy degradation near caustics, we introduce several fast Lagrangian initialization algorithms. We store the dynamic programming plan uncovered by the marcher in order to propagate auxiliary quantities along characteristics. We introduce an approximate origin function which is computed using the dynamic programming plan, and whose 1/2-level set approximates the geometric optic shadow and reflection boundaries. We also use it to propagate geometric spreading factors and unit tangent vector fields needed to compute the amplitude and evaluate the high-frequency edge diffraction coefficient. We conduct numerical tests on a semi-infinite planar wedge to evaluate the accuracy of our method. We also show an example with a more realistic building model with challenging architectural features. Finally, we demonstrate a simple approach to extending the method to handle nonconstant speeds of sound by modifying the semi-Lagrangian updates to account for a varying speed.
\end{abstract}


\end{frontmatter}

\section{Introduction}\label{sec:introduction}

\subsection{Background and motivation: precomputed room acoustics}

The viscosity solution of the eikonal partial differential equation (PDE) models the first-arrival time of a geometric optic (GO) wavefront propagating in a domain. Despite the obvious connection, numerically solving this PDE has found relatively little application in computer graphics, although there are exceptions~\cite{ihrke2007eikonal}. There are a variety of reasons for this, a primary one being the fact that the very short wavelength of visible light renders the effect of diffraction on visual perception insignificant. Nevertheless, there are a wide variety of applications of the eikonal equation, including high-frequency wave propagation (typically in geophysics), computer vision, stereolithography, optimal control, robotics, computational geometry, and more~\cite{sethian1999level}.

Borrowing from the success of raytracing in computer graphics, acousticians interested in simulating sound propagation in complicated built environments have been successful repurposing raytracing methods from computer graphics to the simulation of sound propagation, leading to \emph{geometric acoustics} methods~\cite{savioja2015overview}. These methods have the advantage of being fast enough to be used in realtime. They also provide a reasonable degree of verisimilitude, and are easily interpreted and explained, making it feasible for sound designers to creatively interpret and modify the results. Moreover, incorporating diffraction effects into these models is challenging, and usually involve introducing approximations~\cite{mannall2022perceptual}. Methods based on solving the wave or Helmholtz equations to simulate wave propagation are also used~\cite{hamilton2016finite,gumerov2021fast}. These have the advantage of being more ``physically correct'', but are quite slow---the runtime complexity of the simulation depends at least quadratically on the highest frequency simulated, resulting in solvers which are impractical for use in room acoustics simulations where the product of the wavenumber $k$ and the diameter of the domain is too large. For rooms and buildings with a diameter greater than 10 m, a frequency much higher than 3000 Hz becomes impractical (recall that the audible spectrum typically spans from 20 Hz to 20 kHz).

An important line of research in room acoustic modeling in recent years has been \emph{precomputed room acoustics}. Since physical modeling of sound propagation is simply too expensive to be done in real-time, researchers have developed approaches to the offline precomputation of salient perceptual quantities which can then be accessed as if from a lookup table in order to model the \emph{room impulse response} (i.e., the domain Green's function of the wave equation)~\cite{raghuvanshi2018parametric}. It is interesting to note that these methods essentially extract GO quantities (such as the first arrival time, the first arrival amplitude, early reflection times, etc.) from solutions of the wave equation; consequently, human perception and the attendant psychoacoustic parameters have a GO flavor~\cite{kuttruff2016room}.

\subsection{Prior art: fast eikonal solvers}

Our goal is to investigate methods for simulating wave propagation in complex 3D environments based on solving the eikonal PDE. Ultimately, our goal is to use such a solver as the basis for a system for precomputed room acoustics---this work is a step in that direction. The key observation is that computing the viscosity solution of the eikonal equation naturally incorporates diffraction around obstacles. In room acoustics, diffraction plays an important role: we can hear around obstacles, but not see around them. At the same time, as we mentioned, many psychacoustic quantities are distinctly GO in nature---our view is that computing these quantities directly instead of extracting them after the fact from solutions of the wave equation may be advantageous.

Unfortunately, existing eikonal solvers are not adequate to this task. The classical fast algorithm for solving the eikonal PDE---the \emph{fast marching method} (FMM)---is an $O(h)$ Eulerian solver based on using upwind finite differences to approximate solutions to the eikonal PDE in optimal $O(N \log N)$ time on a regular grid of $N$ nodes with spacing $h$~\cite{sethian1996fast}. A semi-Lagrangian version of the FMM was developed slightly earlier in the optimal control community~\cite{tsitsiklis1995efficient}. Another family of eikonal solvers consists of \emph{fast sweeping methods}~\cite{zhao2005fast}.

Methods based on the classic first-order upwind finite difference discretization are not ideal for propagating GO waves for two major reasons:
\begin{itemize}
\item To compute the amplitude of GO waves consistently (see \Cref{sec:high-frequency-helmholtz}), it is necessary to compute the Hessian of the eikonal consistently. The FMM is roughly $O(h)$ accurate---in fact, the gradient of the eikonal computed by the FMM is also frequently $O(h)$, although the conditions under which this occurs are somewhat unobvious~\cite{potter2021jet}. The Hessian of the eikonal is not approximated consistently by the FMM (its error is $O(1)$).
\item Edge-diffraction plays an important role in high-frequency acoustic wave propagation~\cite{torres2001computation}. If we discretize a domain using a regular grid, it becomes difficult to accurately model reflection and diffraction from a non-grid aligned boundary.
\end{itemize}

Although fast eikonal solvers which are higher-order than $O(h)$ have been developed, they typically use high-order WENO finite difference stencils---hence, they assume a regular grid and use wide stencils~\cite{luo2014high}. Some eikonal solvers with compact finite-difference stencils have been developed, but their order of accuracy is limited to $O(h^2)$~\cite{benamou2010compact}. Note that these solvers are \emph{free-space} eikonal solvers, and have been developed in the context of propagating GO waves on a domain with a rapidly varying speed of sound but without explicitly modeled boundaries, which is the case in exploration geophysics.

In fact, the FMM is only $O(h \log \tfrac{1}{h})$~\cite{zhao2005fast}---the order of accuracy is degraded by the formation of caustics, which are the loci of points where the curvature of the GO wavefront is infinite. Point sources are clearly caustics, but caustics also form at points of edge-diffraction where the GO wave rarefies as it spreads into the shadow zone behind an obstacle. Because of the upwind nature of the FMM and similar solvers, any such inaccuracies will pollute the error in the eikonal beyond a caustic. To this end, \emph{factored} eikonal equations have been developed which fix this deficiency, but generally only in the context of point source singularities~\cite{fomel2009fast}.

To accurately model acoustic scattering, it is necessary to use an unstructured mesh which conforms to the boundary so that the precise locations of the singular parts of the boundary are retained. The fast marching method was extended to triangle meshes, primarily to compute geodesics, but the order of convergence has generally been limited to $O(h)$~\cite{kimmel1998computing,sethian2000fast}. A \emph{gradient-augment} fast marching method on an unstructured mesh with apparent $O(h^2)$ accuracy was investigated but its convergence properties were not carefully explored~\cite{sethian2000fast}. One drawback of these methods is that they make use of complicated ``splitting schemes'' in order to preserve the \emph{causality} (in the sense of Dijkstra's algorithm) of the solver. The fast sweeping method has also been adapted to unstructured meshes, but it is clear in this instance that the iteration count suffers for complicated domains, where a large number of sweeps are required~\cite{qian2007fast}. An alternative iterative scheme known as the \emph{fast iterative method} has also been adapted to unstructured meshes~\cite{fu2013fast}.

To summarize, a fast eikonal solver which suits our purposes has the following properties:
\begin{itemize}
\item At least $O(h)$ accurate approximation of the eikonal Hessian.
\item Runs in $O(N \log N)$ time, where $N$ is the number of degrees of freedom of the spatial discretization (e.g., vertices of an unstructured mesh).
\item Compatible with unstructured meshes conforming to the boundary of a complicated, nonconvex domain.
\end{itemize}
Towards this end, we previously developed the \emph{jet marching method} (JMM)~\cite{potter2021jet,potter2021numerical}, which is a framework for building higher-order semi-Lagrangian eikonal solvers with compact stencils which directly march the \emph{jet} of the eikonal. Here, a jet is the eikonal and some number of its derivatives at a fixed point; that is, the data for a multivariate Taylor polynomial approximating the eikonal. In that work, we developed a variety of JMMs in 2D on regular grids which march the eikonal and its gradient, and compute the amplitude using paraxial raytracing~\cite{popov2002ray}. To ensure locality, Hermite interpolation is used over cells and for the discretizations used for the local updates. The JMM was later applied to compute the quasipotential, an important function in the theory of large deviations~\cite{paskal2022efficient}.

\subsection{Prior art: Eulerian geometric optics, raytracing, and the image method}\label{ssec:prior-art-EGO-and-RT}

In addition to solving the eikonal PDE, we must also compute the amplitude of a GO wave---and for multiple arrivals due to reflection and edge-diffraction, to boot. There are a variety of interesting algorithms for computing multiple arrivals based on solving the eikonal PDE~\cite{engquist2003computational}. Most similar to our approach is the line of research which goes under the name \emph{Eulerian geometric optics}~\cite{benamou2003introduction}, or \emph{big ray tracing} (BRT)~\cite{benamou1996big}. In this method, rays are traced to partition the domain into regions in which the eikonal equation can be solved using finite differences. It was combined with an unstructured mesh eikonal solver based on solving an unsteady version of the eikonal equation~\cite{abgrall1999big}. The main difference between BRT and our approach is that instead of taking a combined Eulerian-Lagrangian approach, we take a semi-Lagrangian approach. This allows us to make our solver more local, which streamlines the use of unstructured meshes. We avoid the need to partition the domain into ``big rays'' by attempting to localize the shadow and reflection boundaries of geometric optics directly---indeed, we attempt to approximately determine visibility information \emph{locally} as opposed to globally, as in raytracing (see the next paragraph). It also makes it possible to design dynamic programming algorithms for applying amplitude boundary conditions based on the uniform theory of diffraction (UTD).

The obvious alternative to our approach is raytracing, where paths are traced by integating the raytracing ODEs directly. This can be done in one of two ways: either by solving initial value problems (IVPs), or two-point boundary value problems (BVPs). For IVPs, ensuring uniform ray coverage is important, and becomes more complicated in the presence of diffraction and reflection. For two-point BVPs, handling diffraction and reflection becomes complicated---e.g., to trace a reflected ray to a particular point, it becomes necessary to optimize the point of reflection to find the minimum travel time, which ultimately leads to the approach taken in the current paper. At the same time, for both IVPs and BVPs, it is necessary to check for each ray whether it intersects the domain boundary---doing so quickly requires a spatial data structure. Available raytracing libraries universally assume straight ray trajectories. To check the visibility of curved trajectories is feasible but then requires special purpose library. All this is to say that while a ``global'' approach to propagating fields of rays based on solving IVPs or two-point BVPs is certainly feasible, it is not clear ahead of time how it compares to the approach proposed in this work. A thoroughgoing comparison is outside the scope of this paper.

We also mention that our work can be viewed as a PDE-based version of the method of images~\cite{kuttruff2016room} or the beam tracing method~\cite{heckbert1984beam}. Beam tracing operates by reflecting a point source over a flat section of boundary in order to generate a secondary source which gives rise to a specular reflection. The obvious shortcoming of such approaches is that it is not straightforward to extend them to handle nonconstant speeds of sound, diffraction, or reflection from curved surfaces. Each of these shortcomings stems from the same physical phenomenon: while the caustic set for a specular reflection is a mirror copy of the source, this is not true for these more complicated cases (reflecting a point over a curved surface produces a one-dimensional locus of images). Instead, our method simply requires one to determine the BCs for scattered fields, which can be done comparatively straightforwardly (see \Cref{sec:eikonal} for more details).

\subsection{Prior art: fast Huygens sweeping method and Gaussian beam summation}

We briefly mention two other directions of research related to high-frequency wave propagation. One approach is the fast Huygens sweeping method~\cite{luo2014high}, which can handle high-frequency multipathing in inhomogeneous media. The key assumptions here are two-fold: 1) there is a primary direction of propagation, and 2) the speed of sound has extremely, even discontinuous, variation. These assumptions are often made in seismic wave propagation---layers of rock and voids lead to a strongly varying speed of sound with discontinuities. At the same time, propagation in the $z$ direction can be assumed to simplify modeling assumptions. Clearly, these assumptions are not appropriate for room acoustics---the speed of sound can be assumed to be constant or slowly varying, and waves can propagate in any direction.

Another direction which has seen recent interest is the Gaussian beam method~\cite{vcerveny1982computation,popov2002ray}. This approach was also developed for seismic wave propagation. The idea is to replace GO rays with rays with a complex ray parameter. The wavefield associated with the ray now has a geometric spreading which varies laterally in the plane normal to direction of ray propagation, unlike in GO, where rays remain ``locally plane'' as they propagate. The result is that Gaussian beams are valid asymptotically as they pass through caustics, which is not true of GO rays. Since GO and Gaussian beams have the asymptotic error, it has been observed that a hybrid method based on GO and Gaussian beams could be developed, where GO is used in regions that are caustic-free, and Gaussian beams are used near caustics~\cite{motamed2015wavefront}. It should be possible to develop such a hybrid approach based on what we present in this article, but we defer it for a future work. The effect of caustics in room acoustics is basically second-order: there are some perceptual phenomena which can only be addressed by modeling caustics, but they are rare and result in only a modest qualitative change in the perception of sound. An example would be a whispering gallery. We view this as a reasonable starting point for our developments.

\subsection{The focus of this work}

\begin{figure}
  \includegraphics[width=\linewidth]{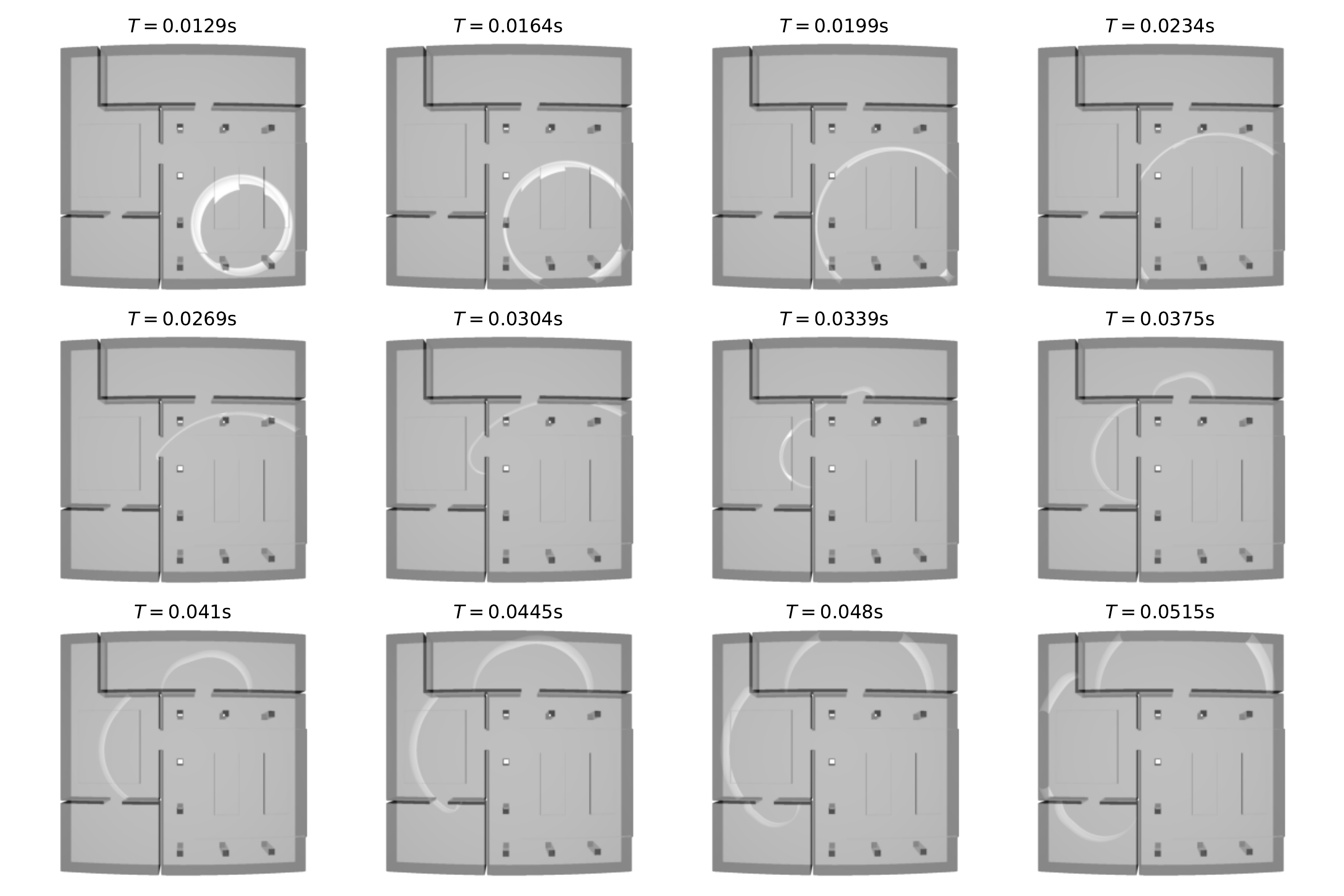}
  \caption{Snapshots of the first arrival time at different time points for the building test problem of \Cref{ssec:building}. A point source is placed in the lower right room. At successive instants, the GO wavefront can be seen to diffract around features of the room and enter into different rooms. The images were constructed using the method described in \Cref{ssec:extracting-level-sets}.}\label{fig:wavefront-plots-top}
\end{figure}

Our goal is to introduce a set of tools which can be used to compute a high-frequency approximation of an acoustic pressure field in a complicated environment, such as the interior of a building. We limit our focus to nonconvex domains with simple, piecewise linear boundaries. We also assume a constant speed of sound. Our strategy is to compute the eikonal for each branch of a GO wave propagating through a domain, including reflection from the linear facets of the boundary and diffraction from the singular parts of the boundary. We do this by recursively generating a tree of eikonal problems. We show how dynamic programming algorithms can be used to propagate information along characteristics in a semi-Lagrangian fashion, which is used to apply boundary conditions, determing visibility, and compute the spreading factor in a post-processing phase. See \Cref{fig:wavefront-plots-top} for several snapshots of the \emph{direct} GO wave (not including reflections) due to a point source placed in a building test model.

There are two key ways in which this model is limited and can be generalized: by allowing a smoothly varying speed of sound and by allowing a curved boundary. If the speed of sound is allowed to vary smoothly, then the minimization problems which need to be solved to do semi-Lagrangian updates become somewhat more complicated, although this does not affect the overall asymptotic complexity of our solver. At the same time, both a varying speed of sound and curved boundary introduce the possibility of freespace caustics, which we largely ignore in this work. It is possible to extend the framework presented here to allow for multiple arrivals due to caustics. Throughout, we have made choices which can be extended in a natural way to account for these generalizations.

The key innovations of this work are as follows:
\begin{itemize}
\item We introduce a fast marcher which computes the eikonal and its gradient with formally $O(h^2)$ accuracy and its Hessian with $O(h)$ accuracy on an unstructured volumetric tetrahedron mesh. To this end, we present data structures and algorithms alternative to previous splitting section schemes which cache and merge previous updates from the valid front.
\item We propagate a variety of derived quantities downwind using the cached dynamic programming plan uncovered by the marcher. This allows us to compute the amplitude in a post-processing phase.
\item We present additional algorithms which allow the boundary conditions (BCs) for the edge-diffracted amplitude to be applied in a semi-Lagrangian setting.
\end{itemize}
We test and evaluate each of these developments using extensive numerical experiments. For a groundtruth test with an analytically available solution, we consider the problem of edge-diffraction from a semi-infinite planar wedge. For a more realistic test case, we simulate sound propagation in a simplified building geometry, typical of architectural environments.

\section{High-frequency approximation of the Helmholtz equation}\label{sec:high-frequency-helmholtz}

In this section, we recapitulate basic facts about the high-frequency approximation of the Helmholtz equation. This is the GO approximation, which has some fundamental limitations, such as the fact that the asymptotic series does not converge at caustics. Nevertheless, since it is the de facto model used in room acoustics, we use it here. Alternatives such as asymptotics based on Babich's ansatz are available~\cite{lu2016babich}, but we leave their exploration for a future work.

\subsection{Geometric acoustics}

The acoustic pressure satisfies the wave equation:
\begin{equation}
  \label{eq:wave-equation}
  p_{tt}(t, \mx) - c(\mx)^2 \Delta p(t, \mx) = 0,
\end{equation}
where $c$ is the speed of sound, which in the general case can be assumed to vary with position, and where $\Delta = \partial_x^2 + \partial_y^2 + \partial_z^2$ is the Laplacian operator on $\RR^3$~\cite{kuttruff2016room}. We assume that $p$ is time-harmonic:
\begin{equation}
  p(t, \mx) = u(\mx) e^{i\omega t},
\end{equation}
where $u$ is the Fourier coefficient of $p$ for frequency $\omega > 0$, and satisfies the Helmholtz equation:
\begin{equation}\label{eq:helmholtz-equation}
  \Delta u(\mx) + \frac{\omega^2}{c(\mx)^2} u(\mx) = 0.
\end{equation}
If we make the GO ansatz---i.e.., the WKB ansatz--- and assume that $u$ can be factored into the product of a complex prefactor $\alpha$, and an exponential which is governed by spatially-varying (real) phase function $\tau$:
\begin{equation}\label{eq:WKB-ansatz}
  u(\mx) = \alpha(\mx) e^{i\omega\tau(\mx)}.
\end{equation}
Plugging ansatz~\eqref{eq:WKB-ansatz} into \eqref{eq:helmholtz-equation} and grouping terms according to powers of $\omega$, we arrive at:
\begin{equation}
  \parens{\Delta\alpha + (i\omega) \parens{\alpha\Delta\tau + 2\grad\alpha\cdot\grad\tau} + \parens{i\omega}^2 \alpha \parens{\grad\tau\cdot\grad\tau - \frac{1}{c^2}}} e^{i\omega\tau} = 0.
\end{equation}
This is satisfied if each of the following PDEs is satisfied:
\begin{align}
  &\norm{\grad\eik(\mx)} = \slow(\mx) := 1/\speed(\mx), \label{eq:eikonal-equation} \\
  &\alpha(\mx) \Delta\eik(\mx) + 2 \grad\alpha(\mx)\cdot\grad\eik(\mx) = 0, \label{eq:amplitude-equation} \\
  &\Delta \alpha(\mx) = 0. \label{eq:laplace-equation}
\end{align}
The spatial phase field $\tau$ is known as \emph{the eikonal}. It satisfies \eqref{eq:eikonal-equation}, which is the \emph{the eikonal equation}. The second PDE is a transport equation which we refer to as \emph{the amplitude equation}, since $\alpha$ corresponds to the amplitude of $u$. The final equation is Laplace's equation: the degree to which $\alpha$ satisfies it can be seen as one measure of the validity of the high-frequency ansatz. For example,  $\Delta \alpha(\mx) \neq 0$ if $\mx$ is a caustic point. Note that we will only solve \eqref{eq:eikonal-equation} and \eqref{eq:amplitude-equation} to obtain $\tau$ and $\alpha$.

If $c$ is
constant, then the rays are straight lines. For a ray which does not pass through a caustic, the eikonal and length of the ray are related by:
\begin{equation}
  c \tau = r.
\end{equation}
Furthermore, if $k$ denotes the wavenumber, we have:
\begin{equation}
  \qquad k = \omega/c, \qquad \omega \tau = \omega r / c.
\end{equation}
As an example, for a point source at $\mxsrc = 0$ in free space, we have:
\begin{equation}
  r = \norm{\mx}, \qquad u(\mx) = \frac{e^{i k r}}{ r}, \qquad \parens{\frac{\omega^2}{c^2} + \Delta} u(\mx) = 4 \pi \delta (\mx).
\end{equation}
Here, $u$ is a scaled Green's function for the Helmholtz equation. From this, we can recover $\tau$ and $\alpha$:
\begin{equation}
  \alpha(\mx) = \frac{\alpha_0}{\norm{\mx}}, \qquad \tau(\mx) = \frac{\norm{\mx}}{c}
\end{equation}
Clearly, the eikonal equation is satisfied. Additionally, we can see that:
\begin{equation}
  \alpha\Delta\tau + 2 \grad\alpha\cdot\grad\eik = \alpha \Delta \tau + 2 \grad\alpha\cdot\grad\eik + \eik\Delta\alpha = \Delta (\alpha \tau) = \Delta(\alpha_0/c) = 0
\end{equation}
In this case, $\alpha$ is also a scaled Green's function for the Laplace equation---hence, $\Delta\alpha(\mx) = \delta(\mx)$, satisfying \eqref{eq:laplace-equation} everywhere except $\mx = \m{0}$, where a caustic resides ($\alpha(\m{0}) = \infty$).

Note that higher-order amplitude corrections are available by using a series form of the WKB ansatz~\cite{mcnamara1990introduction}. This results in a sequence of transport equations, each coupled to the last. The ratio of two successive amplitude terms can be used to check whether the asymptotic series diverges~\cite{popov2002ray}. We do not explore using these higher-order terms in this work.

\subsection{Multiple arrivals}

The eikonal describes the phase of a GO wave at each point, but acoustic waves in a complicated domain will exhibit complicated scattering phenomena, with multiple waves passing through a single point. The model in the previous section does not capture this behavior on its own. More generally, it is necessary to solve multiple instances of \eqref{eq:eikonal-equation} and \eqref{eq:amplitude-equation} and superimpose them. For instance, if $\big\{(\tau_j, \alpha_j)\big\}$ is a set of eikonals and amplitudes, with each $j$ indexing a separate branch (possibly the direct field, a reflection, or an edge-diffracted field), then the total field is given by:
\begin{equation}
  u(\mx) = \sum_j \alpha_j(\mx) e^{i\omega\tau_j(\mx)}.
\end{equation}
It is then a question of how to determine the BCs for these PDEs. If we first solve these PDEs for a point source located somewhere in the domain, this \emph{direct field} should give rise to a recurrence generating a tree of \emph{scattered fields}. We discuss how to generate appropriate BCs for the scattered eikonal and amplitude throughout the remainder of this document.

Classical GO theory predicts specularly refracted and reflected waves satisfying Snell's law---only reflected waves (i.e., specularly reflected waves) concern us here. Keller's \emph{geometric theory of diffraction} (GTD) used a generalized Fermat's principle which allows constraints inferred from the environment, predicting edge-diffracted waves~\cite{keller1962geometrical}. While the BCs for the edge-diffracted eikonal are straightforward, BCs for the edge-diffracted amplitude are not. There are several approaches to modeling edge-diffracted fields, and we opt to employ the asymptotic diffraction coefficients from the \emph{uniform theory of diffraction} (UTD)~\cite{kouyoumjian1974uniform,mcnamara1990introduction,molinet2011acoustic}. In later sections, we show how to apply UTD in a semi-Lagrangian setting.

We also remark that computing all of these fields is fraught, and is not the goal of this work. In computational room acoustics, the direct arrival augmented by several early reflections provides a significant amount of perceptual information---our goal is to initiate the development of a set of tools which can be used systematically in different applications in the future. The most effective way to do this naturally depends on the application and is outside the scope of this work.

\section{The eikonal equation}\label{sec:eikonal}

We now elaborate some of the important properties of the eikonal PDE and its solutions in this section. In particular, we discuss Fermat's principle, which gives the exact solution in a Lagrangian framework, how the curvature of a GO wavefront can be easily extracted in an Eulerian setting, and how to generate BCs for the specularly reflected and edge-diffracted eikonal.

\subsection{Solution by Fermat's principle}

Let $\domain \subseteq \RR^3$ be a domain, and let $\boundary\domain$ be its boundary. We do not require $\domain$ to be convex. We denote the \emph{eikonal} by $\eik:\domain\to\RR$. It satisfies the \emph{eikonal PDE}:
\begin{equation}
  \begin{cases}
    \|\grad\eik(\mx)\| = \slow(\mx), & \mx \in \Omega, \\
    \eik(\mx) = \eik_0(\mx), &\mx \in \BCset \subseteq \Omega.
  \end{cases}
\end{equation}
The set $\BCset$ specifies the points which have Dirichlet boundary conditions, with the boundary values given by $\tau_0:\BCset\to\RR$. The characteristics of the eikonal equation are the rays of geometric optics. Let $\mpsi : [\sigma_0, \sigma_1] \to \Omega$ be a ray segment such that:
\begin{equation}
  \mpsi'(\sigma) = \frac{\grad\eik(\mpsi(\sigma))}{s(\mpsi(\sigma))}, \qquad \sigma_0 \leq \sigma \leq \sigma_1,
\end{equation}
where we assume that $\mpsi$ is parametrized by arc length (i.e., $\|\mpsi'\| \equiv 1$). Evaluating $\tau$ along $\mpsi$ and differentiating with respect to $\sigma$ gives:
\begin{equation}
  \frac{d}{d\sigma} \tau(\mpsi(\sigma)) = \mpsi'(\sigma)^\top \grad\eik(\mpsi(\sigma)) = s(\mpsi(\sigma)).
\end{equation}
Integrating from $\sigma_0$ to $\sigma_1$ gives:
\begin{equation}\label{eq:integrate-tau}
  \tau(\mpsi(\sigma_1)) = \tau(\mpsi(\sigma_0)) + \int_{\sigma_0}^{\sigma_1} s(\mpsi(\sigma)) d\sigma.
\end{equation}
We can determine $\mpsi$ by using \emph{Fermat's principle}:
\begin{equation}
  \eik(\mx) = \min_{\substack{\mpsi : [0, 1] \to \Omega \\ \mpsi(0) \in \BCset \\
      \mpsi(1) = \mx}} \Set{\eik(\mpsi(0)) + \int_0^1 \slow(\mpsi(\sigma))\norm{\mpsi'(\sigma)}d\sigma}.\label{eq:fermat}
\end{equation}
Here, we allow the parametrization of $\mpsi$ to vary, unlike in \eqref{eq:integrate-tau}, where we assume an arc length parametrization. To solve \eqref{eq:fermat}, we fix the endpoints of $\mpsi$, allow it to vary, and minimize the corresponding functional. This can be done using the calculus of variations, leading to the raytracing ODEs (i.e., the Euler-Lagrange equations for the eikonal equation). Alternatively, we can discretize Fermat's principle directly and numerically solve the resulting minimization problem---this is the approach we take.

The eikonal PDE can be rewritten $\grad\eik(\mx)^\top\grad\eik(\mx) = \slow(\mx)^2$, the gradient of which is:
\begin{equation}
  \hess\eik(\mx)\grad\eik(\mx) = \slow(\mx)\grad\slow(\mx),
\end{equation}
where $\hess\eik$ is the Hessian of the eikonal. Hence:
\begin{equation}
  \frac{d}{d\sigma} \grad\eik(\mpsi(\sigma)) = \hess \eik(\mpsi(\sigma)) \mpsi'(\sigma) = \frac{\hess\eik(\mpsi(\sigma))\grad\eik(\mpsi(\sigma))}{\slow(\mpsi(\sigma))} = \grad\slow(\mpsi(\sigma)).
\end{equation}
Evaluating $\grad\eik$ along $\mpsi$, taking the derivative with respect to $\sigma$, and integrating from $\sigma_0$ to $\sigma_1$ gives:
\begin{equation}\label{eq:propagate-grad-eik}
  \grad\eik(\mpsi(\sigma_1)) = \grad\eik(\mpsi(\sigma_0)) + \int_{\sigma_0}^{\sigma_1} \grad\slow(\mpsi(\sigma))d\sigma.
\end{equation}
Once $\mpsi$ has been determined by minimizing the Fermat functional, we can then integrate along $\mpsi$ to propagate $\grad\eik$ along the ray. In this work, we limit ourselves to a constant speed of sound, in which case $\grad\eik$ is constant along ray trajectories, as can be concluded from \eqref{eq:propagate-grad-eik}.

\subsection{Curvature of the wavefront}

Continuation of the initial amplitude along ray trajectories involves determining the attenuation of the amplitude due to wavefront spreading, which in turn relates to the curvature of the wavefront. Since the wavefront at different times is an implicit surface determined by the eikonal, the curvature of the wavefront can be recovered from the Hessian of the eikonal. We elaborate this point in this section.

The $\hat\eik$-level set of the eikonal, $\Set{\mx \in \Omega : \eik(\mx) = \hat\eik}$, describes the set of points  reached by the geometric optics wavefront with eikonal value $\hat\tau$, encoded by the solution of \eqref{eq:eikonal-equation}. If $f(\mx) = 0$ describes an implicit surface (e.g., $f(\mx) = \eik(\mx) - \hat\eik$), and $\mt$ is a vector in the tangent plane at $\mx$, then the corresponding sectional curvature can be recovered from its Hessian using the formula~\cite{goldman2005curvature}:
\begin{equation}
  \kappa_{\mt}(\mx) = -\frac{\mt^\top \hess f(\mx) \mt}{\norm{\grad f(\mx)}}.
\end{equation}
That is, the sectional curvature at $\mx$ for the tangent vector $\mt$ is the curvature of a curve restricted to the wavefront surface which passes through $\mx$ with tangent vector $\mt$.

The unit surface normal is just the normalized eikonal gradient:
\begin{equation}
  \mn(\mx) = \frac{\grad\eik(\mx)}{\norm{\grad\eik(\mx)}}.
\end{equation}
Hence, the first principal direction at $\mx$ is the tangent vector $\m{q}_1 = \m{q}_1(\mx)$ in the tangent plane of the GO wavefront passing through $\mx$:
\begin{equation}
  \m{q}_1(\mx) = \Arg\max_{\substack{\norm{\m{q}} = 1 \\ \m{q}^\top \grad\eik(\mx) = 0}} \mq^\top \hess\eik(\mx) \mq,
\end{equation}
and the second principal direction is $\mq_2$ such that:
\begin{equation}
  \m{q}_i^\top \m{q}_j = \delta_{ij}, \qquad \mn(\mx)^\top \m{q}_i, \qquad i, j = 1, 2,
\end{equation}
where $\delta_{ij}$ is the Kronecker delta. The principal curvatures can be recovered from:
\begin{equation}\label{eq:principal-curvatures}
  \kappa_{\mq_i}(\mx) = \frac{\mq_i^\top \hess\eik(\mx)\mq_i}{\slow(\mx)}, \qquad i = 1, 2
\end{equation}
At the same time, the inner product with the surface normal satisfies:
\begin{equation}
  \frac{\mn(\mx)^\top \hess\eik(\mx) \mn(\mx)}{\slow(\mx)} = \frac{\grad\eik(\mx)^\top\grad\slow(\mx)}{\slow(\mx)^2}.
\end{equation}
Note that:
\begin{equation}
  \frac{\grad\eik(\mx)^\top\grad\slow(\mx)}{\slow(\mx)^2} = \frac{1}{\slow(\mx)} \left.\frac{d\slow}{d\sigma}\right|_{\mx}
\end{equation}
This leads to the following eigenvalue decomposition of $\hess\eik(\mx)/s(\mx)$:
\begin{equation}\label{eq:c-hess-tau-eig-decomp}
  \frac{\hess\eik(\mx)}{\slow(\mx)} = \frac{1}{\slow(\mx)} \left.\frac{d\slow}{d\sigma}\right|_{\mx} \mn(\mx)\mn(\mx)^\top + \sum_{i=1}^2 \kappa_{\mq_i}(\mx)\mq_i\mq_i^\top.
\end{equation}
To remember the sign of the curvature, recall that the sectional curvature gives the leading coefficient of a quadratic function approximating the wavefront in that section. The sign of the curvature depends on the orientation of wavefront. By convention, we assume that the surface normal points in the direction of ray propagation. For an expanding wavefront passing through a point, we expect the wavefront to be ``behind'' the tangent plane. Hence, a negative curvature corresponds to a spreading wavefront, and a positive curvature corresponds to a focusing wavefront.

We remark that the eikonal for a smooth speed of sound on $\Omega$ is globally only $C^0(\domain)$. All derivatives are singular at point sources. Derivatives at caustics due to edge diffraction (i.e., rarefaction fans) can be computed in the direction of ray propagation but are singular otherwise. Downwind of edge diffraction along the \emph{shadow boundary}, the eikonal is $C^1$: the gradient is continuous, although there is a jump discontinuity in the Hessian. The same comment applies to the \emph{reflection boundary} of a specularly reflected GO wave.

\subsection{Boundary conditions for the reflected eikonal}\label{ssec:eik-refl-BCs}

\begin{figure}
  \centering
  \begin{subfigure}[c]{0.5\linewidth}
    \includegraphics{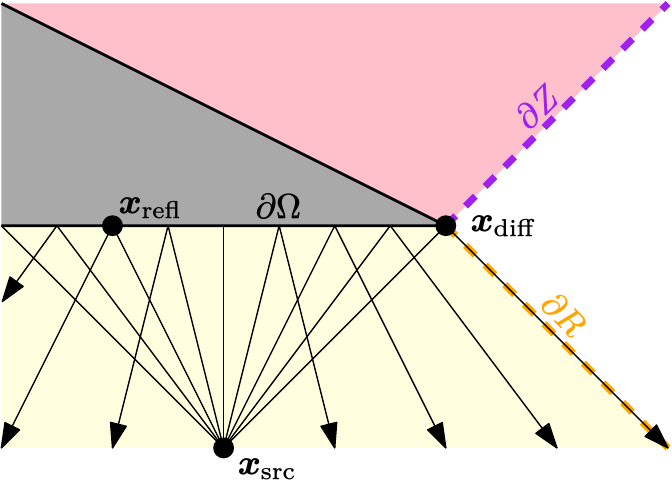}
    \caption{A field of reflected rays in 2D.}\label{fig:refl-BCs}
  \end{subfigure}%
  \begin{subfigure}[c]{0.5\linewidth}
    \includegraphics{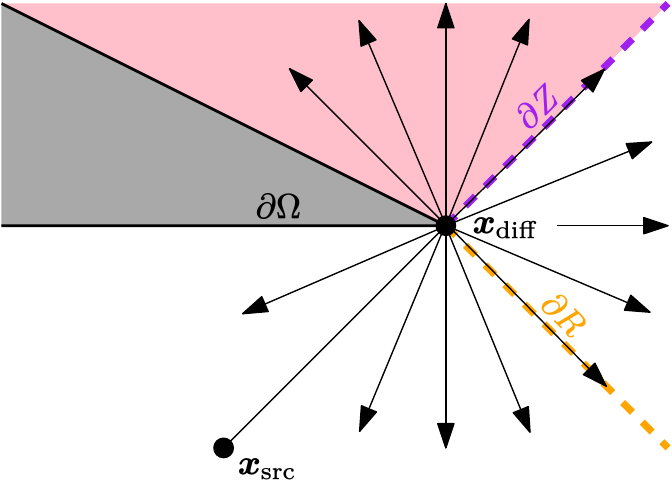}
    \caption{A field of diffracted rays in 2D.}\label{fig:diff-BCs}
  \end{subfigure}
  \caption{Simplified two-dimensional examples of scattered fields  of reflected and diffracted rays generated by an incident field due to a point source position at $\mxsrc$. The rays emitted from the subset of $\partial\Omega$ not incident on the shadow zone $Z$ are reflected specularly from the surface---i.e., the incident and reflected rays make the same angle with the tangent plane at each point of reflection. Edge-diffracted rays make the same angle with the tangent vector of the edge---in 2D, rays are emitted in all directions, while in 3D, edge-diffracted rays emitted from a single point span Keller's cone.}\label{fig:BCs}
\end{figure}

To compute the eikonal for a reflected field, we need to enforce the appropriate BCs for \eqref{eq:eikonal-equation}. These BCs are simple: we require continuity of the eikonal along reflected rays. Hence, for each point from which we would like to reflect a ray, we set $\eikin = \eikout$.

Since we assume that $\boundary\domain$ is a piecewise linear surface, we can associate a distinct reflected field to each facet of $\boundary\domain$. A point on a particular facet can be parametrized as:
\begin{equation}
  \mx_{\mlam} = \mx_0 + \mV\mlam, \qquad \mV = \begin{bmatrix} \mv_1 & \mv_2 \end{bmatrix} \in \mathbb{R}^{3 \times 2},
\end{equation}
where $\mv_1$ and $\mv_2$ are a pair of tangent vectors spanning the facet. Then, if $\mx_{\mlam}$ lies on the reflecting facet:
\begin{equation}
\eikin(\mx_{\mlam}) = \eikout(\mx_{\mlam}).\label{eq:eik-refl-BC}
\end{equation}
We also need BCs for $\grad\eik$. Taking the gradient of \eqref{eq:eik-refl-BC} with respect to $\mlam$ gives:
\begin{equation}
  \mV^\top \grad\eikin(\mx_{\mlam}) = \mV^\top \grad\eikout(\mx_{\mlam}).
\end{equation}
This gives two constraints on $\grad\eikout$. The remaining constraint is given by \eqref{eq:eikonal-equation}:
\begin{equation}
  \norm{\grad\eikin(\mx_{\mlam})} = s(\mx_{\mlam}) = \norm{\grad\eikout(\mx_{\mlam})}.
\end{equation}
Hence, $\grad\eikout(\mx_{\mlam})$ has the same magnitude and makes the same angles with the tangent plane as $\grad\eikin(\mx_{\mlam})$ does. There are two choices of $\grad\eikout(\mx_{\mlam})$ which satisfy these conditions: either the two gradients are equal (in which case transmission occurs) or they are not---the case of reflection. We can conclude from this that:
\begin{equation}
  \grad\eikout(\mx_{\mlam}) = \parens{\eye - 2 \mV \parens{\mV^\top\mV}^{-1} \mV^\top}\grad\eikin(\mx_{\mlam}) = \parens{\eye - 2\mn\mn^\top}\grad\eikin(\mx_{\mlam}),
\end{equation}
where $\mn$ is the surface normal of the reflecting facet. That is, the BCs for $\grad\eik$ are determined by specular (mirror-like) reflection from the facet. See \Cref{fig:refl-BCs}.

% Note: in order to obtain Snell's law describing the angle a ray makes as it refracts through two media with different speeds of sound, we would need to assume a model where $\slow(\mx)$ is discontinuous. In this work, our general assumption is that $\slow$ is constant (hence, smooth) in $\Omega$.

\subsection{Boundary conditions for the edge-diffracted eikonal}\label{ssec:BCs-edge-diffraction}

The eikonal is continuous along diffracted rays. Following a similar procedure as in \Cref{ssec:eik-refl-BCs}, we parametrize a point on the diffracting edge by:
\begin{equation}
  \mx_\lambda = \mx_0 + \lambda \mv,
\end{equation}
where $\mv$ spans the edge. Again, we require:
\begin{equation}
  \eikin(\mx_\lambda) = \eikout(\mx_\lambda)
\end{equation}
to hold along the edge. Differentiating once with respect to the edge parameter $\lambda$ gives:
\begin{equation}
  \mv^\top \grad\eikin(\mx_\lambda) = \mv^\top \grad\eikout(\mx_\lambda).
\end{equation}
For a 2D example, see \Cref{fig:diff-BCs}; for 3D, \Cref{fig:edge-centered-coords}.

The set of points reached by a ray originating from a fixed point $\mx_\lam$ on the edge is the surface of a cone (known as \emph{Keller's cone}~\cite{keller1962geometrical}). The set of cones reached by rays diffracted from each point on the edge together foliate the domain in the vicinity of the edge. Although $\mv^\top \grad\eikin$ is continuous, the gradient in the orthogonal complement of $\mv$ is singular. Consequently, we cannot specify BCs for $\grad\eik$ on the edge itself. However, for a point $\mxhat$ which does not lie on the edge but lies near enough that issues of visibility can be disregarded, we can compute $\tau(\mxhat)$ from:
\begin{equation}
  \eik(\mxhat) = \min_\lambda \left\{\eik(\mx_\lambda) + \int_{\sigma_0}^{\sigma_1} s(\mpsi(\sigma)) \norm{\mpsi'(\sigma)} d\sigma\right\},
\end{equation}
thereby parametrizing $\mpsi$, where $\mxhat = \mpsi(\sigma_1)$. Once we have determined $\mpsi$, we set $\grad\eikout(\mxhat) = s(\mpsi(\sigma_1))\mpsi'(\sigma_1)$.

\section{The amplitude equation}

In this section, we will explain how to continue the initial amplitude data along characteristics. Recall that we use $\sigma$ for the arc length parameter along a ray. So, along $\mpsi$, we have:
\begin{equation}
  \frac{d\alpha}{d\sigma} = \frac{\nabla\tau(\mpsi(\sigma))^\top\nabla\alpha(\mpsi(\sigma))}{\slow(\mpsi(\sigma))},
\end{equation}
allowing us to rewrite \eqref{eq:amplitude-equation} as:
\begin{equation}\label{eq:tmp1}
  \frac{d\alpha}{d\sigma} = \frac{\Delta\eik(\mpsi(\sigma))}{2 \slow(\mpsi(\sigma))}\alpha(\mpsi(\sigma)).
\end{equation}
Let $\sigma_0 < \sigma_1$ be the parameters determining the endpoints of a ray segment of interest, and assume that $\alpha(\mpsi(\sigma_0))$ is known. Then, from \eqref{eq:tmp1}:
\begin{equation}
  \alpha(\mpsi(\sigma_1)) = \alpha(\mpsi(\sigma_0)) \exp\parens{\int_{\sigma_0}^{\sigma_1} \frac{\Delta\tau(\mpsi(\sigma))}{2\slow(\mpsi(\sigma))}d\sigma}.
\end{equation}
Using \eqref{eq:c-hess-tau-eig-decomp}, we have:
\begin{equation}
  \frac{\Delta\eik(\mpsi(\sigma))}{\slow(\mpsi(\sigma))} = \Trace\left(\frac{\hess\eik(\mpsi(\sigma))}{\slow(\mpsi(\sigma))}\right) = \frac{1}{s(\mpsi(\sigma))}\frac{d(\slow\circ\mpsi)}{d\sigma} + \kappa_1(\mpsi(\sigma)) + \kappa_2(\mpsi(\sigma)),
\end{equation}
where $\slow\circ\mpsi$ is the composition of $\slow$ with $\mpsi$. Altogether, integrating from $\sigma_0$ to $\sigma_1$ gives:
\begin{equation}
  \frac{\alpha(\mpsi(\sigma_1))}{\alpha(\mpsi(\sigma_0))} = \sqrt{\frac{s(\mpsi(\sigma_1))}{s(\mpsi(\sigma_0))}} \exp\parens{\frac{1}{2}\int_{\sigma_0}^{\sigma_1}(\kappa_1(\mpsi(\sigma)) + \kappa_2(\mpsi(\sigma)))d\sigma}.\label{eq:amplitude-transport}
\end{equation}
We can check the sense of this expression as follows. If the slowness decreases from $\sigma_0$ to $\sigma_1$ so that $s(\mpsi(\sigma_1))/s(\mpsi(\sigma_0)) < 1$, the speed is increasing---hence, the field of rays is rarefying/spreading and the amplitude is decreasing. At the same time, if the GO wave is expanding, $\kappa_1$ and $\kappa_2$ are negative so that the exponential in \eqref{eq:amplitude-transport} is in the open interval $(0, 1)$.

\subsection{Boundary conditions for the reflected amplitude}

The boundary conditions for the reflected amplitude can be derived from the boundary conditions for Helmholtz equation~\cite{molinet2011acoustic}:
\begin{equation}
  \ampout = R \ampin, \qquad R = \begin{cases}
    1 & \mbox{ for sound hard reflection}, \\
    -1 & \mbox{ for sound soft reflection}.
  \end{cases}
\end{equation}
To allow for some attenuation due to material properties, $R = R(\omega)$ is usually taken to be a complex number such that $|R| \leq 1$ which depends on the frequency. A variety of damping and impedance effects can be modeled this way. Incorporating this type of empirical reflection coefficient into the algorithms presented in this paper is straightforward, so we limit ourselves to considering sound-hard reflection: that is, $R = 1$ so that $\ampout = \ampin$.

\subsection{Boundary conditions for the edge-diffracted amplitude}\label{ssec:BCs-for-edge-diffracted-amplitude}

The boundary conditions for the edge-diffracted amplitude are significantly more complicated. There are several different approaches to incorporating these BCs. We use the edge-difraction coefficient coming from the uniform theory of diffraction (UTD)~\cite{kouyoumjian1974uniform}. UTD gives a GO approximation of the edge-diffraction wave field which is asymptotically accurate for large frequencies. UTD coefficients for a wide variety of canonical scattering geometries are available~\cite{mcnamara1990introduction,molinet2011acoustic}. The edge-diffraction coefficient can also be propagated along characteristics, making it a natural fit for our solver.

UTD is a local theory of edge diffraction. It modifies the usual geometric spreading of a wavefront generated by a vibrating edge with a diffraction coefficient which depends only on the incident and diffracting rays, and local geometry at the point of diffraction. It is unclear whether a theory for a varying speed of sound has been developed (or whether it is necessary). This is one reason we limit our attention to a constant speed of sound. For the remainder of this section, we assume that $\speed = 1/\slow$ is constant.

\begin{figure}
  \centering
  \includegraphics[width=\linewidth]{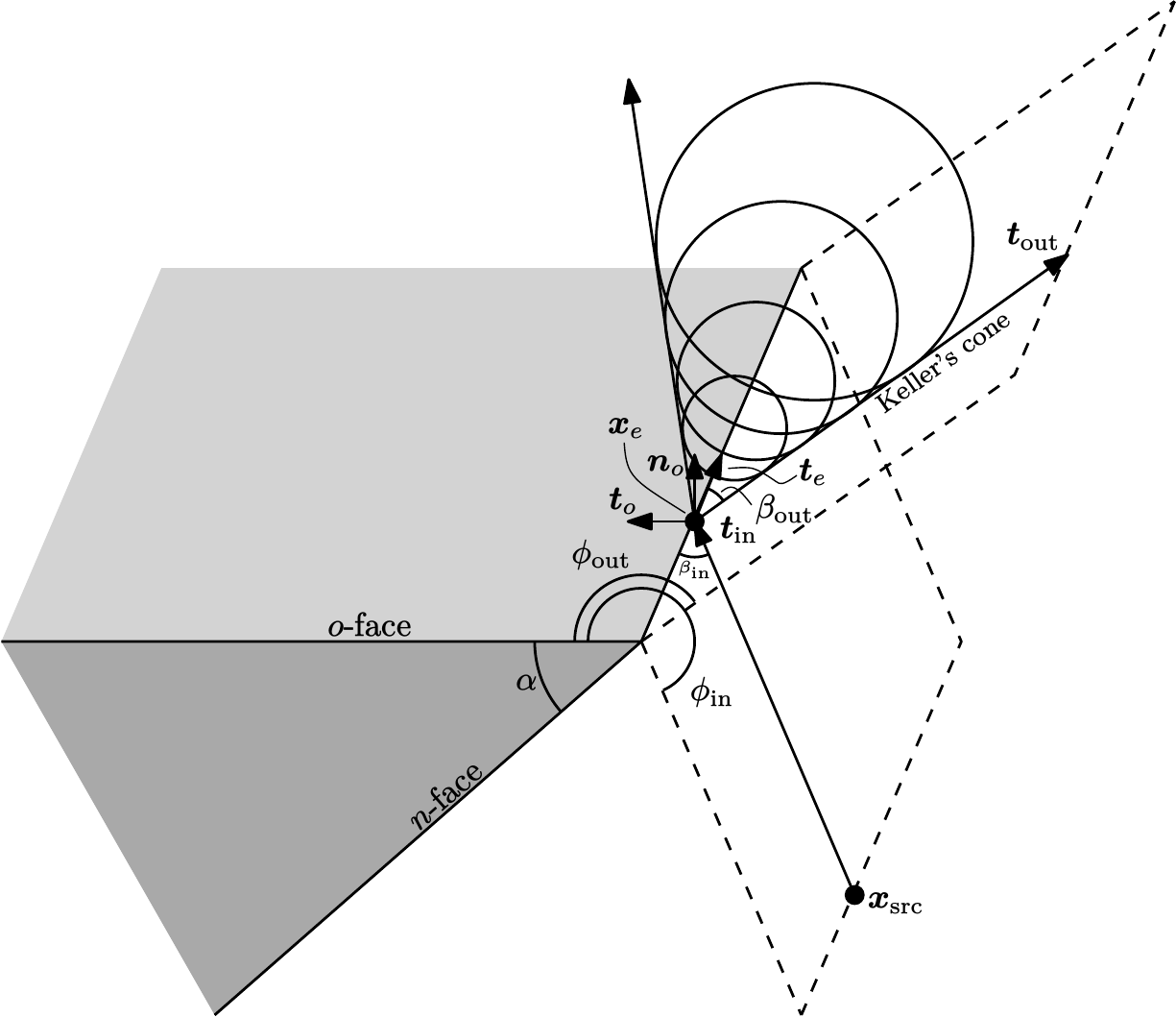}
  \caption{An overview of edge diffraction from a 3D planar wedge. This includes the edge-centered coordinate system for diffraction from a planar wedge (``edge diffraction''). The standard convention in UTD is to let the ``$o$-face'' be the facet of the wedge for which $\mt_o = \mn_o \times \mt_e$ points \emph{into} the wedge. The other face is the ``$n$-face''. The angle of the wedge is $\alpha = (2 - n)\pi$. Also shown is Keller's cone. In this case, since the edge is straight, the axis of Keller's cone is just the edge itself. A sequence of circles lying in the normal plane at points along the edge are also depicted. The edge passes through their centers. The rays diffracted from a point lie on the surface of the corresponding Keller cone, and the cones taken together foliate a subset of $\mathbb{R}^3$.}
  \label{fig:edge-centered-coords}
\end{figure}

If $\udiff$ is the edge-diffracted field satisfying \eqref{eq:helmholtz-equation} and $\uin$ is the incident field, then the edge-diffracted field takes the form:
\begin{equation}
  \udiff(\mx) = D \sqrt{\frac{\rhoe(\mx)}{\rhodiff(\mx)\big(\rhoe(\mx) + \rhodiff(\mx)\big)}} \uin(\mx),
\end{equation}
where $\mx$ is field point, $D$ is the UTD edge-diffraction coefficient, $\rhodiff$ is the distance from the point of diffraction to the field point, and $\rhoe$ is the distance between the field point and its projection onto the edge (note that the projected point may not actually lie on $\boundary\domain$). That is, if $\mx_e$ is the point of diffraction and $\mt_e$ is the unit tangent vector for the edge (see \Cref{fig:edge-centered-coords}):
\begin{equation}\label{eq:xproj}
  \mxproj = \mx_e + \mt_e \mt_e^\top (\mx - \mx_e),
\end{equation}
then $\rhoe(\mx) = \norm{\mx - \mxproj}$. The factor $\sqrt{\rhoe/(\rhodiff(\rhoe + \rhodiff))}$ is a spreading factor which reflects the cylindrical nature of the edge-diffracted field.

The coefficient $D$ takes the form:
\begin{equation}
  D = D_1 + D_2 + R \cdot (D_3 + D_4),
\end{equation}
where $R$ is the material reflection coefficient from the previous section. The diffraction coefficients $D_1$, $D_2$, $D_3$, and $D_4$ correspond to different pieces of the scattered field and are derived in~\cite{kouyoumjian1974uniform} and~\cite{molinet2011acoustic}. We present more relevant details in~\ref{sec:UTD-coef}. As before, we assume that $R \equiv 1$ for simplicity. The edge diffraction coefficient depends on the angle of the wedge, the wavenumber $k$, the direction of the ray as contacts the edge ($\mtin$), and the direction of the ray as it leaves the edge ($\mtout$). We pay particular attention to the last two parameters, $\mtin$ and $\mtout$. Observe carefully that we need to evaluate $D = D(\mx)$ at a generic field point $\mx$ which \emph{does not} lie on the edge, but must still evaluate $\mtin = \mtin(\mx)$ and $\mtout = \mtout(\mx)$---in fact, $\mtout$ is undefined on the edge. Doing so is trivial in a Lagrangian framework, since global information about the ray passing through each point is retained. However, with a semi-Lagrangian method, we ``forget'' the global information. To make up for this, we show how to use a simple dynamic programming algorithm to march $\mtin$ and $\mtout$ throughout $\Omega$ in \Cref{ssec:evaluating-UTD-coefs}.

\section{Simplification for a constant speed of sound}

In this paper, we focus on the case of a point source with a constant speed of sound. Although it is not completely general, it is nevertheless an (possibly \emph{the}) important case:
\begin{enumerate}
\item Simulating spherical waves emitted by point sources allows us to model the basic transfer functions and impulse responses needed in room acoustics.
\item It is possible to compute the analytic solution for simple geometries.
\item It is not an ``easy'' speed function from the perspective of an eikonal solver: at rarefaction fans, the solution possesses singularities which make it challenging to design a solver which gets higher than first-order convergence. This will be true of any speed function which is a constant plus a smooth perturbation.
\item A constant speed of sound is widely used in acoustic modeling. Bending of rays due to temperature variation is an important but second order effect.
\end{enumerate}
All of the algorithms developed in this paper can be generalized to handle a nonconstant speed of sound, but this is a significant enough undertaking to deserve its own treatment.

\subsection{Free-space solution for a point source}\label{ssec:free-space-point-source}

We start by considering a point source with $\Omega = \mathbb{R}^3$. We assume that:
\begin{equation}
  \tau(\mxsrc) = 0,
\end{equation}
so that the eikonal and its gradient are given by:
\begin{equation}
  \eik(\mx) = \norm{\mx - \mxsrc}/c, \qquad \grad\eik(\mx) = \frac{1}{c} \cdot \frac{\mx - \mxsrc}{\norm{\mx - \mxsrc}}.
\end{equation}
The Hessian of the eikonal is:
\begin{equation}
  \hess\eik(\mx) = \frac{1}{c} \cdot \frac{1}{\norm{\mx - \mxsrc}}\parens{\eye - \frac{\mx - \mxsrc}{\norm{\mx - \mxsrc}}\frac{\parens{\mx - \mxsrc}^\top}{\norm{\mx - \mxsrc}}}.
\end{equation}
This choice of eikonal BCs corresponds to a scaled Green's function for the Helmholtz equation. Hence, the amplitude satisfies:
\begin{equation}
  \amp(\mx) \propto \frac{1}{\norm{\mx - \mxsrc}}
\end{equation}
The constant of proportionality here corresponds to the loudness of an impulse emitted from $\mxsrc$, and can be scaled as necessary. Since the wavefronts for this problem (the level sets of the eikonal) are spherical, the principal curvatures of the wavefront passing through the point $\mx$ are:
\begin{equation}
  \kappa_1(\mx) = \kappa_2(\mx) = -\frac{1}{\norm{\mx - \mxsrc}}.
\end{equation}
To continue the amplitude along a ray $\mpsi$, we integrate (recalling that $\sigma$ is the arc length parameter):
\begin{equation}
  \int_{\sigma_0}^{\sigma_1} \frac{\kappa_1(\mpsi(\sigma)) + \kappa_2(\mpsi(\sigma))}{2} d\sigma = -\int_{\sigma_0}^{\sigma_1} \frac{d\sigma}{\sigma} = \log \frac{\sigma_0}{\sigma_1},
\end{equation}
from which we conclude the usual ``$1/r$'' scaling of amplitude (or the inverse squared law for the sound intensity):
\begin{equation}
  \frac{\alpha(\mpsi(\sigma_1))}{\alpha(\mpsi(\sigma_0))} = e^{\log(\sigma_0/\sigma_1)} = \frac{\sigma_0}{\sigma_1}.
\end{equation}
The foregoing considerations give us a way of setting the BCs for the eikonal, its gradient, and the amplitude near a point source.

\section{Numerical discretizations}

Next, we discuss the discretizations made when approximating the continuous scattering problem. In particular, we approximate the domain $\domain$ with a Delaunay tetrahedron mesh, and represent various quantities over the edges, faces, and cells of this mesh using Hermite interpolation. Finally, we discuss how we approach solving the semi-Lagrangian updates using numerical optimization.

\subsection{Spatial discretization: conforming tetrahedron mesh}

We assume that $\boundary\domain$ is piecewise linear. Specifically, we assume that $\boundary\domain$ is comprised of $O(1)$ flat, polygonal facets. The facets and $\Omega$ itself are allowed to be nonconvex. This class of domains is broad and gives rise to a wide variety of important scattering phenomenon which must be modeled accurately to simulate acoustic wave propagation.

Since $\boundary\domain$ is piecewise linear, it is amenable to tetrahedralization~\cite{cheng2013delaunay}. We use the package TetGen~\cite{hang2015tetgen} to discretize $\Omega$ into an exact, constrained Delaunay tetrahedralization. Note that it is not necessary to generate a triangulation of $\boundary\domain$ first. It is sufficient to provide a description of each of the polygonal facets defining $\boundary\domain$.

When we discretize $\Omega$, we specify a maximum tetrahedron volume constraint. We use this parameter to generate a sequence of meshes with different finenesses for our numerical tests. For each tetrahedron mesh, we define the parameter $h$ to be the average edge length of the mesh. We denote the tetrahedron mesh by:
\begin{equation}
  \Omega_h = (\calV_h, \calE_h, \calF_h, \calC_h).
\end{equation}
We denote by $\calV_h$, $\calE_h$, $\calF_h$, and $\calC_h$ the sets of vertices, edges, faces (triangles), and cells (tetrahedra) respectively comprising the tetrahedron mesh. We assume that the mesh is regular so that the intersection of two cells $C, C' \in \calC_h$ is either 1) the empty set, 2) exactly one vertex in $\calV_h$, 3) exactly one edge in $\calE_h$, or 4) exactly one face in $\calF_h$. We also assume that the mesh consists of a single connected component.

It is helpful to think of the graph whose vertex set is $\calV_h$ and whose edge set is $\calE_h$. When we use graph terminology to refer to $\Omega_h$, this is the graph we have in mind. Although we do not dwell overly on computational complexity in this work, we assume roughly that the degree of each vertex in $\calV_h$ is $O(1)$. We also assume that the length of each edge is $O(h)$, that the area of each face is $O(h^2)$, and that the volume of each cell is $O(h^3)$. That is, we assume that the tetrahedron mesh is reasonably uniform.

\subsection{Numerical discretization: Bernstein-\Bezier{} interpolation on triangles}\label{ssec:bernstein-bezier}

\begin{figure}
  \centering
  \hspace{0.025\linewidth}%
  \includegraphics[width=0.975\linewidth]{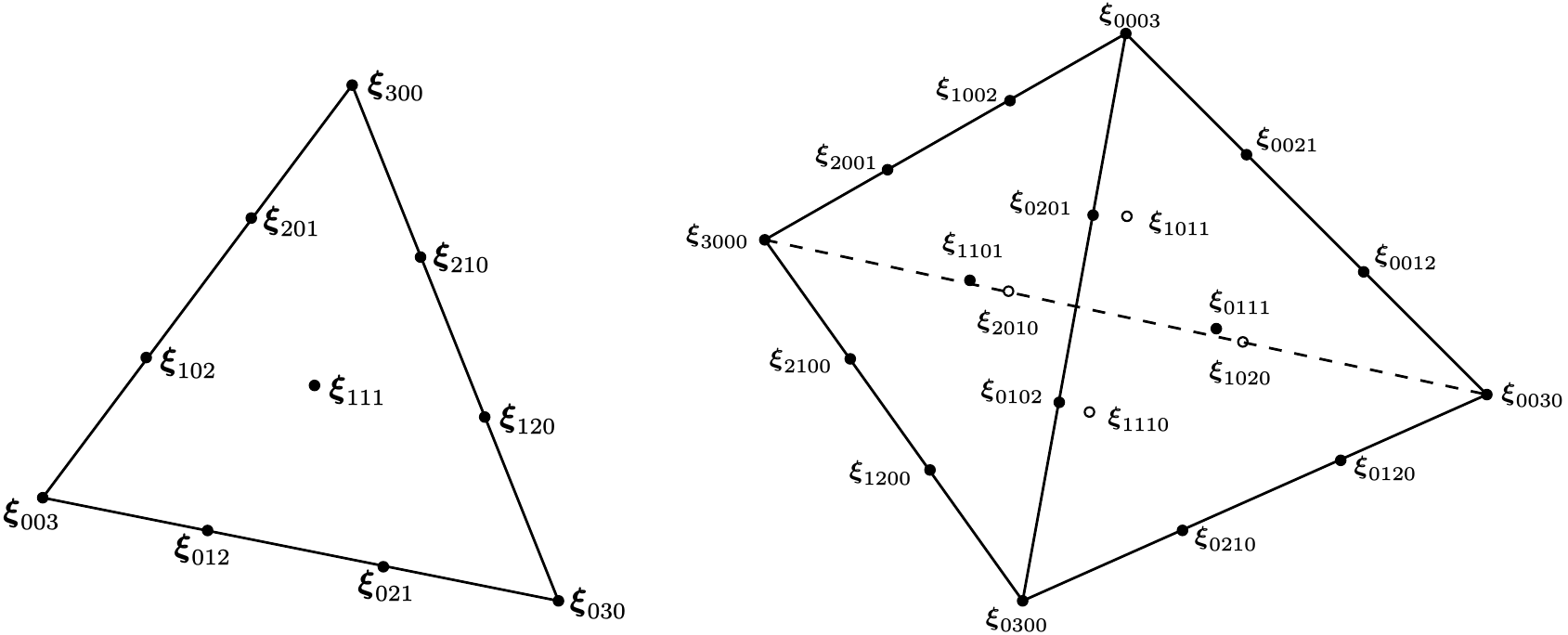}
  \caption{\emph{Left}: the control points for the degree 3 bivariate Bernstein-\Bezier{} triangle. For the ``9-parameter element'', the \Bezier{} coefficient corresponding to $\mx_{111}$ is computed using condensation of parameters. \emph{Right}: the control points for the degree 3 trivariate Bernstein-\Bezier{} tetrahedron. Since all control points lie on the faces of the tetrahedral cell, computing the data for the 9-parameter interpolant for each of the four triangular faces is enough to fully determine the data for the ``20-parameter element''.}\label{fig:9-param-and-20-param}
\end{figure}

Consider a triangle with vertices $\mx_1, \mx_2, \mx_3$, let $\mx_{\m{\gamma}} = \gamma_1\mx_1 + \gamma_2\mx_2 + \gamma_3\mx_3$, where $\m{\gamma} = (\gamma_1, \gamma_2, \gamma_3)$ are the barycentric coordinates of $\mx_{\m{\gamma}}$ such tha  $\gamma_i \geq 0$ ($i = 1, 2, 3$), and $\gamma_1 + \gamma_2 + \gamma_3 = 1$. The degree $n$ Bernstein polynomials $B_{ijk}^n$ span the space of bivariate polynomials of degree $n$, and are given by:
\begin{equation}
  B_{ijk}^n(\gamma_1,\gamma_2,\gamma_3) = \frac{n!}{i!j!k!} \gamma_1^{i}\gamma_2^{j}\gamma_3^{k}, \quad i \geq 0, \quad j \geq 0, \quad k \geq 0 \quad i + j + k = n.
\end{equation}
A polynomial $b$ of degree $n$ is written in the Bernstein basis as:
\begin{equation}
  b(\gamma_1, \gamma_2, \gamma_3) = \sum_{\substack{i + j + k = n \\ i, j, k \geq 0}} b_{ijk} B_{ijk}^n(\mgamma).
\end{equation}
The coefficients $b_{ijk}$ are referred to as \Bezier{} coefficients, and they are associated with a set of control points. For interpolation over the triangle, the grid of control points is:
\begin{equation}
  \m{\xi}_{ijk} = \frac{i}{n} \mx_1 + \frac{j}{n} \mx_2 + \frac{k}{n} \mx_3, \qquad i \geq 0, \quad j \geq 0, \quad k \geq 0, \quad i + j + k = n.
\end{equation}
See \Cref{fig:9-param-and-20-param}. A variety of Bernstein-\Bezier{} ``elements'' exist which can be used for Hermite interpolation over the simplex~\cite{lai2007spline}. Using the \Bezier{} form is stable and allows us to use efficient pyramid algorithms for interpolation, evaluation, and shape interrogation~\cite{farin1986triangular}.

\paragraph{The 9-parameter Hermite interpolant} We use Hermite interpolation to approximate the eikonal over the convex hull of triples of points $\mx_0, \mx_1, \mx_2 \in \calV_h$ to do semi-Lagrangian updates. We discuss the semi-Lagrangian updates more in \Cref{ssec:tetrahedron-updates}, but collect the details of the Hermite element that we use here.

The space of bivariate polynomials of degree $n = 3$ has dimension ${n + 2 \choose 2} = 10$. We use the ``9-parameter'' Hermite interpolant~\cite{farin1986triangular}. See \Cref{fig:9-param-and-20-param}. To interpolate the function values $f(\mx_0), f(\mx_1)$, and $f(\mx_2)$ and gradients $\grad f(\mx_0), \grad f(\mx_1)$, and $\grad f(\mx_2)$, we set the ordinates corresponding to the triangle vertices so that:
\begin{equation}
  b_{300} = f(\mx_0), \qquad b_{030} = f(\mx_1), \qquad b_{003} = f(\mx_2).
\end{equation}
For each control point adjacent to a vertex (six total), we set, e.g.:
\begin{equation}
  b_{210} = b_{300} + \frac{1}{3} \parens{\mx_1 - \mx_0}^\top \nabla f(\mx_0),
\end{equation}
with $b_{201}, b_{120}, b_{021}, b_{102}$, and $b_{012}$ set analogously. This leaves the remaining interior ordinate:
\begin{equation}
  b_{111} = \frac{1}{4} \big(b_{210} + b_{201} + b_{120} + b_{021} + b_{102} + b_{012}\big) - \frac{1}{6} \big(b_{300} + b_{030} + b_{003}\big).
\end{equation}
This gives quadratic precision~\cite{farin1986triangular}. In general, degree $n$ Lagrange interpolation over a triangle using the grid of control points $\{\mx_{ijk}\}_{i+j+k=n}$ is at most $O(h^{n-p})$ accurate, where $n$ is the degree of the space of bivariate polynomials and $p$ is the order of the derivative being approximated~\cite{lai2007spline}. Hence, we can expect $O(h^{3-p})$ accuracy, which is sufficient for our purposes---we find that this is born out numerically.

After running the jet marching method, $\Eik$ and $\grad\Eik$ are known at each vertex $\mx \in \calV_h$. Using this data, we can construct a piecewise Hermite interpolant (tetrahedral spline) over $\Omega_h$. We use the 20-parameter degree 3 Bernstein-B\'{e}zier polynomial on the tetrahedron. This is a straightforward extension of the 9-parameter triangular interpolant to 3D. See \Cref{fig:9-param-and-20-param}. Since the degree is less than the number of vertices, there are no interior nodes for which coefficients must be computed. In fact, if we compute the coefficients for the 9-parameter interpolants for each face of the tetrahedron, we will have ended up computing all coefficients of the 20-parameter interpolant. Note that the coefficients for the spline over $\Omega_h$ can be stored efficiently if we associate the vertex coefficients to the vertices in $\calV_h$, the pairs of edge coefficients to edges in $\calE_h$, and the face coefficients to the faces in $\calF_h$. Altogether, the number of floating point numbers required to store the spline coefficients is $|\calV_h| + 2 |\calE_h| + |\calF_h|$.

\subsection{Tetrahedral splines and extracting level sets}\label{ssec:extracting-level-sets}

Since the level sets of the eikonal correspond to the propagating geometric optics wavefront at different points in time, it is helpful to have an algorithm to extract them for the purposes of visualization. Once we have computed the piecewise Bernstein-B\'{e}zier polynomial interpolating the numerical jets of the eikonal, we have a data structure which encodes a continuum of geometric optic wavefronts. To display a particular wavefront, we must extract it from this data structure and render it. One approach to doing so is raytracing.

A Bernstein-B\'{e}zier patch is contained in the convex hull of its control points~\cite{farin1986triangular}. A consequence of this is that the minimum (resp., maximum) value of the Bernstein-B\'{e}zier polynomial over the set of nonnegative barycentric coefficients is bounded below (resp., above) by the minimum (resp., maximum) values of the control points. This lets us quickly determine which Bernstein-B\'{e}zier elements do not intersect a particular wavefront. Eliminating these polynomials leaves a set of tetrahedra which ``bracket'' the position of the wavefront. We can then test individually whether each ray emitted from a camera intersects a bracket tetrahedra, and then determine the exact position where a ray intersects the wavefront. This procedure can be further accelerated by spatially sorting the bracket tetrahedra using an R-tree~\cite{guttman1984r}.

\subsection{Triangle updates}\label{ssec:triangle-updates}

We model diffraction and two-dimensional ray propagation restricted to $\partial\Omega$ using two-point \emph{triangle updates}. These are local raytracing problems where $\tau$ and $\grad\tau$ are known at two \emph{update points} $\mx_0$ and $\mx_1$, and where we would like to compute $\tau$ and $\grad\tau$ at another point $\mxhat$. We parametrize a point in the interval $[\mx_0, \mx_1]$ by:
\begin{equation}
  \mx_{\lambda} = \mx_0 + \lambda(\mx_1 - \mx_0), \qquad 0 \leq \lambda \leq 1,
\end{equation}
and approximate the eikonal using a cubic Hermite polynomial over $[\mx_0, \mx_1]$. If $H_i^j(\lambda)$ are the cardinal basis functions for Hermite interpolation on $[0, 1]$, where $i = 0, 1$ denotes the endpoint and $j = 0, 1$ denotes the derivative~\cite{prenter2008splines}, we approximate $\tau$ with:
\begin{equation}
  T(\lambda) = T(\mx_0) H_0^0(\lambda) + T(\mx_1) H_0^1(\lambda) + (\mx_1 - \mx_0)^\top \grad T(\mx_0) H_0^1(\lambda) + (\mx_1 - \mx_0)^\top \grad T(\mx_1) H_1^1 (\lambda).
\end{equation}
Then, the triangle update corresponds to the nonlinear constrained optimization problem:
\begin{equation}
  \begin{split}
    \mbox{minimize} \qquad & T(\lambda) + \frac{1}{c} \norm{\mxhat - \mx_{\lam}} \\
    \mbox{subject to} \qquad & 0 \leq \lambda \leq 1.
  \end{split}
\end{equation}
This problem can be easily solved using a hybrid rootfinder~\cite{stewart1998afternotes}. For a varying speed of sound, the main difference is that the action function in Fermat's principle must be approximated using a suitable quadrature rule~\cite{potter2021jet}.

% We note the following special case which occurs with boundary triangle updates. There can be a degradation in accuracy when updating a node reached by a diffracting ray with a triangle update which has exactly one of $\mx_0$ and $\mx_1$ incident on the diffracting edge. It is clear that $\grad\eik(\mx_0)$ and $\grad\eik(\mx_1)$ may point in incompatible directions. For example, if $\mx_0$ is incident on the diffracting edge, $\grad\eik(\mx_0)$ will point in the direction of the ray incident on $\mx_0$, but $\grad\eik(\mx_1)$ will point in the diffracted direction. In this case, we manually rotate $\grad\eik(\mx_0)$ to point in the diffracted direction. If $\mt_e$ is the unit tangent vector for the diffracting edge, and $\mt_f$ is the unit tangent vector for the face containing $\mx_1$ (and pointing towards $\mx_1$), then we set:
% \begin{equation}
%   \grad\eik(\mx_0) = \cos(\beta) \mt_e + \sin(\beta) \mt_f,
% \end{equation}
% where $\cos\beta = \mt_e \cdot \grad\eik(\mx_0)$.

\subsection{Tetrahedron updates}\label{ssec:tetrahedron-updates}

Similar to triangle updates, we have \emph{tetrahedron updates}. These are used to model free-space ray propagation and involve three points where the eikonal jet is known: $\mx_0$, $\mx_1$, and $\mx_2$. We let $\mlam = (\lam_1, \lam_2)$ such that $\lam_1 \geq 0, \lam_2 \geq 0$, and $1 - \lam_1  - \lam_2 \geq 0$ parametrize the standard triangle in $\RR^2$. Then, the convex hull of $\mx_0, \mx_1$, and $\mx_2$ is parametrized by:
\begin{equation}
  \mx_{\mlam} = \mx_0 + \lam_1 (\mx_1 - \mx_0) + \lam_2 (\mx_2 - \mx_0).
\end{equation}
If we let $T(\mlam)$ be the 9-parameter element (described in \Cref{ssec:bernstein-bezier}) interpolating $\eik$ and $\grad\eik$ at $\mx_0$, $\mx_1$, and $\mx_2$, then the solving the local raytracing problem corresponds to solving the following optimization problem:
\begin{equation}
  \begin{split}
    \mbox{minimize} \qquad & T(\mlam) + \frac{1}{c} \norm{\mxhat - \mx_{\mlam}}, \\
    \mbox{subject to} \qquad & \lam_1 \geq 0, \quad \lam_2 \geq 0, \quad 1 - \lam_1 - \lam_2 \geq 0.
  \end{split}
\end{equation}
Solving this optimization problem is more involved than the 1D minimization problem for the triangle update. There are a wide variety of approaches that could be taken. We chose to use sequential quadratic programming (SQP)---this is essentially a projected Newton's method~\cite{nocedal1999numerical}. If $f(\mlam)$ denotes the cost function, at each step we form the Taylor approximation:
\begin{equation}
  f(\mlam + \delta\mlam) = f(\mlam) + \nabla f(\mlam)^\top \delta\mlam + \frac{1}{2} \delta\mlam^\top \nabla^2 f(\mlam) \delta\mlam + O(\norm{\delta\mlam}^2).
\end{equation}
At each iteration, we use the active set method to solve:
\begin{equation}
  \begin{split}
    \mbox{minimize} \qquad & f(\mlam) + \nabla f(\mlam)^\top \delta\mlam + \frac{1}{2} \delta\mlam^\top \nabla^2 f(\mlam) \delta\mlam \\
    \mbox{subject to} \qquad & \mlam + \delta\mlam \mbox{ being in the standard triangle.}
  \end{split}
\end{equation}
The next iterate is then $\mlam + \delta\mlam^*$, where $\delta\mlam^*$ is optimum. Doing this repeatedly ($f$ varies from iteration to iteration) comprises SQP. Note that the minimization problem for each iteration is a quadratic program with inequality constraints, which can be solved using the active set method~\cite{nocedal1999numerical}. Since the dimension of these optimization problems is very small, this is efficient.

\subsection{Tolerances}\label{ssec:tolerances}

When we do the updates described in the preceding sections, we must choose a tolerance to use to terminate each optimization algorithm. To do so in a manner which respects the uneven size of the elements in the mesh, we set a tolerance for $\delta\mlam$, which is scale independent. Roughly speaking, we need $\mx^*$ to approximate the optimum of the ``true'' minimization problem with $O(h^3)$ error. To this end, for a tetrahedron update, if we let $\max\ell$ be the minimum edge length of the edges comprising the update triangle $\{\mx_0, \mx_1, \mx_2\}$, and let $\operatorname{diam}(\Omega)$ denote the diameter of the domain, then we set the tolerance for the tetrahedron update to be $(\max\ell/\operatorname{diam}(\Omega))^2$.

\section{The jet marching method (JMM)}

In a previous work, we developed a 2D JMM for solving the eikonal PDE with a smoothly varying speed of sound on a regular grid~\cite{potter2021jet}. In this section, we explain how we adapt these ideas to computing the eikonal for point source BCs on an unstructured domain in 3D which has been discretized into a uniform tetrahedron mesh, but for a constant speed of sound. We also discuss how to reinitialize the eikonal to propagate reflected and edge-diffracted GO waves. Special attention is paid to edge-diffraction---a caustic forms at the site of edge diffraction, introducing a singularity which limits the overall accuracy of $T$ at the nodes downwind of the edge to $O(h \log \tfrac{1}{h})$.

\subsection{Marching methods for solving the eikonal equation}

We label the points in $\calV_h$ (the vertices of the tetrahedron mesh) with one of three states: $\mathtt{far}$, $\mathtt{trial}$, and $\mathtt{valid}$. We let $T : \calV_h \to \mathbb{R}$ denote the numerical eikonal, and let $\grad\Eik : \calV_h\to\RR$ be the numerical gradient. Nodes with BCs start with the state $\mathtt{trial}$ and have $T$ and $\grad T$ set accordingly, with the state of all other nodes set to $\mathtt{far}$ and $T$ set to $+\infty$. We sort all $\mathtt{trial}$ nodes in increasing order of $T$ into a priority queue, which is typically implemented as an array-based binary heap~\cite{sedgewick2001algorithms}.

At each step, we pop the first node from the priority queue. Label this node $\mx_0$. We set the state of $\mx_0$ to $\mathtt{valid}$, set the state of each $\mathtt{far}$ neighbor to $\mathtt{trial}$, insert those nodes into the priority queue, and then update all of $\mx_0$'s $\mathtt{trial}$ neighbors. There are two key details which must be specified: 1) how to determine the neighbors of $\mx_0$ which need to be updated, and 2) how to update them. The first point is straightforward. We call a node $\mxhat$ a neighbor of $\mx_0$ if the edge $\{\mxhat, \mx_0\}$ is in $\calE_h$. That is, we use the graph topology of the tetrahedron mesh. Updating $\mxhat$ is more delicate. We discuss how to do this in the next several sections.

\subsection{Tetrahedron update fans on the $\mathtt{valid}$ front}

At each step of a marching method, we can make the following observations:
\begin{itemize}
\item The set of $\mathtt{valid}$ nodes forms a connected component of the tetrahedron mesh.
\item The $\mathtt{far}$ nodes are only adjacent to $\mathtt{trial}$ nodes, and $\mathtt{valid}$ nodes will only be adjacent to $\mathtt{trial}$ nodes.
\item Hence, a layer of $\mathtt{trial}$ nodes separates the $\mathtt{valid}$ nodes from the $\mathtt{far}$ nodes.
\end{itemize}
We call the subset of $\mathtt{valid}$ nodes which is adjacent to $\mathtt{trial}$ nodes the \emph{$\mathtt{valid}$ front}.

\begin{figure}
  \centering
  \includegraphics[width=0.6\linewidth]{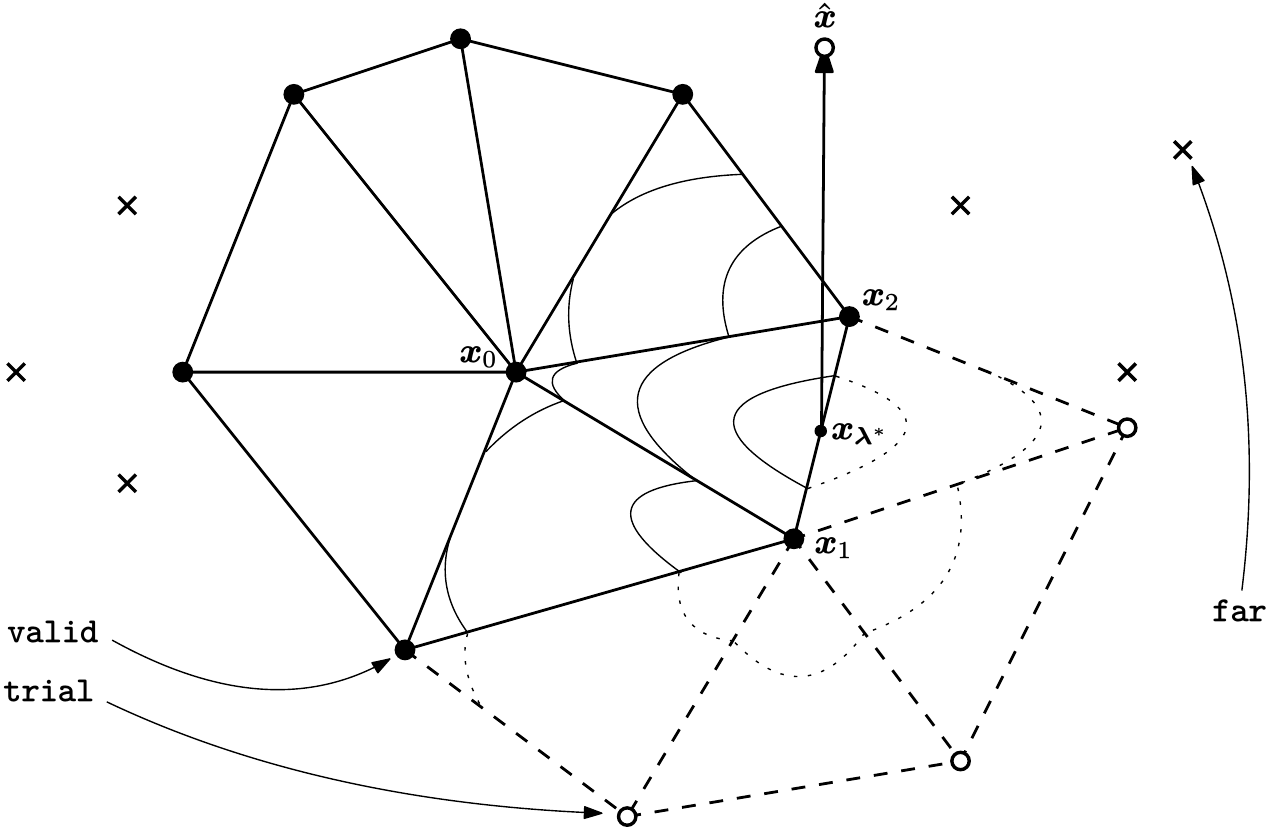}
  \caption{A depiction of the piecewise optimization problem modeling freespace ray propagation (the ``update fan'') solved at each step to update $\mxhat$. We numerically minimize Fermat's integral over each triangular facet on the $\mathtt{valid}$ front to find the minimum travel time ray. Since the travel time is in effect a $C^0$ function defined on the $\mathtt{valid}$ front, we must take care with optima which lie on common edges and vertices. In the diagram, $\mathtt{valid}$, $\mathtt{trial}$, and $\mathtt{far}$ nodes are shown---note that we only depict several triangles on the $\mathtt{valid}$ front, hence the $\mathtt{far}$ nodes appear to be detached from mesh. It is to be understood that there is a background tetrahedron mesh which is not shown.}\label{fig:update-fan}
\end{figure}

Consider the $\mathtt{valid}$ front for the solver at some intermediate stage. Let $\mx_0$ be the most recent $\mathtt{valid}$ node, and consider a neighboring $\mathtt{trial}$ node, denoted $\mxhat$ (we generally denote a node which is to be updated by $\mxhat$). Fermat's principle indicates that $\eik(\mxhat)$ can be computed exactly by applying \eqref{eq:fermat} for the original BCs of the problem. At the same time, since the $\mathtt{valid}$ nodes are now fixed for all time, we can consider them as providing BCs for a new eikonal problem with the same solution as the original, so we need only consider the BCs of this new problem restricted to the $\mathtt{valid}$ front. Hence, at each stage, we can allow the domain of the Fermat minimization problem given by \eqref{eq:fermat} to range over the $\mathtt{valid}$ front.

If we take this approach on a regular grid, then it is always sufficient to solve \eqref{eq:fermat} locally over the part of the $\mathtt{valid}$ front which is immediately adjacent to $\mx_0$. For more general static Hamilton-Jacobi equations, this is not the case, and larger stencils may be required~\cite{sethian2001ordered,sethian2003ordered,mirebeau2014anisotropic,mirebeau2014efficient}. On the other hand, for a tetrahedron mesh, the highly irregular structure of the mesh makes determining which updates to do at each step unintuitive. In the past, methods based on \emph{splitting sections} were used to ensure the \emph{causality} of each update: namely, that $T(\mxhat) > \max_i T(\mx_i)$, where the maximum is taken over each update node~\cite{kimmel1998computing,sethian2000fast}. We now describe an alternative approach.

Observe that the $\mathtt{valid}$ front can be triangulated by first considering the subset of $\calC_h$ consisting of tetrahedra whose vertices are all in a $\mathtt{valid}$ state. Contained within this subset is another subset whose members are each adjacent to a cell with three $\mathtt{valid}$ vertices and one $\mathtt{trial}$ vertex. The triangle face shared by these two cells is on the $\mathtt{valid}$ front, with the set of all such triangles giving a triangulation of the $\mathtt{valid}$ front. When we accept a node $\mx_0$, the triangulation of the $\mathtt{valid}$ front will change. There will be a set of triangles incident on $\mx_0$ which will become part of the triangulation, while a nearby set of triangles will be removed. This set of newly $\mathtt{valid}$ triangles is a \emph{triangle fan} incident on $\mx_0$ which is simple to locate in the mesh in $O(1)$ operations. If we consider the tetrahedron updates corresponding to each of these triangles, we can see that none of these updates will have been done before, since the node $\mx_0$ will have only just become $\mathtt{valid}$. So, at each step, we pull a new set of tetrahedron updates from this \emph{update fan} and execute each. See \Cref{fig:update-fan}.

Just doing each update from the fan blindly is not enough to ensure consistency. If one were to implement a solver which takes this approach, it would quickly be seen that because of the irregularity of a tetrahedron mesh, there will be updates with boundary minimizers which result in $O(1)$ error in the update ray direction, limiting the update's accuracy to $O(h)$. Frequently, such a solver will fail to compute a consistent solution, since this error will propagate downwind to $O(N)$ other nodes. In the next sections, we address this problem---and several others---which need to be handled in order to ensure consistency.

\subsection{Physical updates}

In \Cref{ssec:tetrahedron-updates}, we described how we propagate a ray locally by solving a constrained, nonlinear minimization problem. Upon solving this optimization problem, we parametrize a point $\mx_{\mlam^*} \in \operatorname{conv}(\mx_0, \mx_1, \mx_2)$. This is the start of the update ray. We can see that the ray arriving at $\mxhat$ travels first to $\mx_{\mlam^*}$, and then proceeds from $\mx_{\mlam^*}$ to $\mxhat$. Provided that the speed of sound is smooth, we expect the ray trajectory to be smooth as well. This criterion suggests that we should reject optima $\mx_{\mlam^*}$ which lie on the boundary of $\operatorname{conv}(\mx_0, \mx_1, \mx_2)$, since otherwise we would allow the ray to artifically diffract in freespace, polluting the error of $T$ for all nodes downwind of $\mxhat$. This can be done by checking the Lagrange multipliers corresponding to the active inequalities---if they are nonzero, then the ray is bending.

Because of the irregularity of the tetrahedron mesh, it is also important for us to check the orientation of each update. Note that $\grad\Eik(\mx_0)$, $\grad\Eik(\mx_1)$, and $\grad\Eik(\mx_2)$ together determine an orientation for the tetrahedron update: we can pick a surface normal for the base of the update such that it points in the same direction as each gradient vector. We do not do tetrahedron updates where $\mxhat$ lies on the ``wrong'' side of the base of the update.

\subsection{Caching updates}

The previous section gives us a criterion for accepting an individual tetrahedron update, but it turns out this is insufficient on its own. We must consider multiple adjacent updates together in order to determine when to accept boundary optima. Consider an optimum $\mx_{\mlam^*}$ lying on the boundary of multiple adjacent updates and assume that the corresponding ray is physical. If we reject boundary updates entirely, we might reject a local ray which is needed to propagate the solution in a particular direction. In order to accept this update, we consider each update producing $\mx_{\mlam^*}$ together. The reason for this is that the bases of the incident updates might not be coplanar---hence, $\mx_{\mlam^*}$ is optimum but none of the Lagrange multipliers for the updates are zero. See \ref{sec:common-minima}. In this section, we describe our approach to this problem.

The essential challenge is as follows: since we only search over the fan of updates immediately surrounding a newly \texttt{valid} node $\mx_0$, if an update $(\mxhat, \mx_0, \mx_1, \mx_2)$ has a constrained minimum in the interior of $[\mx_1, \mx_2]$, then we must remember that update in order to match it up with another incident update on the other side of $[\mx_1, \mx_2]$ which we might encounter later. To that end, consider an update fan incident on a newly $\mathtt{valid}$ $\mx_0$. For each
$\texttt{trial}$ node $\mxhat$ neighboring $\mx_0$, and for each fan
triangle $(\mx_0, \mx_1, \mx_2)$, we do the tetrahedron update
$(\mxhat, \mx_0, \mx_1, \mx_2)$. Of course, if the corresponding
Lagrange multipliers are all zero, then we've found a physical ray and
can update the eikonal jet at $\mxhat$ accordingly. But now consider
the case where at least one Lagrange multiplier is nonzero. We have argued that these tetrahedron updates may still
produce a physical ray, but we lack the information in the moment to decide
whether this is the case---we may need to wait until a future update is discovered which matches, allowing us to conclude that the ray is indeed physical. To this end, we employ an \emph{update
  cache}.

The update cache is simply a dictionary whose keys are \texttt{trial}
nodes, and whose values are tetrahedron updates that we have done in
the past and which have nonzero Lagrange multipliers. Specifically, if
$\mxhat$ is a \texttt{trial} node, then:
\begin{equation}
  \begin{split}
    \texttt{update\_cache}[\mxhat] &= \big\{(\mx_0, \mx_1, \mx_2) : (\mxhat, \mx_0, \mx_1, \mx_2) \mbox{ is an old tetrahedron update} \\
    & \qquad\qquad\mbox{with at least one nonzero Lagrange multiplier}\big\}.
  \end{split}
\end{equation}
To populate the update cache, we proceed as follows. Let
$(\mxhat, \mx_0, \mx_1, \mx_2)$ be a tetrahedron update which we have
just done, and assume that it has at least one nonzero Lagrange
multiplier. We first search $\texttt{update\_cache}[\mxhat]$ for any
tetrahedron updates which are incident on
$(\mxhat, \mx_0, \mx_1, \mx_2)$---i.e., which have a vertex or edge in common. If we
find another tetrahedron update which is edge adjacent, we check whether both tetrahedron updates have a minimizer on the common edge which is bracketed by the opposing vertices: if they do, we update the eikonal jet at $\mxhat$ using the ray (which will be the same regardless of which tetrahedron update is used to compute it) and delete the tetrahedron update we found from $\texttt{update\_cache}[\mxhat]$.

The situation for vertex adjacent updates is similar, but with one caveat. Let $\mx_i$ be the optimum vertex. Before using this ray to update the eikonal jet at $\mxhat$, we must ensure that we have found a collection of tetrahedron updates which form a fan around $\mx_i$. These updates can be thought of as bracketing $\mx_i$.

Finally, to keep the size of \texttt{update\_cache} in check, when we
accept a node $\mx_0$ and mark it \texttt{valid}, we delete
$\texttt{update\_cache}[\mxhat]$. As the solver runs,
$|\cup_{\mxhat} \texttt{update\_cache}[\mxhat]|$ will always be at most the number of
\texttt{trial} nodes. It can be fewer. Plenty of data structures exist
which implement this type of dictionary so that each update operation takes $O(\log |\cup_{\mxhat}\texttt{update\_cache}[\mxhat]|)$ time. Finding all of the tetrahedron
updates that are incident on the newly completed update then takes
$O(1)$ time by our assumption that each node has $O(1)$ neighbors
(and, hence, will only ever be updated by $O(1)$ different tetrahedron
updates during a solve). Because of this, it is clear that the extra time required to manage the update cache during a solve is $O(|\calV_h| \log |\calV_h|)$. Similar arguments can be made to establish that the extra space required is also $O(|\calV_h| \log |\calV_h|)$.

\subsection{Triangle updates: shed rays and diffracted rays}

The foregoing discussion for tetrahedron updates applies also to the triangle updates used to model rays that are shed into the shadow zone and edge-diffracted rays. We briefly describe how these are handled.

\begin{itemize}
\item \textbf{Shed rays}. If we consider the restriction of $T$ to $\partial\Omega_h$, we can see that the corresponding restriction of the valid front to $\partial\Omega_h$ is a polygonal curve imbedded in $\partial\Omega_h$. All of the preceding discussion goes through with the essential details unchanged or simplified (e.g., there can be at most two boundary triangle updates incident on a newly valid node $\mx_0$).
\item \textbf{Edge-diffracted rays}. We handle these triangle updates similarly, but skip any tetrahedron updates for which $\mx_0$ is incident on a diffracting edge. This allows us to give precedence to the triangle updates modeling edge diffraction, ``snapping'' the ray to an edge when diffraction has occurred.
\end{itemize}

\noindent We maintain a separate update cache for each type of update: free space tetrahedron updates, boundary triangle updates (for shed rays), and edge triangle updates (for edge-diffraction).

\subsection{Local visibility: cones of feasible ray directions}

% \begin{figure}
%   \centering
%   \includegraphics[width=0.75\linewidth]{feasible-rays-edge.png}
%   \caption{The set of feasible ray directions at the edge of a tetrahedron.}\label{fig:feasible-rays-edge}
% \end{figure}

A consequence of our choice of tetrahedron updates is that $\mxhat$ (the $\mathtt{trial}$ update point) may not be adjacent to $\mx_1$ or $\mx_2$. Because of this, the local update ray $[\mxhat, \mx_{\mlam^*}]$ could potentially leave the tetrahedron mesh. The purpose of using an unstructured mesh is to have our spatial discretization conform to $\partial\Omega$, in order to accurately approximate the characteristics of the eikonal PDE. Then, $O(1)$ errors in $\grad\Eik$ can occur if  these updates are not rejected. One could reject these updates by carefully walking the cells of the tetrahedron mesh, from $\mx_{\mlam^*}$ to $\mxhat$, checking that each piece of the ray connecting $\mx_{\mlam^*}$ to $\mxhat$ is contained in $\Omega_h$. This approach is expensive and overly cautious.

As an alternative, we simply check whether $-\grad\eik(\mxhat)$ and $\grad\eik(\mx_{\mlam^*})$ are contained in the cones of feasible ray propagation directions at $\mxhat$ and $\mx_{\mlam^*}$. Consider a direction vector $\mt$ centered at a point $\mx$. We can check whether $\mt$ points into the mesh (i.e., there exists some $\epsilon > 0$, such that $[\mx, \mx + \epsilon \mt] \subseteq \Omega$) by handling the following three cases:
\begin{itemize}
\item \textbf{Case 1} ($\mx$ is a vertex of the mesh). Enumerate the cells incident on $\mx$. For each cell, let $\mx_1, \mx_2, \mx_3$ denote the remaining three vertices of the cell. Let $\mV = \begin{bmatrix} \mx_1 - \mx & \mx_2 - \mx & \mx_3  - \mx \end{bmatrix}$, and let $\malpha = \mV^{-1} \mt$. Then, $\mt$ points into the mesh if $\malpha$ has all nonnegative components.
\item \textbf{Case 2} ($\mx$ lies in the interior of an edge). In this case, we consider the fan of cells incident on the edge. Each cell defines a cone of feasible directions centered at $\mx$. For this case, we project $\mt$ and the cells incident on the edge into the edge's normal plane and then check whether $\mt$'s projection points into one of the projected cones.
\item \textbf{Case 3} ($\mx$ lies in the interior of a face). There are at most two cells incident on the face. For each cell, we let $\mz$ be the opposite vertex. If $(\mz - \mx) \cdot \mt \geq 0$, then $\mt$ points into the mesh.
\end{itemize}
Of course, if $\mx$ lies in the interior of a cell, then all ray directions $\mt$ point into the mesh.

\subsection{Accurate initialization near caustics}\label{ssec:accurate-initialization-near-caustics}

Solving the eikonal PDE with caustic BCs limits the accuracy of the scheme to $O(h \log \tfrac{1}{h})$~\cite{zhao2005fast}. Special first-order factoring schemes were initially designed to addressed this problem~\cite{fomel2009fast}, and were later extended to high-order schemes~\cite{luo2014high}. The source of this problem can be understood as follows. The accuracy of an update depends on the curvature of the wavefront over the base of the update. This curvature is singular in the vicinity of a caustic. As $h \to 0$, if we do an update in the immediate neighborhood (with respect to the mesh graph) of the caustic, this error grows rapidly, resulting in $O(h)$ error overall~\cite{potter2022numerical}. This wrecks the higher-order accuracy of the local updates used in this work. However, if we only do updates whose bases are at least $O(1)$ distance from the caustic, the original order of accuracy is maintained. We note that a first order local factoring has been developed in the context of complicated environments, but only in 2D~\cite{qi2019corner}.

As a simple alternative to factoring, we use the following approach. We use a Lagrangian scheme in the vicinity of the caustic to initialize each point within a small, constant radius of the caustic. We explain how to do this to initialize nodes near a point source, near a diffracting edge, and near a reflecting facet.

\paragraph{Initialization near a point source} We let the point source be denoted $\mxsrc$ and define a parameter $\rfac > 0$ which scales as a constant with respect to $h$. To simplify bookkeeping, we insert $\mxsrc$ into $\calV_h$ (we just request that TetGen insert this point into the mesh as a Steiner point~\cite{hang2015tetgen}). We initially set $T(\mxsrc) = 0$ and $\grad\Eik(\mxsrc) = \mathtt{undefined}$. We use a breadth-first search (BFS) to initialize nodes in the vicinity of $\mxsrc$. Each time we pop a node $\mxhat$ from the queue backing the BFS, we set $T(\mxhat) = \norm{\mxhat - \mxsrc}/c$ and $\grad\Eik(\mxhat) = (\mxhat - \mxsrc)/(c\norm{\mxhat - \mxsrc})$. If $c T(\mxhat) < \rfac$, then we insert each of $\mxhat$'s neighbors into the queue.

Afterwards, for each node $\mxhat$ which we have initialized this way, we check whether $\mxhat$ is only adjacent to other initialized nodes. If it is, we mark it $\mathtt{valid}$. Otherwise, we mark it $\mathtt{trial}$ and insert it into the heap. In this way, we leave the JMM in a logically consistent state so that we can begin marching.

This is a form of Lagrangian ray marching. Although we have presented the idea here for constant $c$, it can be adapted to a varying speed of sound by solving two-point BVPs at each step. The reason for using a BFS is that it allows one to check whether the expanding front contacts some part of $\partial\Omega$. A drawback of this approach is that if the domain is too irregular, it may not be possible to march the solution through a sufficiently large neighborhood to control the error for practically sized meshes.

\paragraph{Initialization near a diffracting edge} As before, we use a BFS where we update each node $\mxhat$ using a fully Lagrangian update. The main difference is how this update is done. To initialize $T$ around a diffracting edge, we clearly have some knowledge of the eikonal of a GO wavefront incident on the edge, call it $\Eikin$. For each node $\mx$ on the diffracting edge (which could itself consist of a number of edges in $\calE_h$), we set $T(\mx) = \Eikin(\mx)$. If $\mt$ is a unit tangent vector for the edge, we can see from \Cref{ssec:BCs-edge-diffraction} that $\mt \cdot \grad\Eik(\mx)$ is well-defined for each $\mx$ on the diffracting edge. However, the directional derivative of $\Eik$ is undefined for any other direction. To get around this issue, we initialize $T$ using a piecewise cubic Hermite polynomial on the edge, with data given by $T(\mx)$ and $\mt \cdot \grad\Eik(\mx)$. This puts us in a position to do triangle updates (\Cref{ssec:triangle-updates}). So, in order to update each node $\mxhat$ as we run the BFS, we simply do triangle updates from each $(\mx_0, \mx_1) \in \calE_h$.

\paragraph{Initialization near a reflection} When we compute a reflection, we will also need to initialize $\Eik$ on some reflecting facet in $\partial\Omega_h$. In this case, we can extend the algorithm for initializing a diffracting edge by doing tetrahedron updates from the triangular facets comprising the reflector. Note that this means that both reflection and diffraction updates are taken together to ensure that rarefaction fans are modeled correctly.

\paragraph{Fast algorithms for initialization} In the case of a diffracting edge, if we do every triangle update from the diffracting edge to update a node $\mxhat$, we can expect this to take $O(N^{4/3})$ time, since there will be $O(N^{1/3})$ nodes on the diffracting edge, and $O(N)$ nodes in the tubular neighborhood of radius $\rfac$ surrounding it. Along the same lines, the cost to initialize from a reflecting face is $O(N^{5/3})$. We note that although these costs seem modest, they quickly dominate the overall runtime of the algorithm as the mesh is refined. To mitigate these costs, we keep track of the ``parent'' of each node inserted into the BFS, as well as the update which led to its final value---we use this update as a warm start in the search for a new update. This can be combined with another graph search to efficiently locate a new update. In practice, we observe that this significantly reduces the cost of accurately initializing $\Eik$ and $\grad\Eik$ in a tubular neighborhood.

\subsection{Jet marching on an unstructured tetrahedron mesh}\label{sec:JMM-algorithm}

Finally, following the discussion in the preceding sections, we arrive at the following algorithm:
\begin{enumerate}
\item Follow \Cref{ssec:accurate-initialization-near-caustics} to initially set the state and jet of nodes involved in setting the BCs.
\item Sort the $\mathtt{trial}$ nodes into a priority queue in increasing order of $T$.
\item While not all nodes are $\mathtt{valid}$:
  \begin{enumerate}
  \item Pop the first node (call it $\mx_0$) from the prioritiy queue. Set the state of $\mx_0$ to $\mathtt{valid}$, and free the old entries in the update caches associated with $\mx_0$.
  \item Set the state of each $\mathtt{far}$ neighbor of $\mx_0$ to $\mathtt{trial}$ and insert it into the priority queue.
  \item For each $\mathtt{trial}$ neighbor of $\mx_0$:
    \begin{enumerate}
    \item If $\{\mx_0, \mxhat\}$ is an edge of $\boundary\domain_h$, do each boundary triangle update.
    \item If $\mx_0$ is on a diffracting edge, do the incident edge-diffraction triangle updates.
    \item If $\mx_0$ isn't on a diffracting edge, do the fan of freespace tetrahedron updates.
    \end{enumerate}
    In each case, for each physical update, if the corresponding minimization problem has an interior point solution, accept it. Otherwise, search through the relevant update cache and look for a matching update---if one is found, accept the update. Otherwise, cache the current update.
  \end{enumerate}
\end{enumerate}
This is a high-level view of the actual implementation used. There are a number of low-level details which must be accounted for which would obscure the presentation if they were included, and which could easily vary with the implementation of the algorithm.

\paragraph{Algorithm complexity} We give a rough estimate of the complexity of this algorithm. In what follows, we let $N_h = |\calV_h|$. We observe that:
\begin{enumerate}
\item The cost of maintaining the priority queue is $O(N_h \log N_h)$.
\item The cost of maintaining the update caches is also $O(N_h \log N_h)$.
\item The cost of doing individual updates is $O(1)$. The reason for this is that the cost of the individual updates does not increase as $N_h \to \infty$. SQP converges linearly or quadratically; hence, 5--20 iterations if machine precision is sought. Nevertheless, we cannot do better than machine precision, so the total number of iterations per update is capped. The amount of work per iteration depends on the dimension of $\Omega$, which is constant---hence, $O(1)$.
\item Each node will be updated $O(1)$ times. This follows from our assumption that each node has $O(1)$ neighbors.
\item Altogether, $O(N_h)$ updates will be done.
\end{enumerate}
Hence, the asymptotic complexity stems from three major sources: 1) priority queue maintenance, 2) update cache maintenance, 3) doing updates. The first two points are both $O(N_h \log N_h)$, while the latter is $O(N_h)$. The asymptotic prefactor of doing updates is by far the largest, followed by maintaining the update caches, and finally the priority queue. Because of this, for most problems of practical size, the scaling is essentially linear in the number of nodes, until the problem becomes extremely large~\cite{potter2019ordered}. Since the purpose of this solver is at least in part to reduce the memory costs of an $O(h)$ eikonal solver by using a more accurate discretization, it is unlikely that such large problem sizes will be encountered in practice.

\subsection{Storing the dynamic programming plan and computing origins}\label{ssec:approximate-origins}

As we run the JMM, we track each node $\mxhat$'s \emph{parent}---the set of vertices which updated it---and the convex coefficient $\mlam^*$ giving the optimum $\mx_{\mlam^*}$. This information together with the order in which nodes were relaxed constitutes the dynamic programming plan discovered by the label-setting algorithm backing the JMM. A useful fact is that it can be used to transport quantities downwind from a set of nodes with BCs. We will use this technique more in \Cref{sec:amplitude}, but we start with a first application in this section: computing what we refer to as the \emph{origin} of a node.

We define a new function $\mathtt{org} : \calV_h \to [0, 1]$. Our intent is for the $1/2$-level set of $\mathtt{org}$ to approximate the reflection and shadow zone GO boundaries, where if $\mathtt{org}(\mx) < 1/2$, it lies inside the \emph{shadow zone}---the set of points reached indirectly by a diffracted ray. By setting $\mathtt{org}(\mx) = 1$ for each node $\mx$ with BCs (a point source, a vertex on a diffracting edge, or a vertex in the interior of a reflecting facet), and setting $\mathtt{org}(\mx) = 0$ for each node $\mx$ without BCs but which lies on a diffracting edge. We then \emph{replay} the dynamic programming plan, iterating over nodes $\mxhat\in\calV_h$ in order of acceptance, setting $\mathtt{org}(\mxhat) = \sum_i \lam_i \mathtt{org}(\mx_i)$ if $\mxhat$ hasn't received a value yet, where the $\lam_i$'s and $\mx_i$'s are the parent convex coefficients and vertices of $\mxhat$. Afterwards, for each node $\mx$ on a diffracting edge with $\mathtt{org}(\mx) = 0$ initially, we set $\mathtt{org}(\mx) = 1/2$ if any of $\mx$'s parent nodes have an $\mathtt{org}$ value greater than $1/2$ (i.e., lie outside the shadow zone).

\subsection{Reinitialization in the shadow zone}\label{ssec:reinit-shadow-zone}

For the reasons outlined in \Cref{ssec:accurate-initialization-near-caustics}, if a node lies in the shadow zone, its order of accuracy can be at most $O(h)$ since it will be reached by a chain of updates which passes through a caustic. To maintain the order of accuracy when computing the first arrival time, we propose the following simple strategy:
\begin{enumerate}
\item After solving the eikonal PDE using the algorithm of \Cref{sec:JMM-algorithm}, reset the value of each node in the shadow zone. That is, reset its state to $\mathtt{far}$ and set $T$ to infinity.
\item Add diffraction BCs to each edge incident on the shadow zone.
\item Use the BFS algorithm in \Cref{ssec:accurate-initialization-near-caustics} to do Lagrangian initialization of nodes in an $O(1)$ neighborhood around the diffracting edge.
\item To maintain the correctness of the $\mathtt{valid}$ front, for each $\mathtt{valid}$ node which is now adjacent to a $\mathtt{far}$ node, reset its state to $\mathtt{trial}$ and reinsert it into the heap.
\item Continue marching.
\end{enumerate}
This algorithm can be used to incrementally correct the accuracy of the first arrival time field deeper and deeper into shadow zone. We leave developing an optimal complexity algorithm which avoids the need for multiple passes for future work.

\section{Post-processing: computing the amplitude}\label{sec:amplitude}

\subsection{Computing the Hessian using cell averaging}

% \begin{figure}
%   \includegraphics[width=\linewidth]{cell-averaging.png}
%   \caption{A schematic depiction of cell-averaging in 2D. We approximately determine the locating of $\boundary\shadow$ as the $\tfrac{1}{2}$-level set of $\mathtt{origin}(\mx)$. Hence, we in principle possess the information necessary to initialize $\hess\Eik$ correctly for averaging on either side of $\boundary\shadow$.}\label{fig:cell-averaging}
% \end{figure}

The 20-parameter Bernstein-\Bezier{} element has the virtue of being cheap and simple. It is also $O(h^{3-p})$ accurate, which is sufficient for our purposes. Unfortunately, it is only a $C^0$ spline. The amplitude can be determined from the curvature of the wavefront at each point, where the wavefront at time $\hat\tau$ corresponds to the $\hat\tau$-level set of $\tau$. Because of the low degree of continuity of the spline defining $T$, we need a way to evaluate $\hess\Eik(\mx)$ for each $\mx\in\calV_h$. Our approach is to average the value of $\hess\Eik(\mx)$ obtained from the 20-parameter element defining $T$ over each adjacent cell. However, there are several corner cases that can occur which can degrade the accuracy of the average. Specifically:
\begin{itemize}
\item Since $T$ is singular at caustics, we must manually set the value of $\hess\Eik(\mx)$ for $\mx\in\calV_h$ in the vicinity of a caustic. The algorithms in \Cref{ssec:accurate-initialization-near-caustics} provide the set of nodes for which this must be done as a byproduct. We note that caustics can also form for a constant speed of sound if focusing occurs (for instance, if a wave is reflected from the interior of a convex surface), but we leave the treatment of this type of caustic for a future work.
\item If $\mx$ was updated from a diffracting edge, we do not average over cells which are incident on the diffracting edge. Contrariwise, if $\mx$ is incident on a diffracting edge, we do not use cells which contain vertices which were updated from that diffracting edge. The logic here is that we should respect the fact that the Hessian is discontinuous as we traverse a ray into a shadow zone. If we average values of $\hess\Eik$ over nodes which are inside and outside of the shadow zone, the approximation will be inaccurate.
\end{itemize}

\subsection{Computing the spreading factors}

After approximating $\hess\Eik(\mx)$ for each $\mx\in\calV_h$, we can approximate the spreading factor at each point using the cached dynamic programming plan described in \Cref{ssec:approximate-origins}. For a point $\mxhat \in \calV_h$, assume that we have an approximation of the amplitude at each parent of $\mxhat$; i.e., we have approximated $A(\mx_i)$ for each $\mx_0, \mx_1$, and $\mx_2$. Let $\mlam^*$ be optimum. To approximate $A(\mx_{\mlam}^*)$, we first linearly interpolate the logarithm of the amplitude:
\begin{equation}
  \log A(\mx_{\mlam^*}) = (1 - \lam_1^* - \lam_2^*)\log A(\mx_0) + \lam_1^* \log A(\mx_1) + \lam_2^* \log A(\mx_2),
\end{equation}
and then exponentiate. Equivalently, we set:
\begin{equation}
  A(\mx_{\mlam^*}) = A(\mx_0)^{1 - \lam_1^* - \lam_2^*} A(\mx_1)^{\lam_1^*} A(\mx_2)^{\lam_2^*}.
\end{equation}
Approximating the amplitude this way is better motivated than linearly interpolating the amplitude directly, since perception of amplitude is logarithmic.

To propagate the amplitude along the update ray, we approximate the principal curvatures of the wavefront passing through $\mxhat$ with the first two largest magnitude eigenvalues of $\hess\Eik(\mxhat)$, giving approximations of $\kappa_1(\mxhat)$ and $\kappa_2(\mxhat)$ in \eqref{eq:amplitude-transport}. We then approximate the integral in \eqref{eq:amplitude-transport} using the left-hand rule, giving:
\begin{equation}
  A(\mxhat) = A(\mx_{\mlam^*}) \exp\parens{\frac{L}{2}(\kappa_1 + \kappa_2)}, \qquad L = \norm{\mxhat - \mx_{\mlam^*}},
\end{equation}
where the principal curvatures $\kappa_1$ and $\kappa_2$ are computed from \eqref{eq:principal-curvatures}.

After computing $A(\mx)$ for each $\mx\in\calV_h$, we can evaluate $A$ for a point $\mx\in\domain_h$ which is not in $\calV_h$ by first determing the cell $C = (\mx_0,\mx_1,\mx_2,\mx_3) \in \calC_h$ which contains $\mx$, computing the coordinates $\mlam = (\lam_1, \lam_2, \lam_3)$ of $C$ with respect to its vertices:
\begin{equation}
  \mx_{\mlam} = \mx_0 + \lam_1 (\mx_1 - \mx_0) + \lam_2 (\mx_2 - \mx_0) + \lam_3 (\mx_3 - \mx_0),
\end{equation}
and setting:
\begin{equation}
  A(\mx_{\mlam}) = A(\mx_0)^{1-\lam_1-\lam_2-\lam_3}A(\mx_1)^{\lam_1} A(\mx_2)^{\lam_2}A(\mx_3)^{\lam_3}.
\end{equation}

\subsection{Evaluating the UTD edge diffraction coefficients}\label{ssec:evaluating-UTD-coefs}

At each point in the observation field, the parameters of the UTD edge diffraction coefficients are either global constants, or functions which depend on the point of diffraction. Our strategy is to use the precomputed dynamic programming plan to transport these parameters from the diffracting edge along the rays to each field point. We only consider planar wedges, in which case the wavenumber $k$ and the wedge parameter $n$ are global constants. The remaining parameters---namely, $\mtin$ and $\mtout$---must be transported. For a point $\mx$ reached by an edge-diffracted ray, if $\mx_e$ is the point of diffraction, then $\mtin(\mx)$ is the unit tangent vector of the ray which originally impinged on $\mx_e$, evaluated at $\mx_e$; $\mtout(\mx)$ is the unit tangent vector of the ray giving the direction in which the ray initially diffracted from the edge.

To transport $\mtin$ and $\mtout$, we again make use of the precomputed dynamic programming plan discussed in \Cref{ssec:approximate-origins}. If $\mx_0, \mx_1$, and $\mx_2$ are the parents of $\mxhat$, with $\mlam^* = (\lam_1, \lam_2)$ optimum, and $\mt = \mt(\mx)$ is a unit vector field we could simply set:
\begin{equation}
  \mt(\mxhat) = \frac{(1 - \lam_1^* - \lam_2^*) \mt(\mx_0) + \lam_1^* \mt(\mx_1) + \lam_2^* \mt(\mx_2)}{\norm{(1 - \lam_1^* - \lam_2^*) \mt(\mx_0) + \lam_1^* \mt(\mx_1) + \lam_2^* \mt(\mx_2)}}
\end{equation}
to interpolate at each step. Unfortunately, this sort of normalized linear interpolation is inaccurate near caustics. Instead, we compute a spherical weighted average of the parents' unit vectors~\cite{buss2001spherical} (i.e., we compute their Riemannian mean). Namely, if $\mq_1, \hdots, \mq_n$ are unit vectors, we think of them as points on the sphere, and let $d_g$ be the corresponding geodesic distance. This is just the angle subtended by a pair of points, which is well-defined if they lie in the same hemisphere. In general, if $\eta_1, \hdots, \eta_d$ are convex coefficients such that $\sum_{i=1}^d \eta_i = 1$ and $\eta_i \geq 0$ for $i = 1, \hdots, n$, then we consider the energy defined by:
\begin{equation}
  f(\mq) = \frac{1}{2}\sum_{i=1}^n \eta_i d_g(\mq, \mq_i)^2 = \frac{1}{2} \sum_{i=1}^n \eta_i \cos^{-1}\parens{\mq_i\cdot\mq}^2.
\end{equation}
Its gradient is given by:
\begin{equation}
  \grad f(\mq) = -\sum_{i=1}^n \eta_i \frac{\cos^{-1}(\mq_i\cdot\mq)}{\sqrt{1 - \parens{\mq_i\cdot\mq}^2}}\mq_i = -\sum_{i=1}^n \eta_i \frac{\theta_i(\mq)}{\norm{\mq - \operatorname{proj}_{\mq_i}\mq}}\mq_i
\end{equation}
The weighted average is then the unit vector $\mq$ solving:
\begin{equation}
  \begin{split}
    \mbox{minimize} \qquad& f(\mq) \\
    \mbox{subject to} \qquad& \norm{\mq} = 1.
  \end{split}
\end{equation}
We use a simple projected gradient descent algorithm to solve this optimization problem. We use normalized linear interpolation to compute a warm start, take a descent step, and then normalize to project back onto the sphere. The iteration takes the form:
\begin{equation}
  \mq_0 = \frac{\sum_{i=1}^n \eta_i \mq_i}{\norm{\sum_{i=1}^n \eta_i\mq_i}}, \qquad \tilde{\mq}_{i+1} = \mq_i + \sum_{i=1}^n \eta_i  \frac{\cos^{-1}(\mq_i\cdot\mq)}{\sqrt{1 - \parens{\mq_i\cdot\mq}^2}}\mq_i, \qquad \mq_{i+1} = \frac{\tilde{\mq}_{i+1}}{\norm{\tilde{\mq}_{i+1}}}.
\end{equation}
This iteration converges reliably, typically with a linear rate. More sophisticated optimization algorithms which attain a quadratic rate of convergence are available~\cite{buss2001spherical,absil2009optimization}.

\paragraph{Special case: $n = 2$} For a triangle update, the spherical weighted combination is simple to compute, and is given by the ``slerp'' formula~\cite{shoemake1985animating}:
\begin{equation}
  \mt(\mxhat) = \frac{\sin\big((1 - \lam^*)\theta\big)}{\sin(\theta)}\mt(\mx_0) + \frac{\sin(\lam^*\theta)}{\sin(\theta)}\mt(\mx_1),
\end{equation}
where $\theta$ is the angle between $\mt(\mx_0)$ and $\mt(\mx_1)$.

\paragraph{Initializing $\mtout$} Note that while $\mtin$ is well-defined on the diffracting edge---it is simply the normalized gradient of the incident field---$\mtout$ is not. The field of rays diverges in multiple directions along Keller's cone at each point on the diffracting edge. However, if we initialize $\mtout$ in a tube of radius $\rfac$ around the diffracting edge, analogous to what is done in \Cref{ssec:accurate-initialization-near-caustics}, we can avoid this problem. In practice, we initialize both $\mtin$ and $\mtout$ in this way if we can, since it increases the accuracy of the approximation.

\section{Test case: a semi-infinite sound hard wedge}\label{sec:wedge}

\subsection{Numerical solution}

High-frequency scattering from a semi-infinite sound hard wedge is a test case with an analytically known solution which we can use to test our solver. We let $w > 0$ and $h > 0$, and initially set $\Omega = [-w/2, w/2] \times [-w/2, w/2] \times [-h/2, h/2]$. If $\phi$ is the polar angle of $(x, y)$ in the $xy$-plane, then we further require $0 \leq \phi \leq n \pi$, where $n$ is the parameter defining the wedge. That is, the diffracting edge is the $z$-axis. We let $\mt_e = (0, 0, 1)$ be a unit tangent vector for the edge. The faces of the wedge are vertical half-planes aligned with the angles $\phi = 0$ and $\phi = n\pi$ in the $xy$-plane, respectively. Note that in the UTD literature, the $\phi = 0$ is referred to as the \emph{o-face} and the $\phi = n\pi$ face is referred to as the \emph{n-face}, for the obvious reason.

In exterior Helmholtz problems, it is standard to consider an incident field $u_{\In}$, to let $u = u_{\In} + u_{\operatorname{scat}}$, and to solve for the scattered field $u_{\operatorname{scat}}$. The approach here is somewhat different. Again, we assume an incident field $u_{\In}$, but we must instead specify the BCs for $\Eik$ and $A$ on $\BCset_h$. We restrict our attention to point sources in our numerical experiments. For a point source at $\mxsrc \in \Omega$, we pick a radius $\rfac > 0$ and initialize $\Eik$ and $\grad\Eik$ following \Cref{ssec:accurate-initialization-near-caustics}, and $A$ following \Cref{sec:amplitude}. After solving the point source problem, we correct the solution in the shadow zone following \Cref{ssec:reinit-shadow-zone}. We also compute reflected and edge-diffracted fields, with BCs set according to \Cref{ssec:accurate-initialization-near-caustics}.

\subsection{Exact solution}\label{ssec:wedge-exact-solution}

% \begin{figure}
%   \centering
%   \includegraphics[width=0.75\linewidth]{n-refl-zones.png}
%   \caption{The angles involved in edge diffraction. The wedge angle $\alpha > 0$ is such that $\alpha = (2 - n)\pi$. We assume that the wedge parameter $n$ is in the range $(1, 2)$. The dashed line indicates the $\mathtt{origin}(\mx) = 1/2$ level set. This corresponds to the shadow boundary $\boundary\shadow$ for the direct field, the reflection boundary $\partial R$ for the $o$-face reflection, and is not present for the $n$-face reflection for this problem, since the $n$-face reflection is a pure edge diffraction problem in this case. Different rows correspond to different tetrahedron meshes (coarsest to finest from top to bottom). Left: the direct field. Middle: the $o$-face reflection. Right: the $n$-face reflection.}\label{fig:refl-zones}
% \end{figure}

We explain how to evaluate $\Eik, \grad\Eik$, and $\hess\Eik$ for the incident field due to the point source, reflected fields, and edge-diffracted fields.

\paragraph{Direct field due to point source} There are two cases to check. If $\mx$ is directly reached by a ray from the point source, then set the solution according to \Cref{ssec:free-space-point-source}. Otherwise, $\mx$ lies in the shadow zone. Let $\mx_e$ be a point on the diffracting edge, and let:
\begin{equation}\label{eq:tmp2}
  T(\mx) = \frac{1}{c} \min_\lambda \left\{\norm{\mx - \mx_e - \lambda \mt} + \norm{\mx_e + \lambda \mt - \mxsrc}\right\}.
\end{equation}
We can then set:
\begin{equation}
  \grad\Eik(\mx) = \frac{1}{c} \frac{\mx - \mx_e - \lambda^* \mt}{\norm{\mx - \mx_e - \lambda^* \mt}},
\end{equation}
where $\lambda^*$ is the minimizing argument in \eqref{eq:tmp2}. To determine $\hess\Eik$, recall that $\hess\Eik\grad\Eik = 0$ if $c$ is constant. Then, from \eqref{eq:c-hess-tau-eig-decomp}, we have:
\begin{equation}
  c \hess\Eik(\mx) = \kappa_1 \mq_1 \mq_1^\top + \kappa_2 \mq_2 \mq_2^\top,
\end{equation}
where $\kappa_1$ and $\kappa_2$ are the principal curvatures of the GO wavefront passing through $\mx$, and $\mq_1$ and $\mq_2$ are the associated principal directions. For edge diffraction, we can determine one of the principal directions by first letting the projection of $\mx$ onto the diffracting edge be $\mxproj$ given by \eqref{eq:xproj}, and then considering the circle of radius $\rhoe = \norm{\mx - \mxproj}$ lying in the plane perpendicular to the diffracting edge passing through $\mx_e$. One of the principal directions is the tangent vector of this circle, with radius of curvature given by $\rhoe$~\cite{mcnamara1990introduction}. Call this direction $\mq_2$, since it will in general be the smaller of the two. Consequently, the other principal direction will be orthogonal to $\mq_2$ and $\grad\Eik$, and the radius of curvature is just $c T(\mx)$. Hence, $\kappa_1 = 1/(cT(\mx))$ and $\kappa_2 = 1/\rhoe$.

\paragraph{Scattered fields} Handling edge-diffraction is essentially the same as the case of a field point lying in the shadow zone for the point source. For a reflection, we can use the method of images. We reflect $\mxsrc$ over the face, which defines a reflection boundary. If a point lies on the near side of the reflection boundary, we set its data according to \Cref{ssec:free-space-point-source}. Otherwise, it lies in the shadow boundary, and we follow the approach outlined in the preceding paragraph.

\subsection{Results}

\begin{figure}
  \centering
  \includegraphics[width=\linewidth]{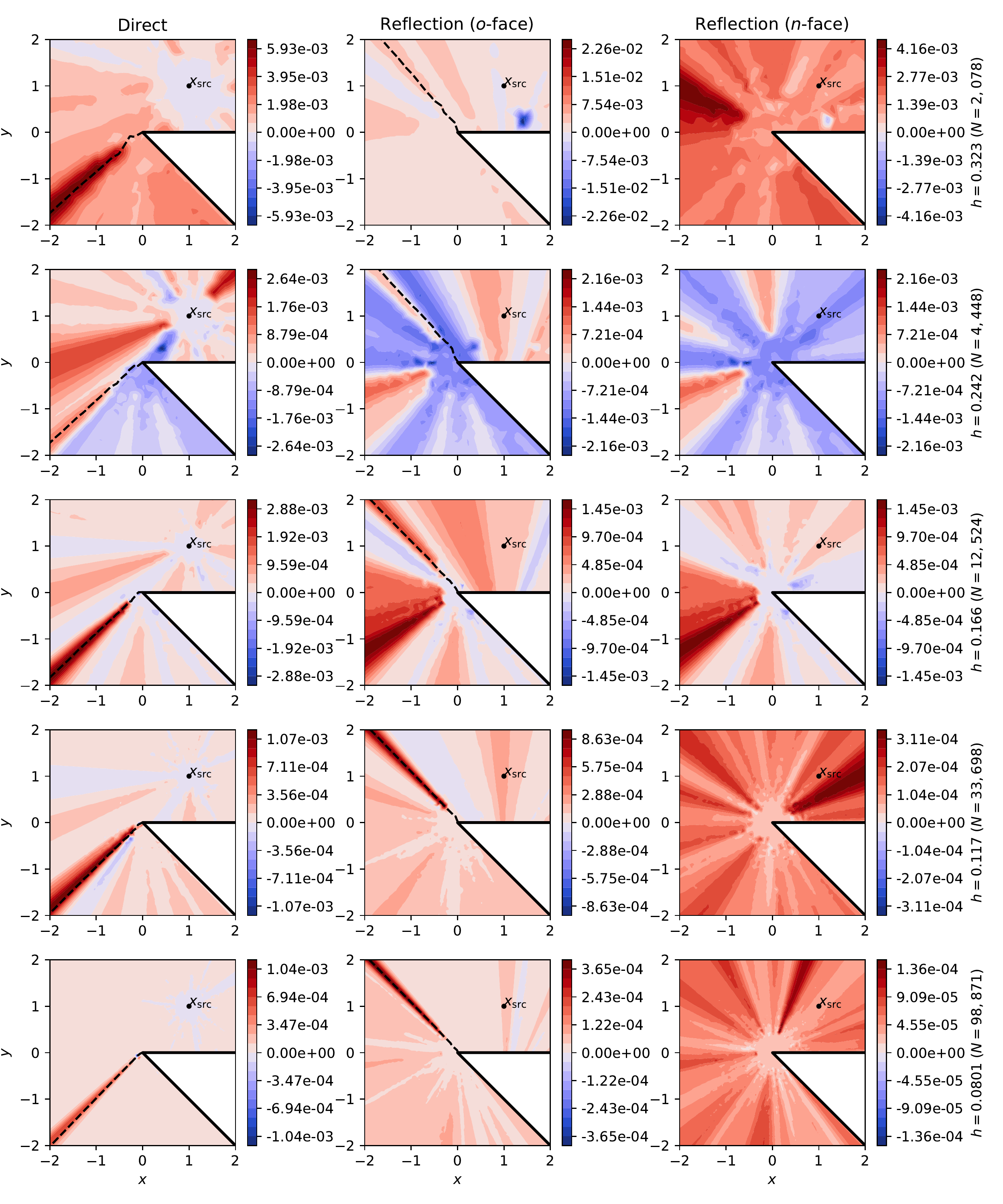}
  \caption{Pointwise error plots (horizontal slices with $z=0$) for $T(\mx)$ for the edge-diffraction problem with $n = 1.75$, $w = 4$, $h = 2$, and $\mxsrc = (1, 1, 0)$. Errors are computed by building 20-parameter Bernstein-\Bezier{} splines with nodal function values given by $\tau(\mx) - T(\mx)$ and nodal gradient values given by $\grad\eik(\mx) - \grad\Eik(\mx)$ for each $\mx\in\calV_h$. The dashed lines correspond to the $\mathtt{org}(\mx) = 1/2$ level set. In the left column, for the direct field, it approximates the shadow boundary, and in the middle column it approximates the reflection boundary. There is no GO boundary for the rightmost column corresponding to pure edge diffraction, since $\mathtt{org} \equiv 1$ in that case.}\label{fig:error-field-plots}
\end{figure}

\begin{figure}
  \centering
  \includegraphics[width=\linewidth]{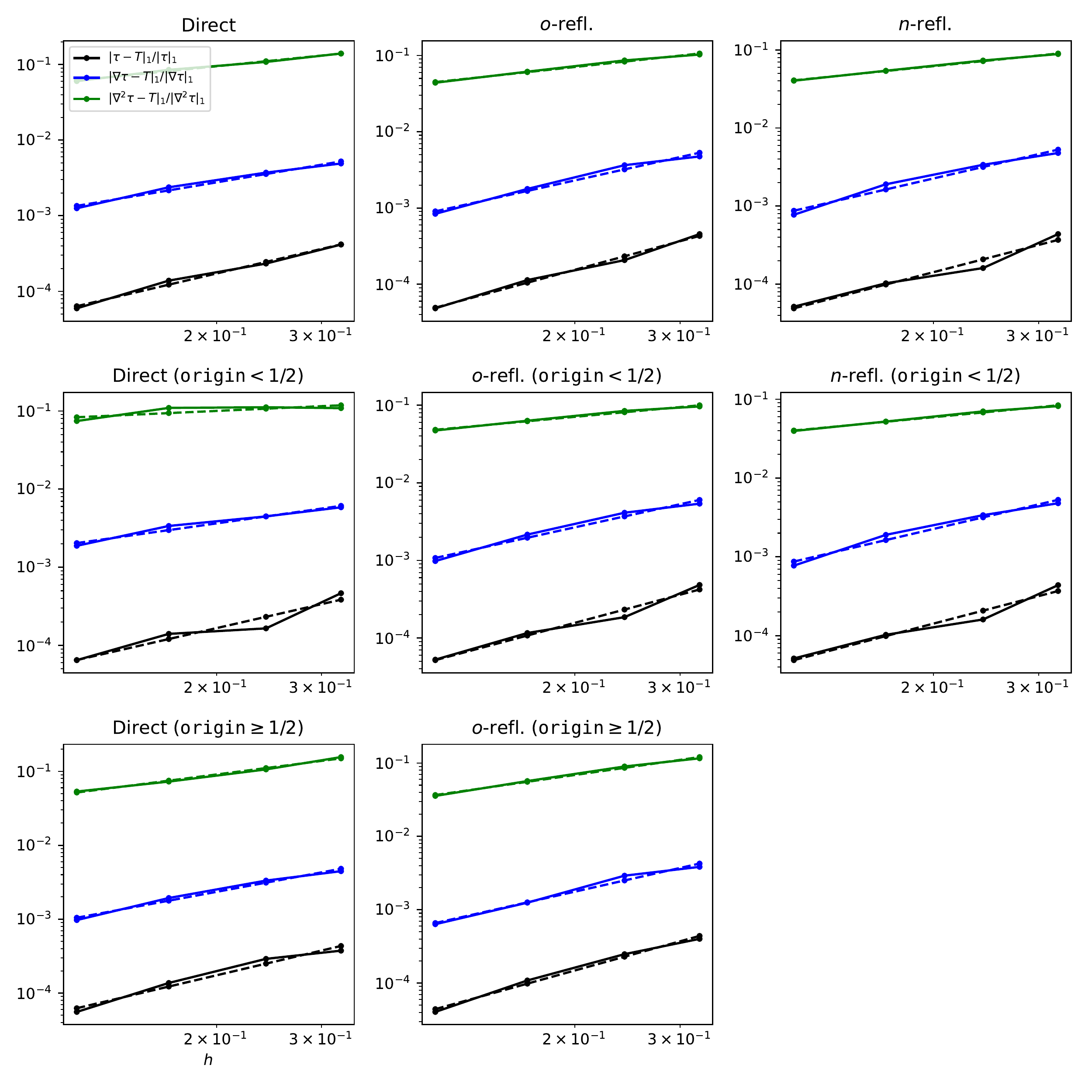}
  \caption{Relative $\ell_1$ errors for varying mesh finenesses for the edge-diffraction problem considered in \Cref{fig:error-field-plots}. Left column: the direct field. Middle column: the reflection from the $o$-face. Right column: the reflection from the $n$-face. Top row: the entire field. Middle row: the diffracted part of the field ($\mathtt{origin}(\mx) < 1/2$). Bottom row: the direct part of the field ($\mathtt{origin}(\mx) \geq 1/2$).}\label{fig:wedge-errors}
\end{figure}

\begin{figure}
  \includegraphics[width=\linewidth]{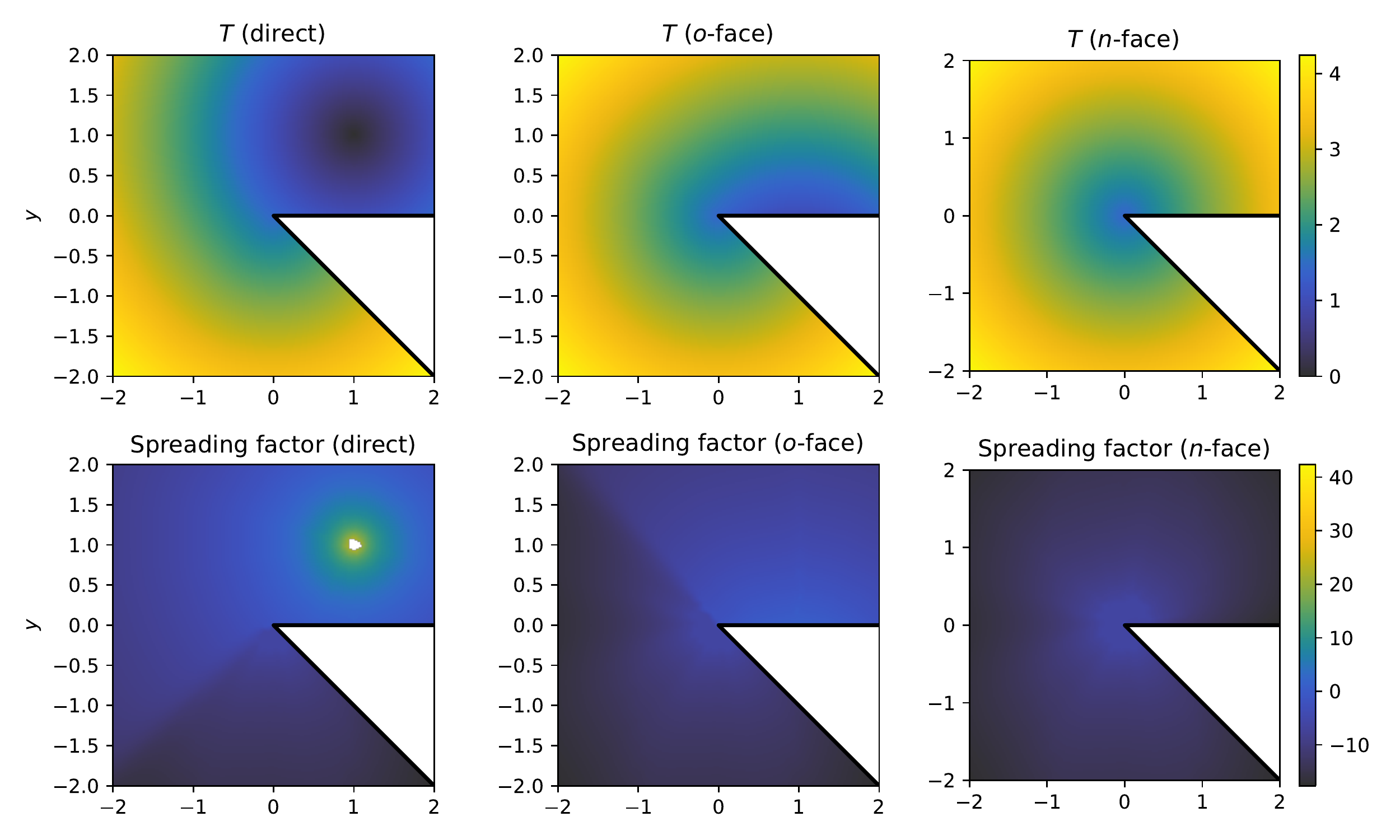}
  \caption{Top row: the eikonal for each branch of the wedge problem. Bottom row: the spreading. Since the spreading is singular at one vertex for the direct branch, we do not compute the spreading over that tetrahedron (hence the white patch in the bottom left subplot). Near the point source, a local model can be used instead.}\label{fig:T-and-amp-fields}
\end{figure}

% \begin{figure}
%   \includegraphics[width=\linewidth]{u-fields.pdf}
%   \caption{\sfp{...}}
% \end{figure}

\begin{figure}
  \centering
  \includegraphics[width=0.6\linewidth]{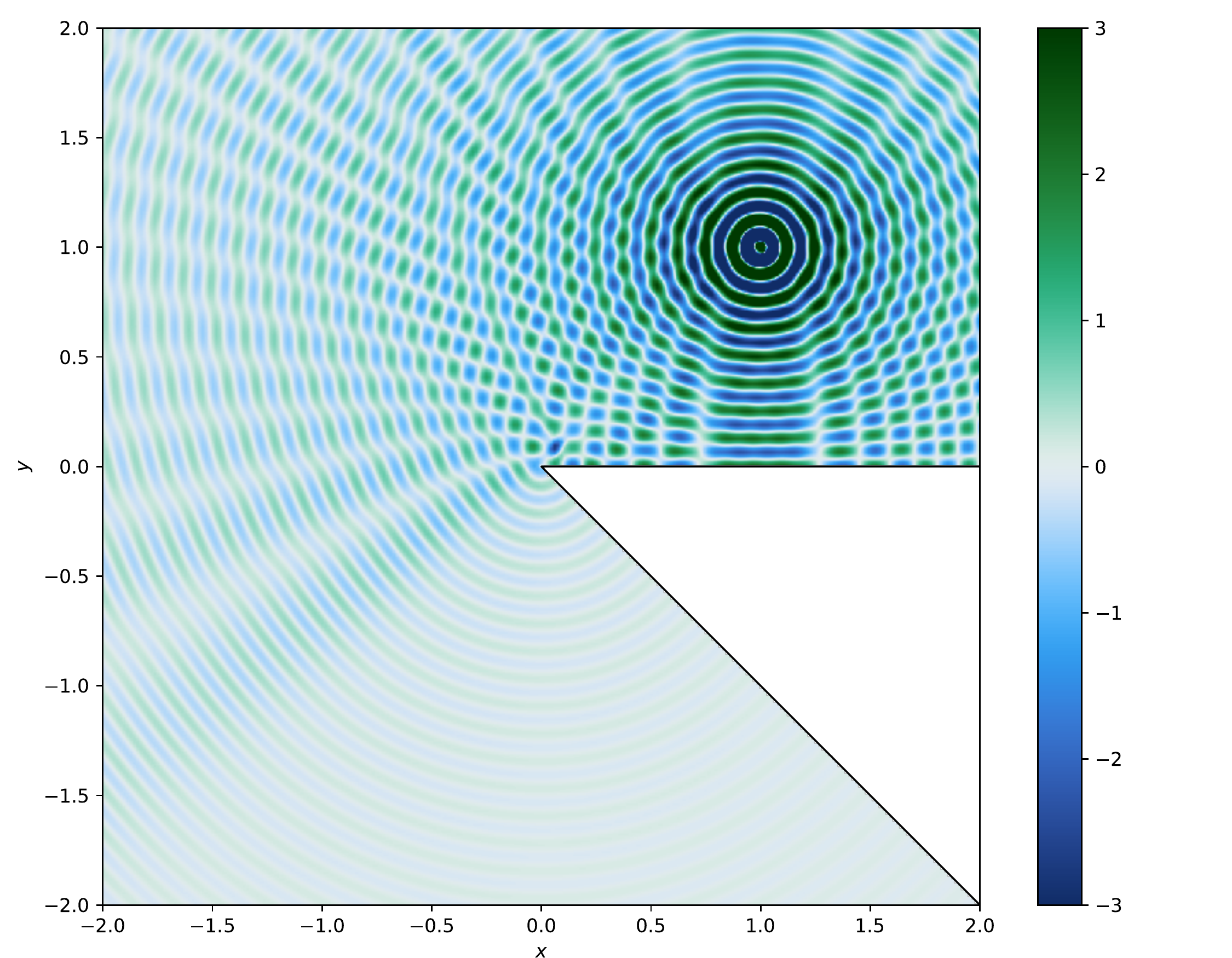}
  \caption{The total field ($u = u_i + u_o + u_n$) computed numerically using the jet marching method for the wedge diffraction problem shown in \Cref{fig:error-field-plots}.}\label{fig:total-field}
\end{figure}

\begin{table}
  \centering
\begin{tabular}{c|ccc|ccc|ccc}
  & \multicolumn{3}{c|}{All} & \multicolumn{3}{c|}{Direct} & \multicolumn{3}{c}{Diffracted} \\
  & $T$ & $\nabla T$ & $\nabla^2 T$ & $T$ & $\nabla T$ & $\nabla^2 T$ & $T$ & $\nabla T$ & $\nabla^2 T$ \\
  \midrule
  Direct & 1.90 & 1.42 & 0.80 & 1.92 & 1.52 & 0.97 & 1.84 & 1.25 & 0.47 \\
  Reflection ($o$-face) & 2.23 & 1.99 & 1.18 & 2.54 & 2.33 & 1.82 & 2.06 & 1.75 & 0.69 \\
  Reflection ($n$-face) & 2.00 & 1.92 & 0.75 & & & & 2.00 & 1.92 & 0.71 \\
\end{tabular}
\caption{Least squares fits for the relative $\ell_1$ norm error for $\Eik$, $\grad\Eik$, and $\hess\Eik$, measured for $h = {0.0801, 0.117, 0.166, 0.242, 0.323}$ ($N = {98871, 33698, 12524, 4448, 2078}$). First column (``All''): the error measured in all of $\domain_h$. Second column (``Direct''): for $\mx \in \calV_h$ with $\mathtt{origin}(\mx) > 0.5$. Third column (``Diffracted''): for $\mx \in \calV_h$ with $\mathtt{origin}(\mx) \leq 0.5$.}\label{table:errors}

\end{table}

We set $w = 4$, $h = 2$, $n = 1.75$, $\mxsrc = (1, 1, 0)$, and $\rfac = 0.3$. To conduct an experimental error analysis, we discretize the wedge into a sequence of 5 tetrahedron meshes with maximum volume constraints given by $0.01, 0.00316, 0.001, 0.000316, 0.0001$. We then compute the direct arrival, a reflection from the $o$-face, and an edge-diffracted field. In this case, the $n$-face lies in the shadow zone, and doesn't excite a scattered field.

To check the quality of the pointwise error, we construct the tetrahedral spline approximating $T$ from $T$ and $\grad\Eik$ (\Cref{ssec:bernstein-bezier}) and compare it on the $xy$-plane with the exact solution (\Cref{ssec:wedge-exact-solution}). See \Cref{fig:error-field-plots,fig:T-and-amp-fields}. We  also compute least squares fits of the relative $\ell_1$ error in $\Eik$, $\grad\Eik$, and $\hess\Eik$ (\Cref{fig:wedge-errors}) and the high-frequency approximation of the solution of the Helmholtz equation computed by superimposing the part of the field due to each branch (\Cref{fig:total-field}).

\section{Test case: a building interior with simplified geometry}\label{ssec:building}

\begin{figure}
  \includegraphics[width=0.6\linewidth]{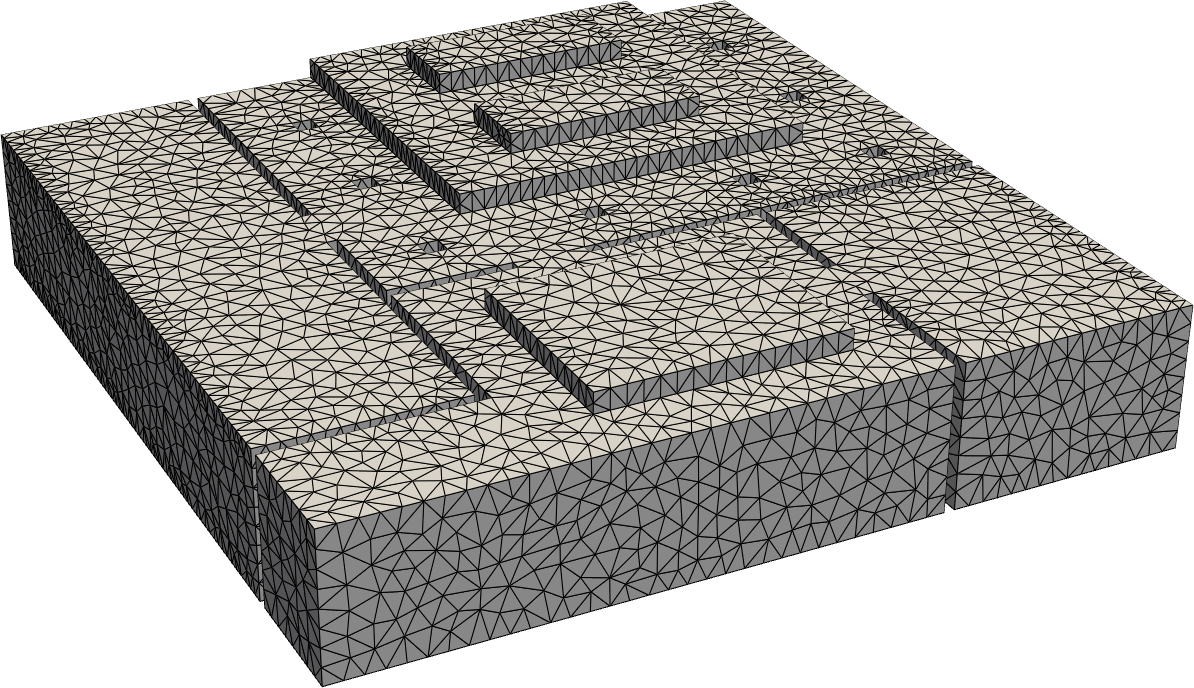}%
  \hspace{0.025\linewidth}%
  \includegraphics[width=0.35\linewidth]{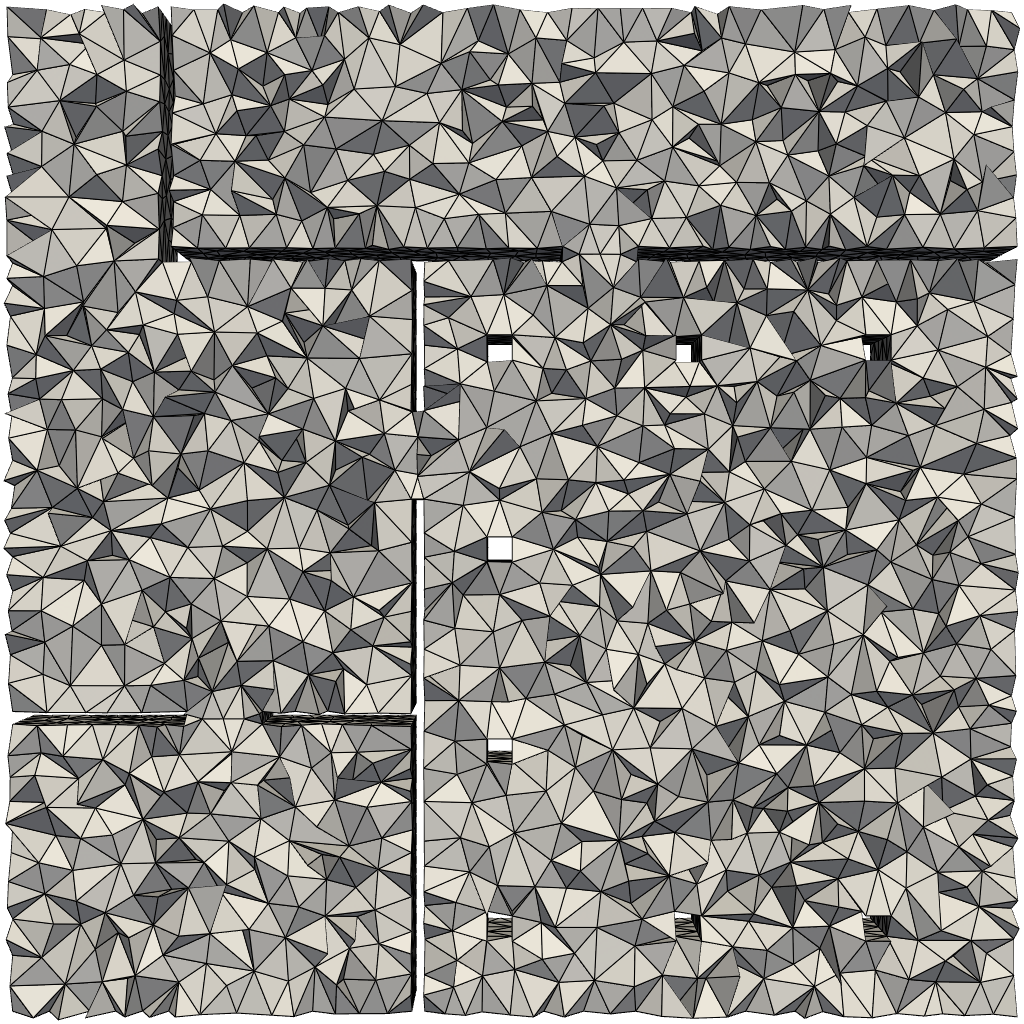}
  \caption{A view of the simple building interior that we use for a more complicated test domain. The building includes simplified architectural features, such as recessed ceilings, columns, and multiple rooms connected by shallow portals, emulating doorways.}\label{fig:simple-building}
\end{figure}

\begin{figure}
  \centering
  \includegraphics[width=0.5\linewidth]{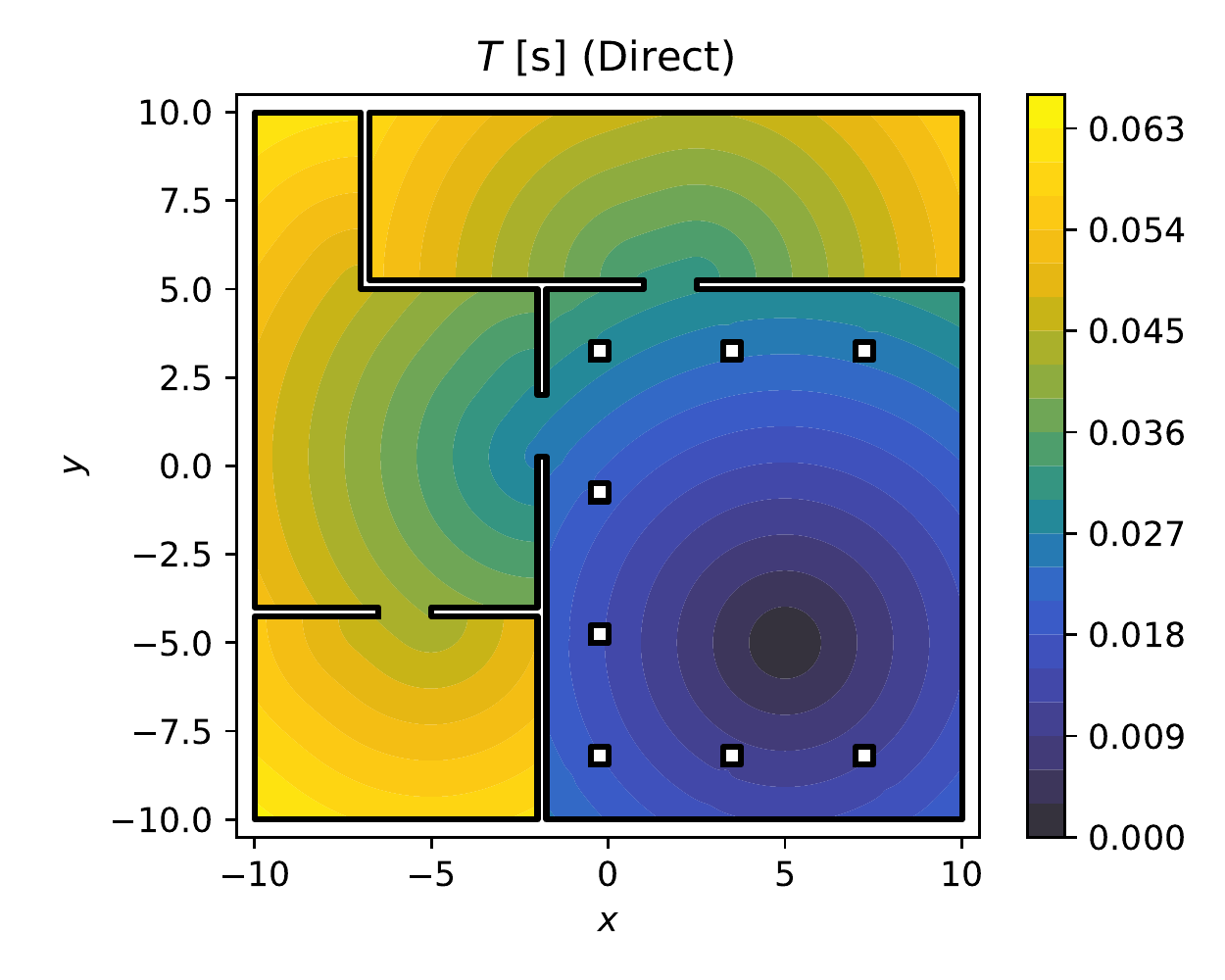}%
  \includegraphics[width=0.5\linewidth]{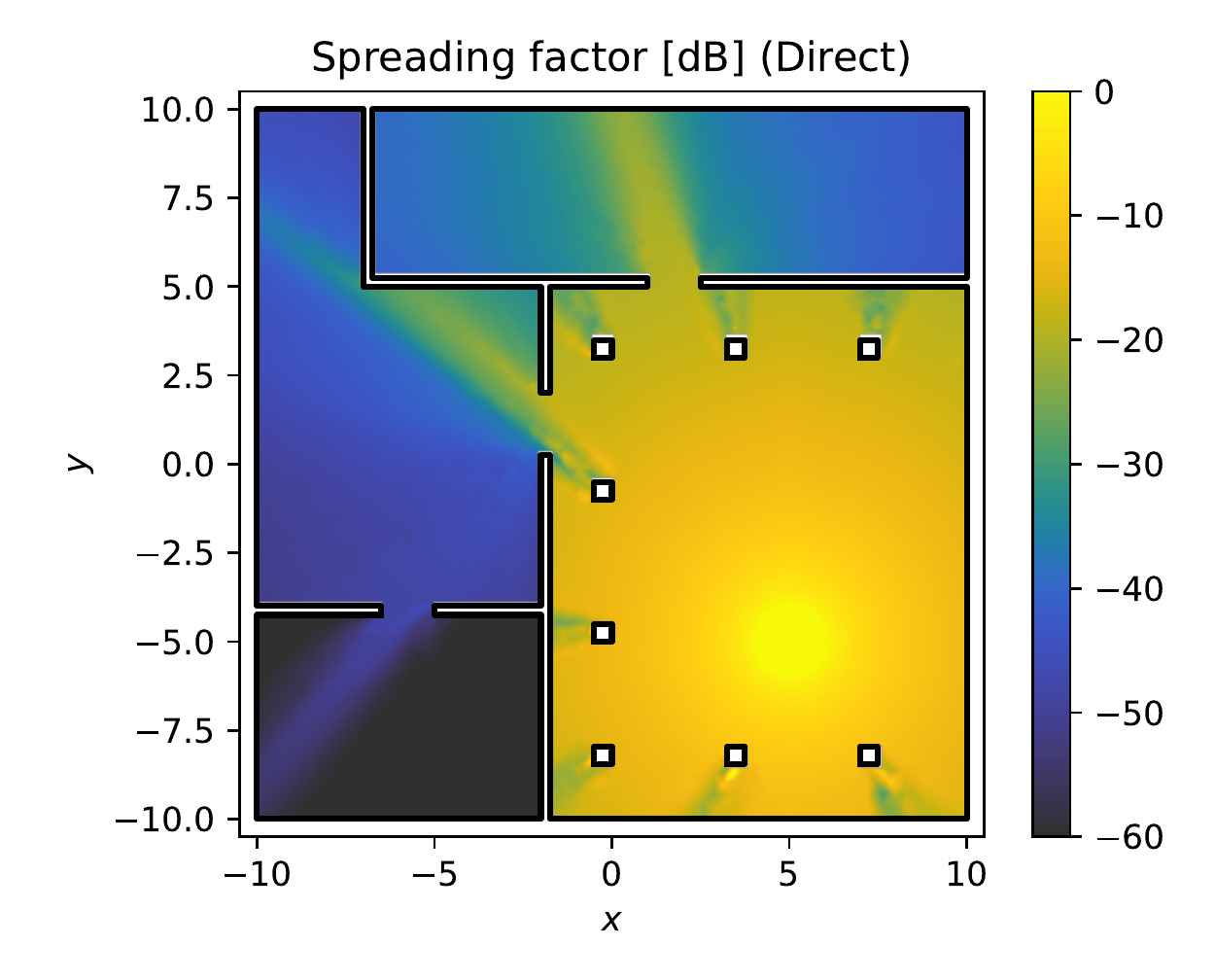}
  \includegraphics[width=0.5\linewidth]{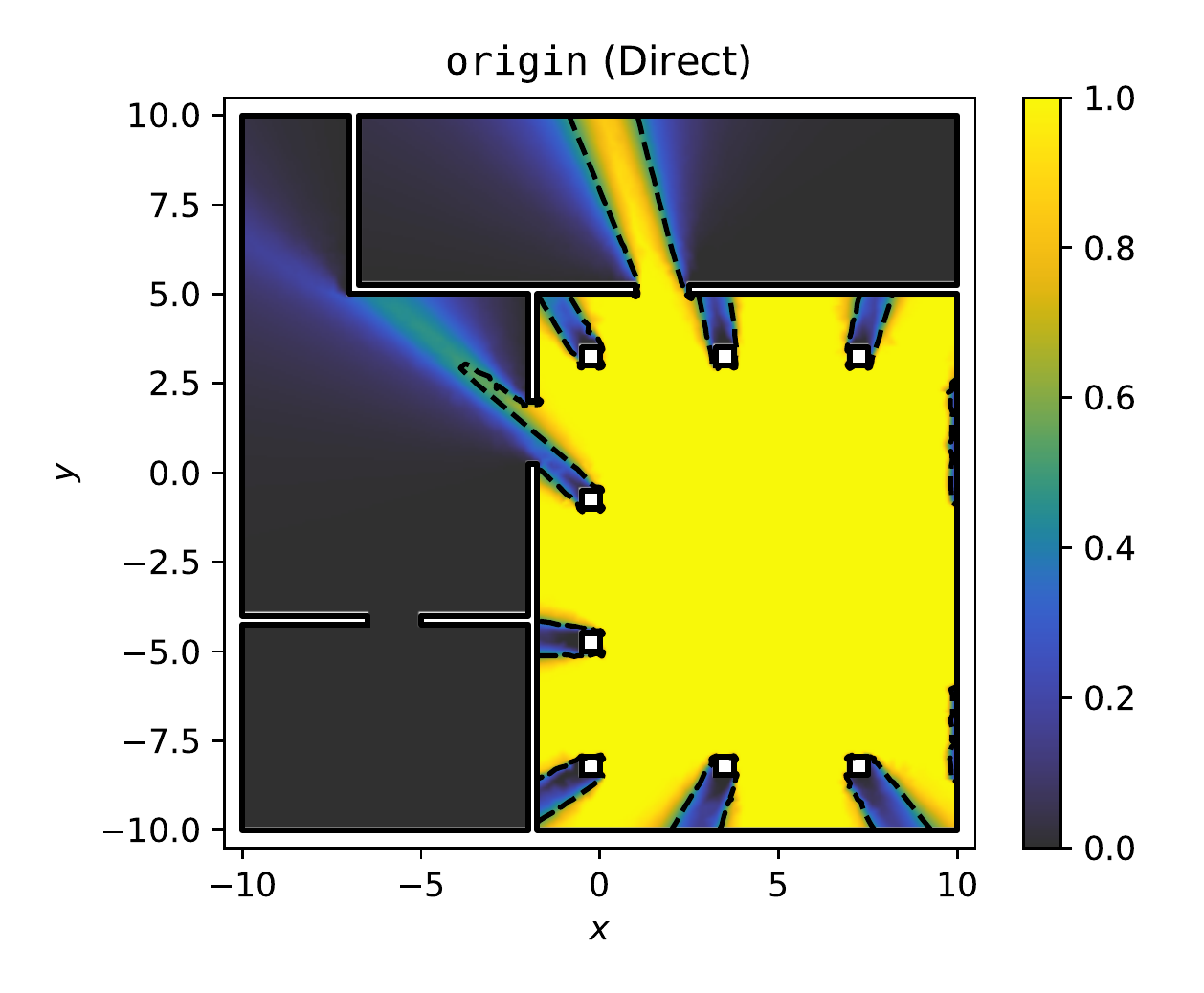}%
  \includegraphics[width=0.5\linewidth]{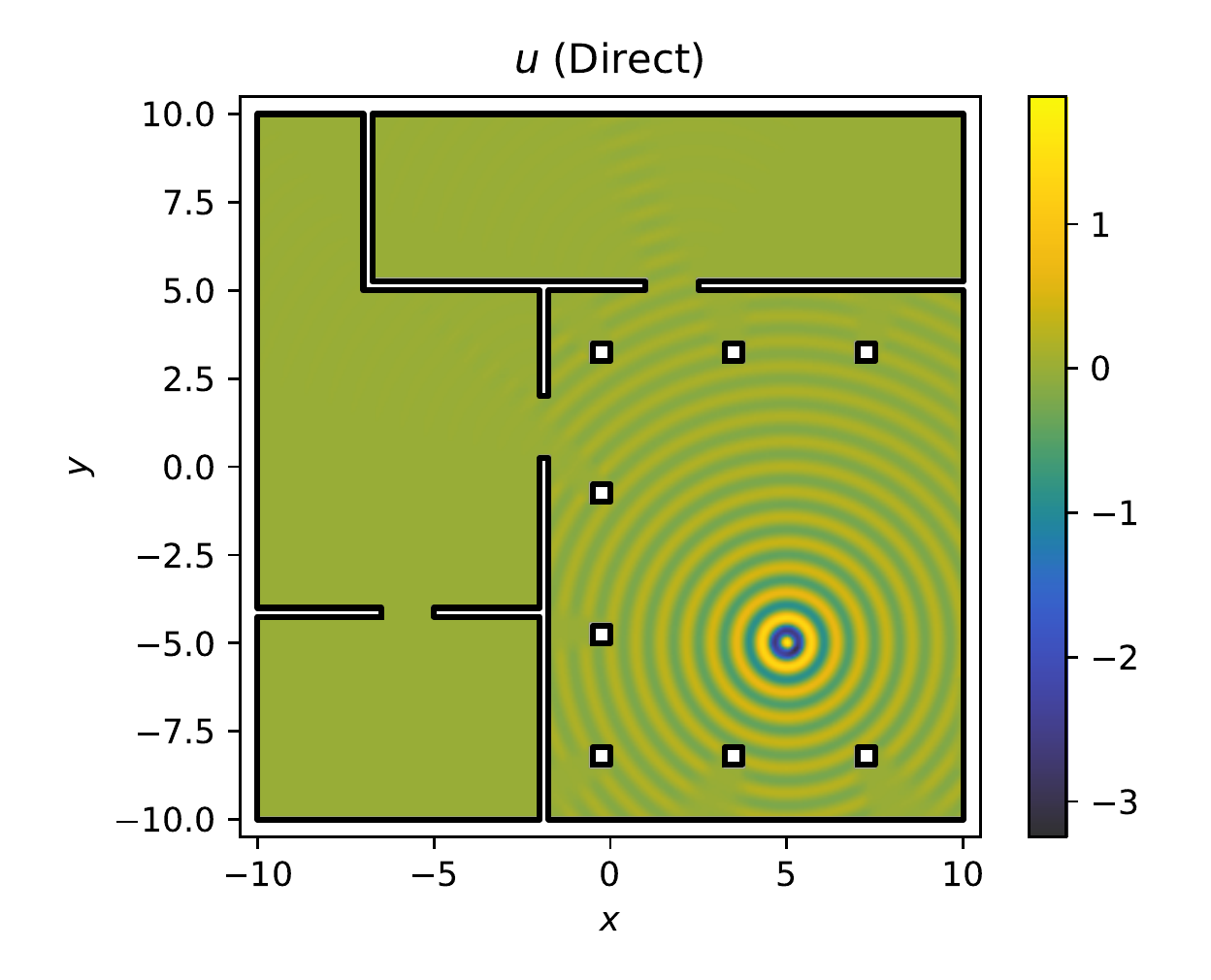}
  \caption{Horizontal slices ($z=1$) of fields corresponding to the
    first arrival for a point source placed at $\mxsrc = (5, -5, 1)$
    for the simple building geometry. \emph{Top left}: the first arrival time in seconds. \emph{Top right}: the spreading factor in decibels (dB). \emph{Lower left}: the approximate origin field, with $\partial Z$ (the $1/2$-level set) outlined. \emph{Lower right}: the heuristic solution $u$.}\label{fig:building-slices}
\end{figure}

\begin{figure}
  \centering
  \includegraphics[width=0.333\linewidth]{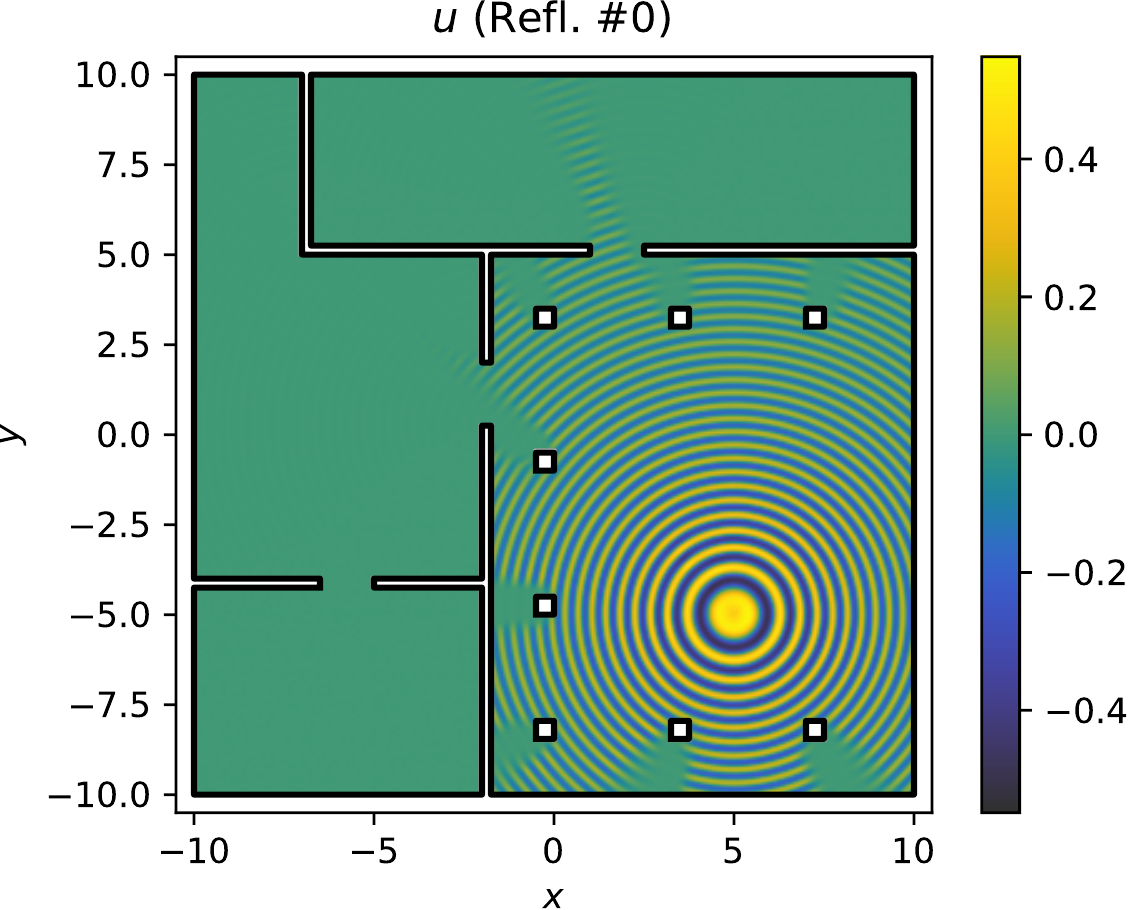}%
  \includegraphics[width=0.333\linewidth]{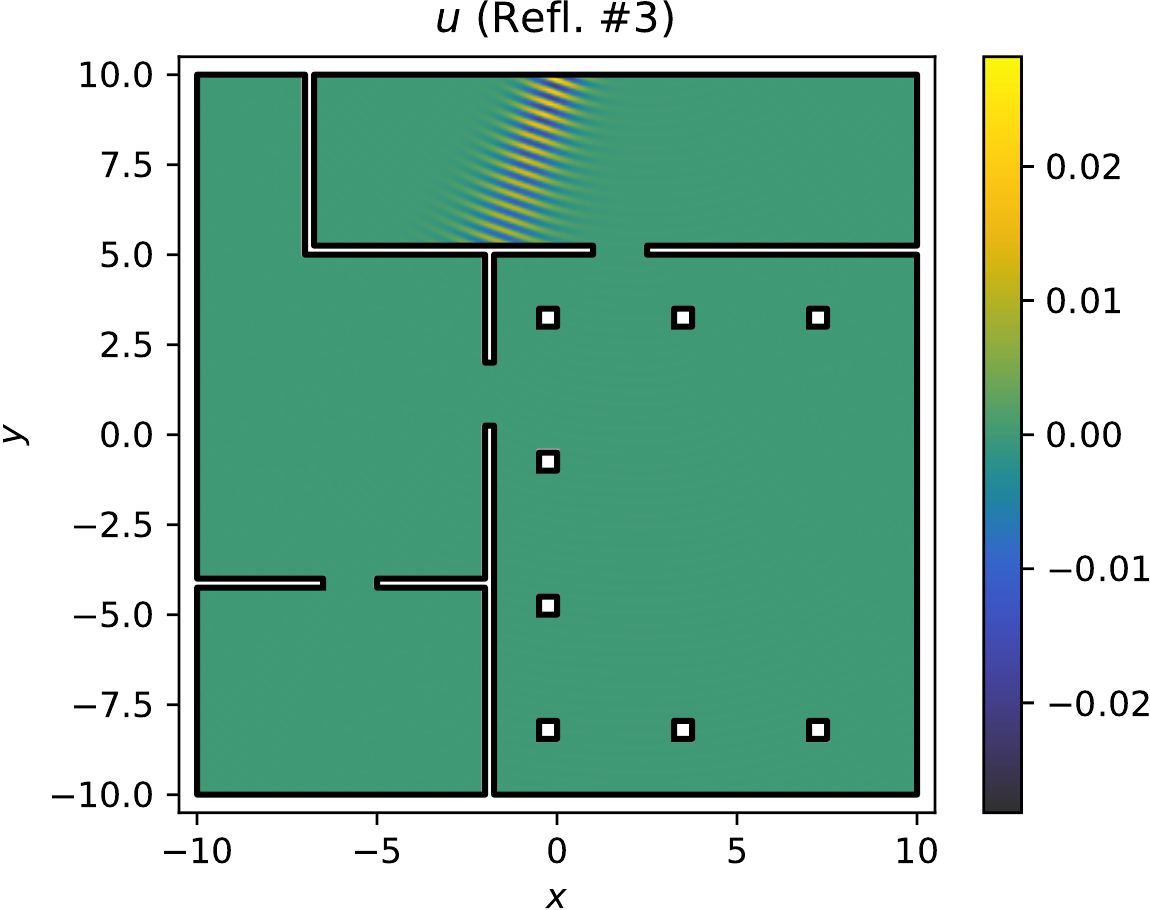}%
  \includegraphics[width=0.333\linewidth]{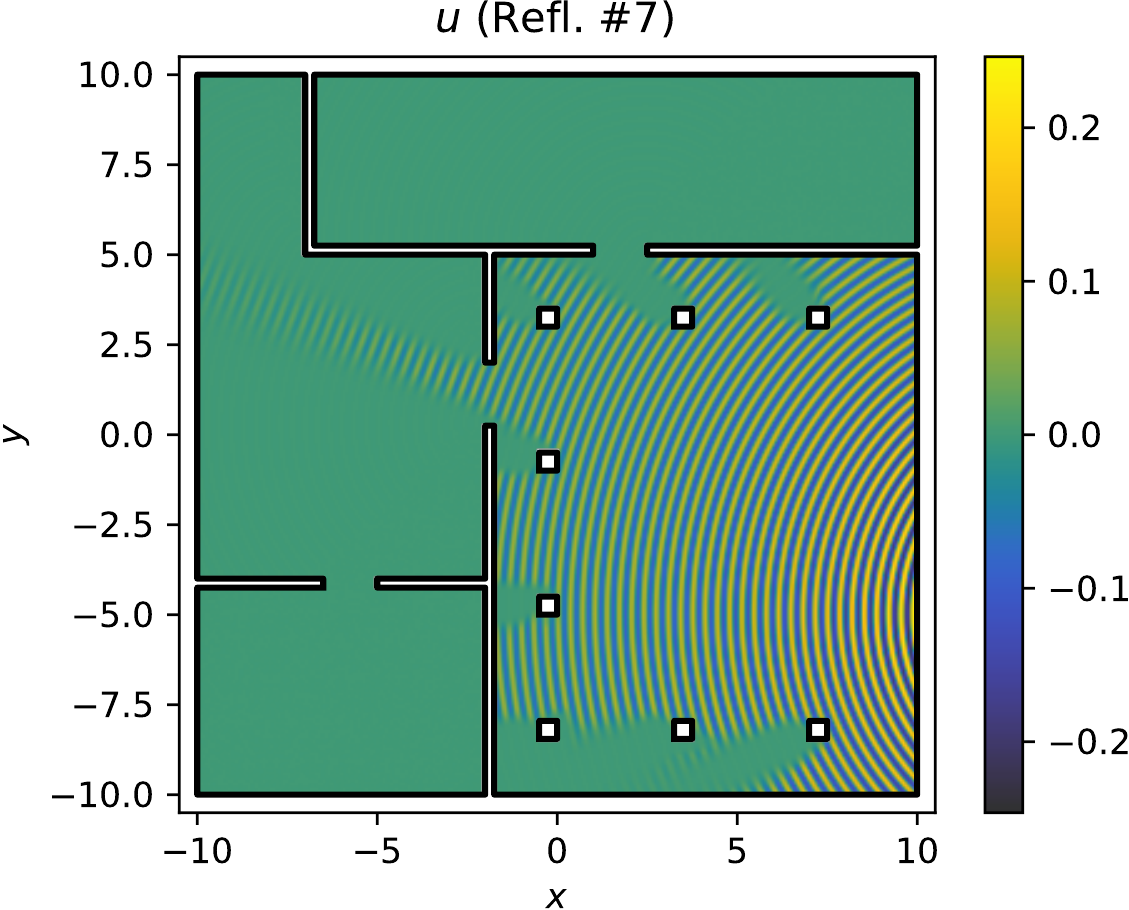} \\
  \includegraphics[width=0.333\linewidth]{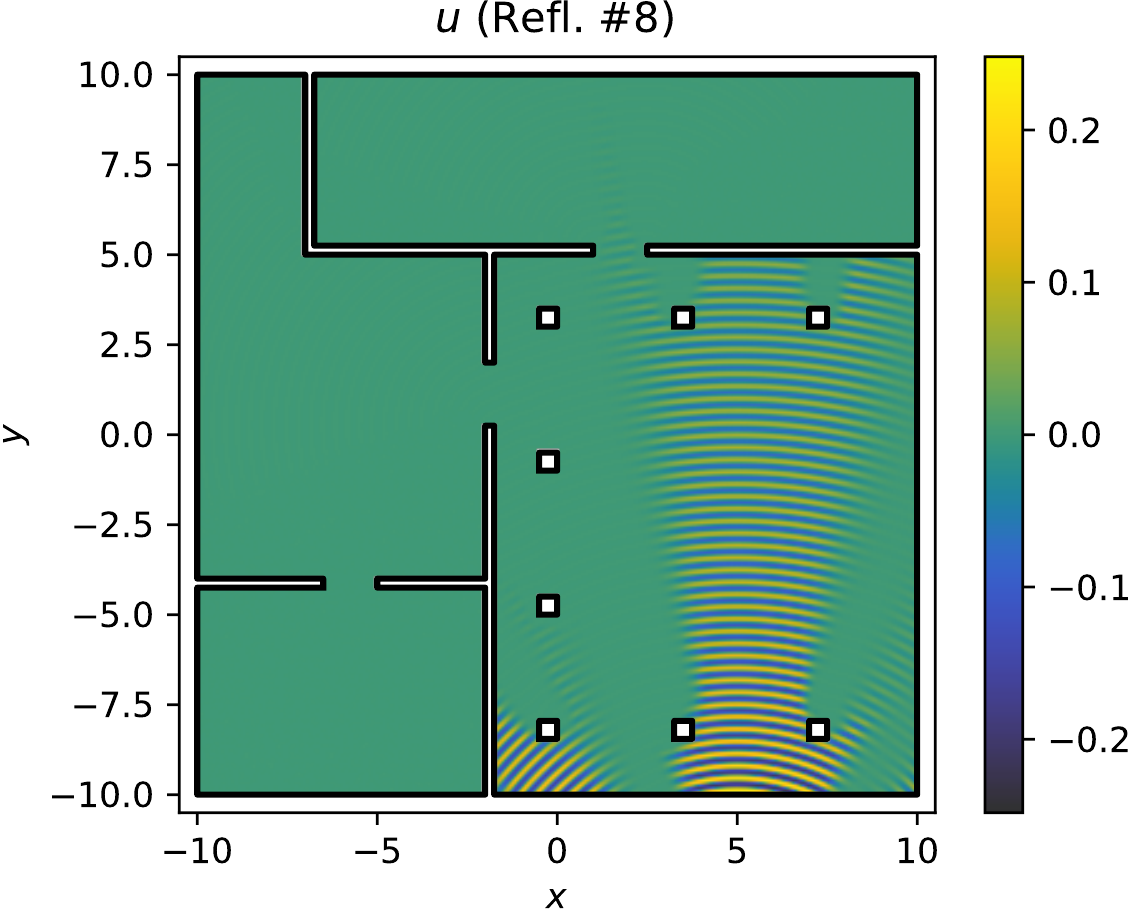}%
  \includegraphics[width=0.333\linewidth]{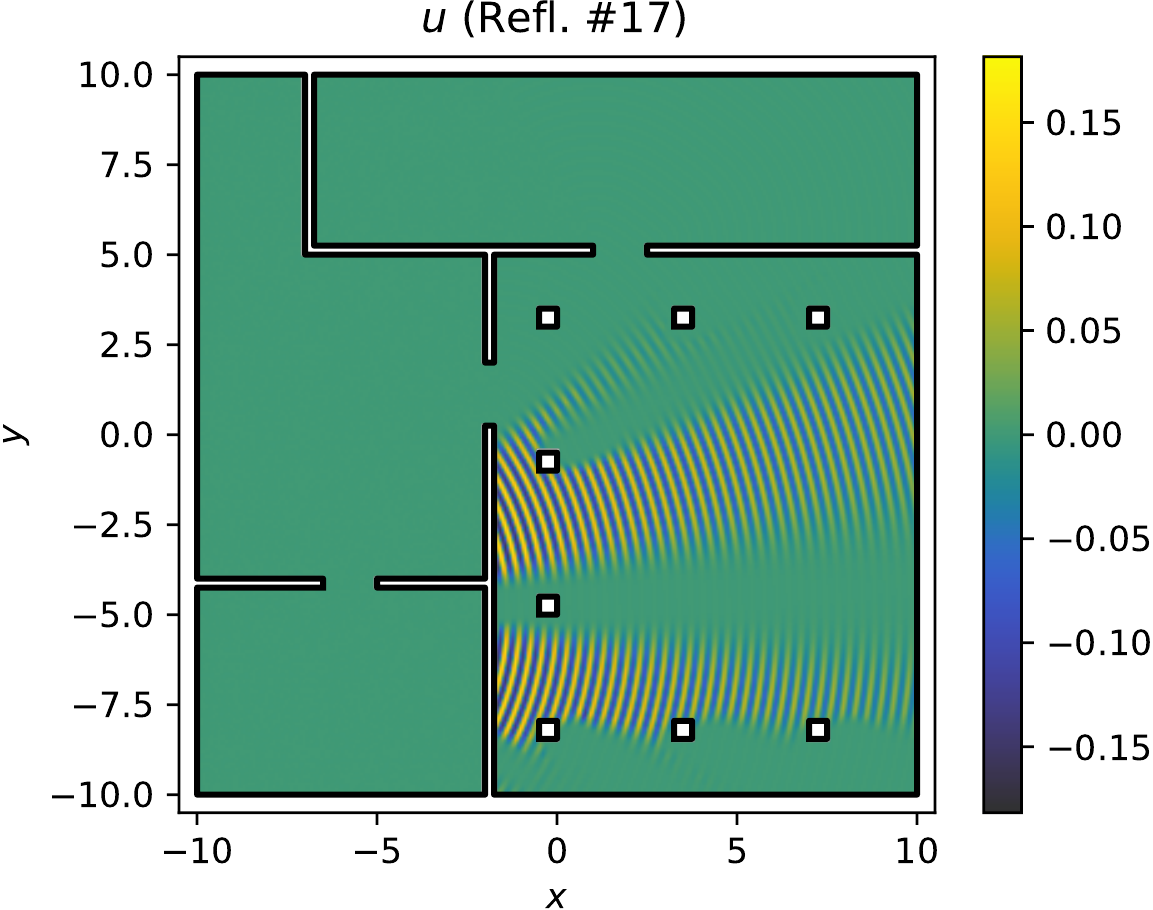}%
  \includegraphics[width=0.333\linewidth]{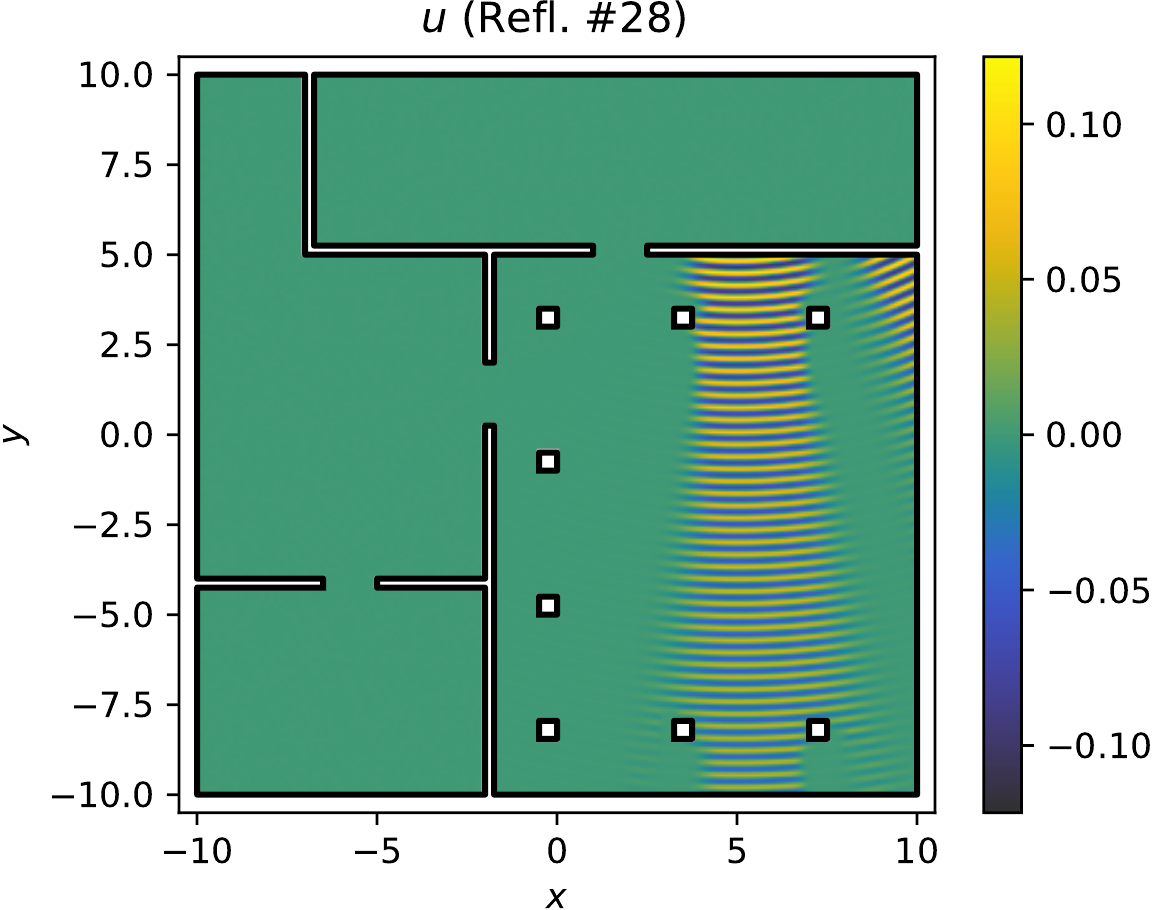} \\
  \includegraphics[width=0.333\linewidth]{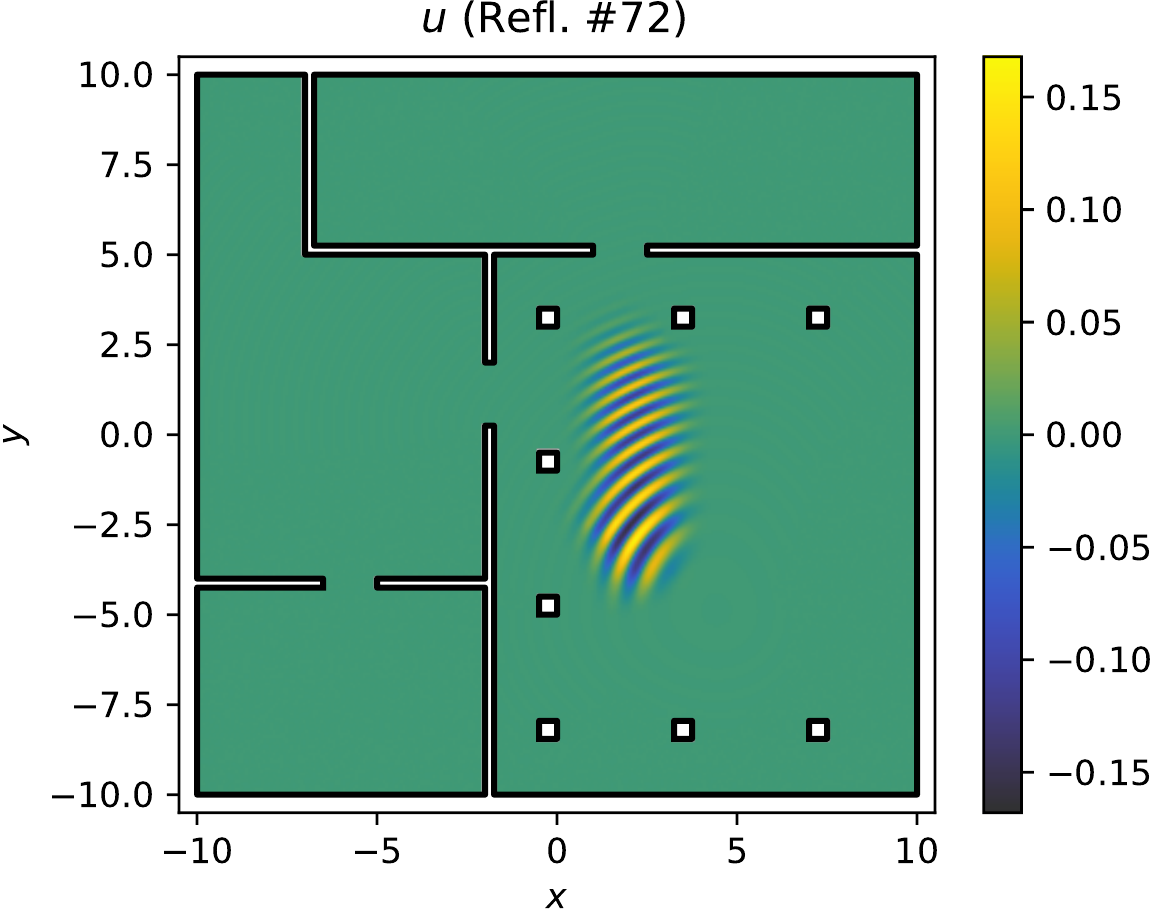}%
  \includegraphics[width=0.333\linewidth]{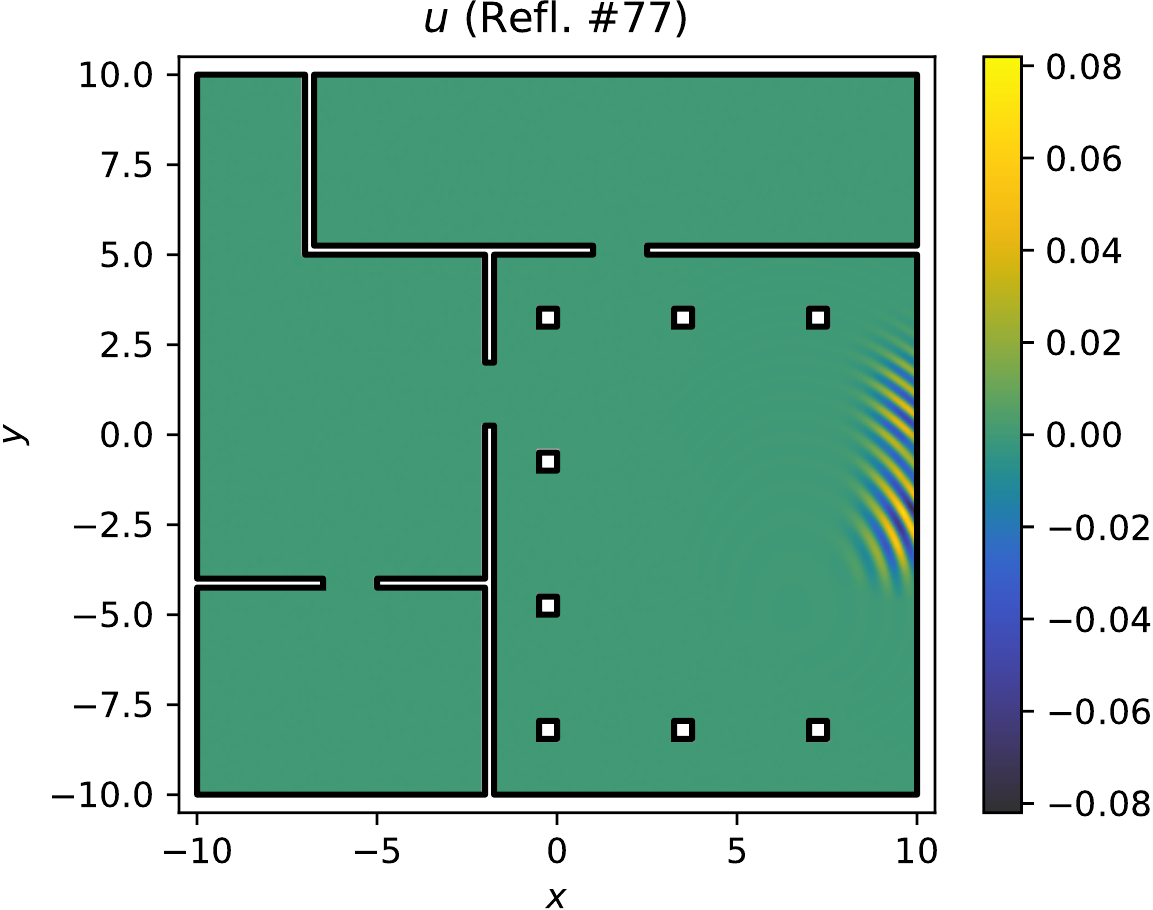}%
  \includegraphics[width=0.333\linewidth]{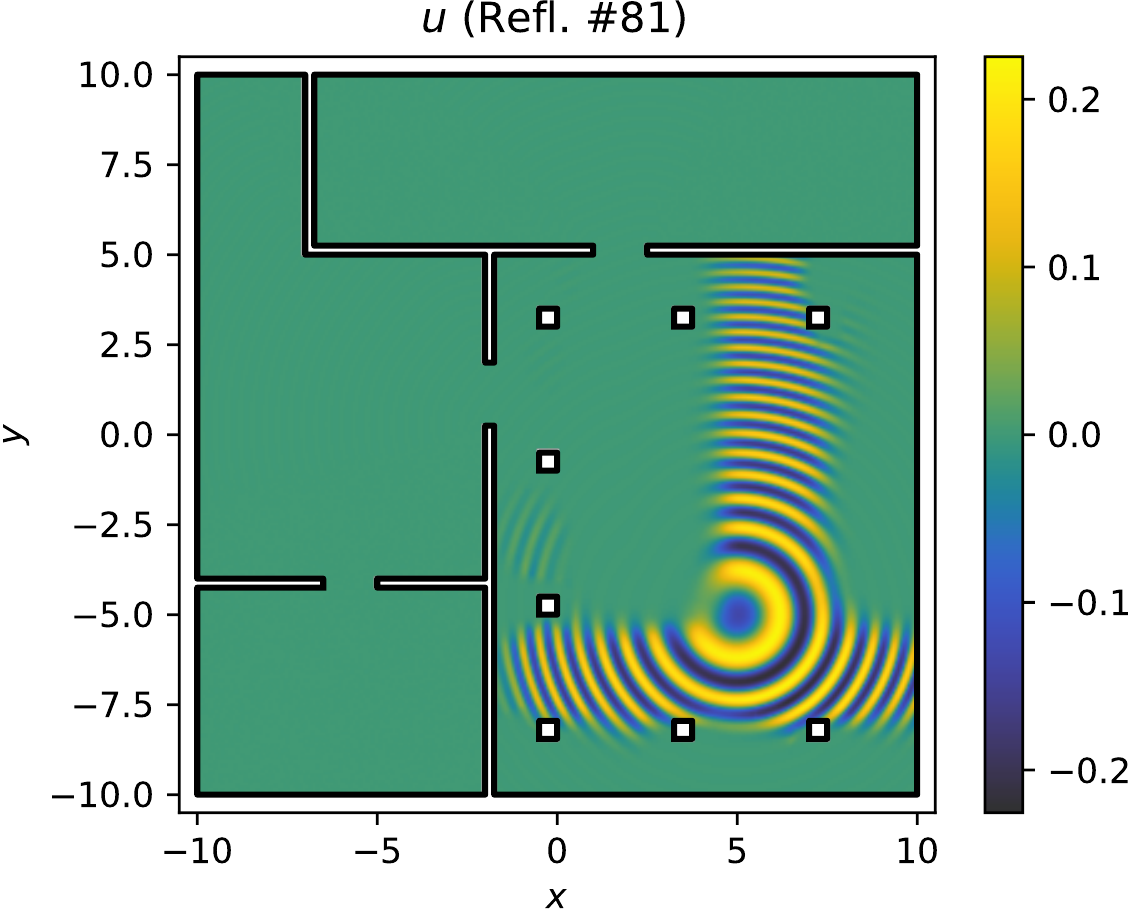}
  \caption{A variety of individual ``reflections'' propagated downwind
    from a direct field due to a point source at
    $\mxsrc = (5, -5, 1)$. Each reflecting facet in the mesh is given
    an index to distinguish among them, and the corresponding number
    is shown in the subplot title to differentiate the reflect
    fields. In these plots, the reflected fields are tapered using a
    scaled version of $\mathtt{org}(\mx)$.}\label{fig:simple-building-reflections}
\end{figure}

As an example of a more realistic test problem, we modeled a simple 3D domain with a piecewise linear boundary which has features which are typical of a simple built environment. See \Cref{fig:simple-building}. There are four rooms of differing shape and size, interconnected by door-shaped portals. Two of the rooms have recessed ceilings, and the main room additionally has features shaped like skylights along with several columns. The bounding box of the domain is $[-10, 10] \times [-10, 10] \times [0, 4.5]$, where we assume units are in meters. In this example, we fix a single maximum tetrahedron volume constraint of 0.05, resulting in $N_h = |\calV_h|$ = 10,835 vertices and $|\calC_h| =$ 51,755 tetrahedra in the mesh. This is a very coarse tetrahedron mesh.

We could imagine computing the groundtruth eikonal for this problem by repeatedly applying the tricks used in \Cref{sec:wedge} to find the analytic eikonal, but this quickly becomes too cumbersome to be practical. This illustrates one of the benefits of our approach---semi-analytic approaches based on geometric algorithms and optimization are unwieldy. For this reason, we do not attempt to directly evaluate the numerical accuracy in a more complicated domain like this. By solving the eikonal equation, we introduce some approximation error into our solution in order to use the viscosity solution as a tool to automatically find the various branches of the multivalued eikonal describing the propagating GO wavefront.

To get a sense of how the GO wavefront propagates throughout the domain, we compute the solution for a point source placed in the main room, with $\mxsrc = (5, -5, 1)$. We then build the tetrahedral spline approximating $T$ (\Cref{ssec:bernstein-bezier}) and raytrace the GO wavefront at different time points by intersecting camera rays with level sets drawn from the tetrahedral spline (\Cref{ssec:extracting-level-sets}). See \Cref{fig:wavefront-plots-top}. In this case, we do not attempt to correct the solution and improve it incrementally by applying the reinitialization algorithm of \Cref{ssec:reinit-shadow-zone}---we allow some error to accumulate as the GO wave diffracts around obstacles. Nevertheless, the wave appears quite reasonable and smooth. See \Cref{fig:wavefront-plots-top}. We also extract horizontal slices of $T$, the $A$ (in decibels), $\mathtt{org}$, and the direct field $u$ and plot them in \Cref{fig:building-slices}.

Finally, to demonstrate multiple arrivals, we restrict $T$ to each ``visible'' reflecting facet and solve a new instance of the eikonal PDE, transport the various quantities, and compute a wavefield $u$ for each reflection. In this case, since the pillars are close to the walls, they make reinitializing the eikonal in a tube near each reflecting facet difficult. As an alternative, we simply give each reflecting facet node with BCs a $\mathtt{trial}$ state and begin marching directly. We plot horizontal slices of a number of reflections from the walls, floor, and ceiling of the main room with pillars and skylights in \Cref{fig:simple-building-reflections}.

From this example, we can see that the main challenge is complicated local geometry, in the form of door portals and pillars. We include these features to give a proper sense of the outstanding difficulties with the method. There are a number of approaches outside the scope of this paper which could be taken to improve matters. For example, adaptive mesh refinement could be used to ensure that diffracting edges are spaced apart from other diffracting edges by at least, say, 5--10 edges. Adaptively refining the mesh near singular features would have the added benefit of increasing accuracy where it is needed most. We leave the exploration of this idea for a future work.

\section{A simple extension to handle nonconstant speeds of sound}

So far, we have primarily concerned ourselves with a constant speed of sound, $c \equiv 1$. In this section, we will show that it is feasible to extend to proposed approach to a nonconstant speed of sound. We will focus on giving a proof of concept and leaving its further development for future work. An extension of our solver to moderately varying nonconstant speeds of sound can be attained by replacing the triangle and tetrahedron updates with versions that account for nonconstant $c$. In this section, we describe a simple approach to making this modification based on minimizing a sequence of quadratic surrogate functions, avoiding the complications of SQP that are presented in this case.

\paragraph{A line update for a nonconstant speed of sound} Our modified updates are based on having access to a ``line update'' for nonuniform $c$. We will approximately minimize Fermat's principle by approximating the ray with a quadratic. Essentially, we solve the two-point BVP for raytracing in a variational form over a very short interval where we can expect the ray to be smooth and not vary too much. Note that this is compatible with the assumptions we make throughout the rest of this work and presents no real limitation.

Let $\mxhat$ be the update point and $\mx_0$ be the source point. We assume that $L := \|\mxhat - \mx_0\| = O(h)$. We let the basis functions for degree 2 Lagrange interpolation on the interval $[0, L]$ be:
\begin{equation}
  K_0(\sigma) = 1 - 3 \frac{\sigma}{L} + 2 \frac{\sigma^2}{L^2}, \quad K_1(\sigma) = 4 \frac{\sigma}{L} - 4 \frac{\sigma^2}{L^2}, \quad K_2(\sigma) = -\frac{\sigma}{L} + 2 \frac{\sigma^2}{L^2},
\end{equation}
whose derivatives are given by:
\begin{equation}
  K_0'(\sigma) = -3 \frac{1}{L} + 4\frac{\sigma}{L^2}, \quad K_1'(\sigma) = 4 \frac{1}{L} - 8 \frac{\sigma}{L^2}, \quad K_2'(\sigma) = -\frac{1}{L} + 4\frac{\sigma}{L^2}
\end{equation}
We let $\mx_m$ be the midpoint of the quadratic and define:
\begin{equation}
  \mphi(\sigma) = \mx_0 L_0(\sigma) + \mx_m L_1(\sigma) + \mxhat L_2(\sigma).
\end{equation}
Applying Simpson's rule to the Fermat integral and letting $\mx_m$ be the variable we will optimize over gives the cost function:
\begin{equation}
  \begin{split}
    F(\mx_m) &= T(\mx_0) + \frac{L}{6} \bigg(s(\mx_0) \norm{\mphi'(0)} + 4 s(\mphi(L/2)) \norm{\mphi'(L/2)} + s(\mxhat)\norm{\mphi'(L)}\bigg) \\
    &= T(\mx_0) + \frac{L}{6} \bigg(s(\mx_0) \norm{\mphi'(0)} + 4 s(\mx_m) + s(\mxhat) \norm{\mphi'(L)}\bigg),
  \end{split}
\end{equation}
where we have used the fact that $\mphi(L/2) = \mx_m$, $\mphi'(L/2) = (\mxhat - \mx_0)/L$, and:
\begin{equation}
  \begin{split}
    \mphi'(0) = \frac{-3\mx_0 + 4\mx_m - \mxhat}{L}, \qquad \mphi'(L) = \frac{\mx_0 - 4\mx_m + 3\mxhat}{L},
  \end{split}
\end{equation}
to simplify the first expression. Note that:
\begin{equation}
  \frac{\partial}{\partial \mx_m} \norm{\mphi'(0)} = \frac{1}{\norm{\mphi'(0)}} \frac{\partial \mphi'(0)}{\partial \mx_m}^\top \mphi'(0) = \frac{4}{L} \frac{\mphi'(0)}{\norm{\mphi'(0)}}.
\end{equation}
Likewise:
\begin{equation}
  \frac{\partial}{\partial \mx_m} \norm{\mphi'(L)} = -\frac{4}{L} \frac{\mphi'(L)}{\norm{\mphi'(L)}}.
\end{equation}
Hence, the gradient of the cost function is:
\begin{equation}
  \frac{\partial F}{\partial \mx_m} = \frac{2}{3} \bigg(s(\mx_0) \frac{\mphi'(0)}{\norm{\mphi'(0)}} + L \nabla s(\mx_m) - s(\mxhat) \frac{\mphi'(L)}{\norm{\mphi'(L)}}\bigg).
\end{equation}
To do a line update, we apply any numerical optimization method to minimize $F$ with respect to $\mx_m$. After doing the line update, if we let $\mx_m^*$ be the minimizing argument, we have $T(\mxhat) = F(\mx_m^*)$ and $\nabla T(\mxhat) = s(\mxhat) \mphi'(L; \mx_m^*)/\|\mphi'(L; \mx_m^*)\|$.

\paragraph{A tetrahedron update for a nonconstant speed of sound} Developing a tetrahedron update based on using SQP to minimize a cost function based on Simpson's rule (analogous to the updates in \cite{potter2021jet}) is unwieldy. The expressions and programming become rather complicated. As an alternative, we explore iteratively minimizing quadratic surrogate functions. This is where we use the line update presented in this section.

Let $\mxhat$ be the update point and let $\mx_0, \mx_1$, and $\mx_2$ be the $\mathtt{valid}$ nodes. We parametrize points in $\conv(\mx_0, \mx_1, \mx_2)$ as $\mx_{\mlam} = \mx_0 + \mlam_1 (\mx_1 - \mx_0) + \mlam_2 (\mx_2 - \mx_0)$, where $\mlam_1 \geq 0$, $\mlam_2 \geq 0$, and $\mlam_1 + \mlam_2 \leq 1$. We define an initial grid of Lagrange nodes with the coordinates:
\begin{equation}
  \m{\Lambda}^{(\operatorname{Lag})}_0 = \{(0, 0), (1/2, 0), (1, 0), (0, 1/2), (1/2, 1/2), (1, 1)\}.
\end{equation}
These are the uniform nodes for degree 2 Lagrange interpolation on the standard triangle (i.e., the domain of $\mlam$ above). We form a degree 2 bivariate polynomial by using a line update to evaluate $T(\mlam)$ for each $\mlam \in \m{\Lambda}^{(\operatorname{Lag})}_0$ and doing Lagrange interpolation. Call this polynomial $p_0(\mlam)$.

Our iterative optimization procedure follows. For each $n \geq 0$, we solve:
\begin{equation}
  \begin{split}
    \mbox{minimize}\qquad & p_n(\mlam)  \\
    \mbox{subject to}\qquad & \mlam_1 \geq 0, \quad \mlam_2 \geq 0, \quad \mlam_1 + \mlam_2 \leq 1.
  \end{split}
\end{equation}
This can be done analytically very easily since $p_n$ is quadratic and since the feasible set is very simple. Specifically, we first analytically compute the unconstrained optimum by minimizing the quadratic $p_n(\mlam)$ over $\mathbb{R}^2$. If the solution is inside the triangle, we're done. Otherwise, we minimize $p_n$ restricted to each edge of the triangle and take the minimizing argument, which can also be done analytically. Now, let the optimum be $\mlam^*_n$. After computing $\mlam^*_n$ we contract the Lagrange grid closer to $\mlam^*_n$ be setting, e.g.:
\begin{equation}
  \m{\Lambda}^{(\operatorname{Lag})}_{n+1} = \{(\mlam - \mlam^*_n)/\beta_n + \mlam^*_n : \mlam \in \m{\Lambda}^{(\operatorname{Lag})}_n\},
\end{equation}
where $\beta_n > 1$ ($n \geq 0$) is a sequence of scaling factors. These can be determined heuristically using a variety of methods. When $\mlam_n^*$ and $\mlam_{n-1}^*$ are sufficiently close (see \Cref{ssec:tolerances}), we terminate the iteration. In practice, we have found this algorithm to be simple and reliable, and to converge with a roughly linear rate of convergence. It may be possible to choose the scaling factors so that this obtains superlinear convergence, but we defer this topic for future study.

\paragraph{A triangle update for a nonconstant speed of sound} We recall again that triangle updates must be done to account for creeping rays and diffracted rays. These can be handled with a version of the optimization algorithm described in the foregoing paragraphs.

\begin{figure}
  \includegraphics[width=\linewidth]{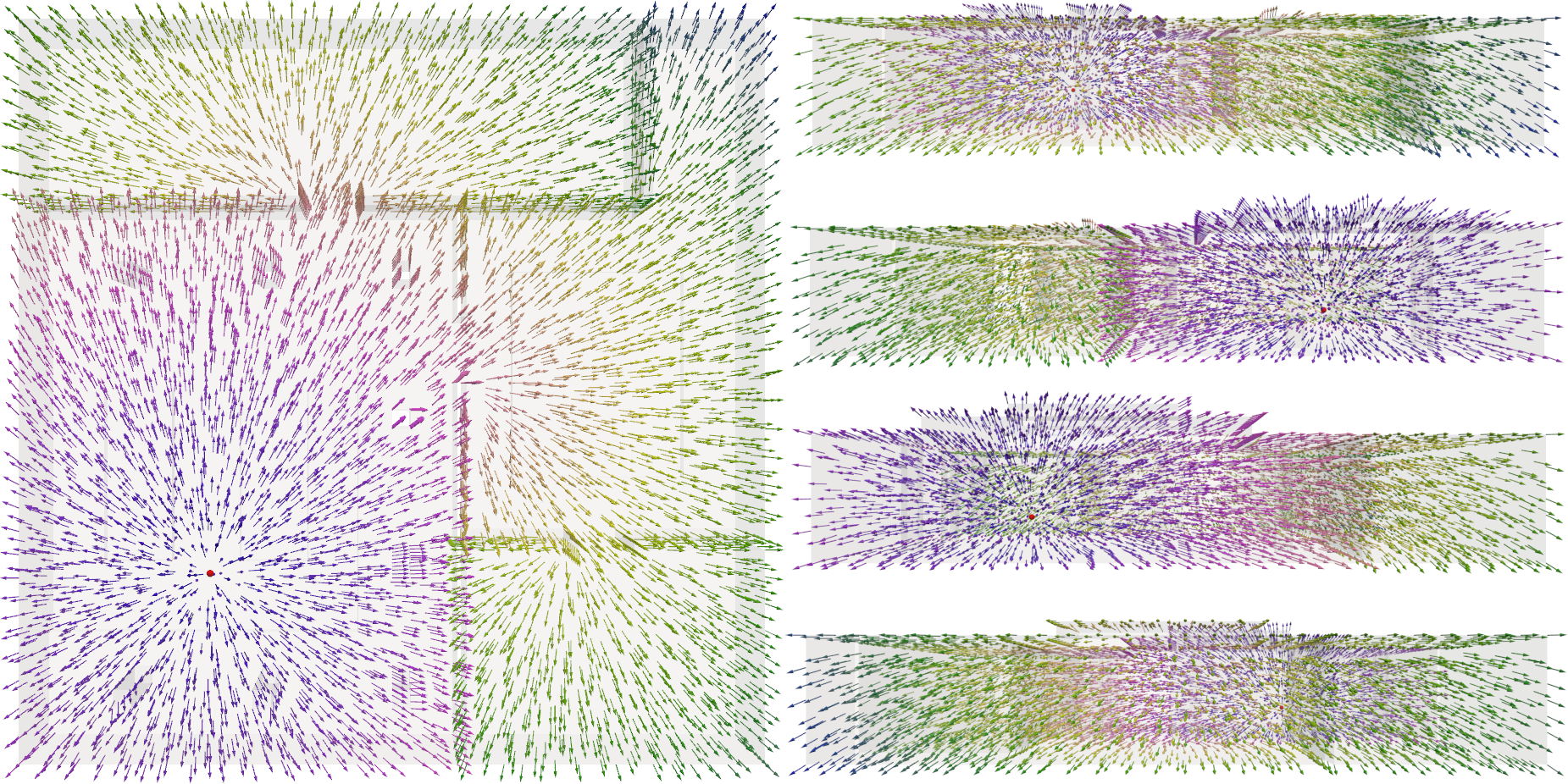}
  \caption{A quiver plot of the eikonal and its gradient for the nonconstant speed of sound given by \eqref{eq:nonconstant-c}. The vectors show the direction of the gradient while the color shows the different eikonal values. We do not show a colorbar since the scaling here is irrelevant and would make the images visually cluttered. \emph{Left}: top view. \emph{Right, top to bottom}: left view, right view, front view, and back view.}\label{fig:nonconstant-slowness}
\end{figure}

\begin{figure}
  \includegraphics[width=\linewidth]{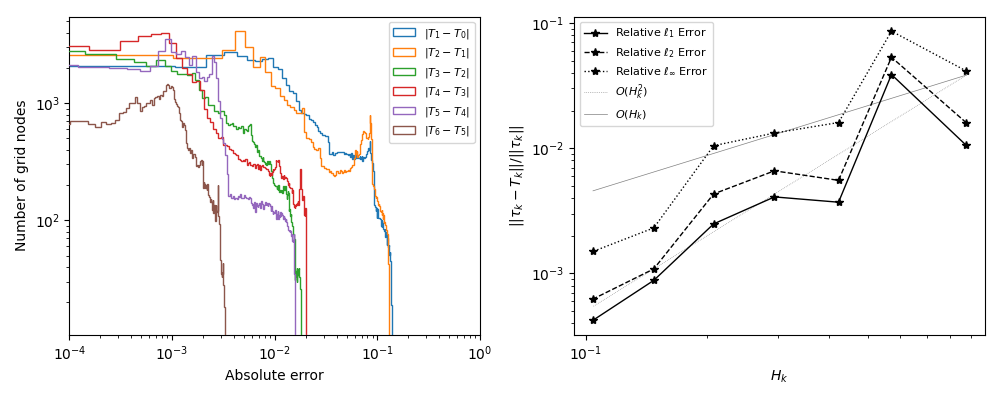}
  \caption{\emph{Left}: a loglog histogram of the computed eikonal for a varying speed of sound on the cube $\Omega = [-1, 1]^3$. The data is sampled on a $256 \times 256$ grid discretizing the horizontal plane intersecting $\Omega$ at $z = 0.5$. \emph{Right}: the relative $\ell_1$ error of the eikonal values on the $256\times 256$ grid versus the average tetrahedron mesh edge length $H_k$.}\label{fig:nonconstant-slowness-hist-and-loglog}
\end{figure}

\begin{figure}
  \includegraphics[width=\linewidth]{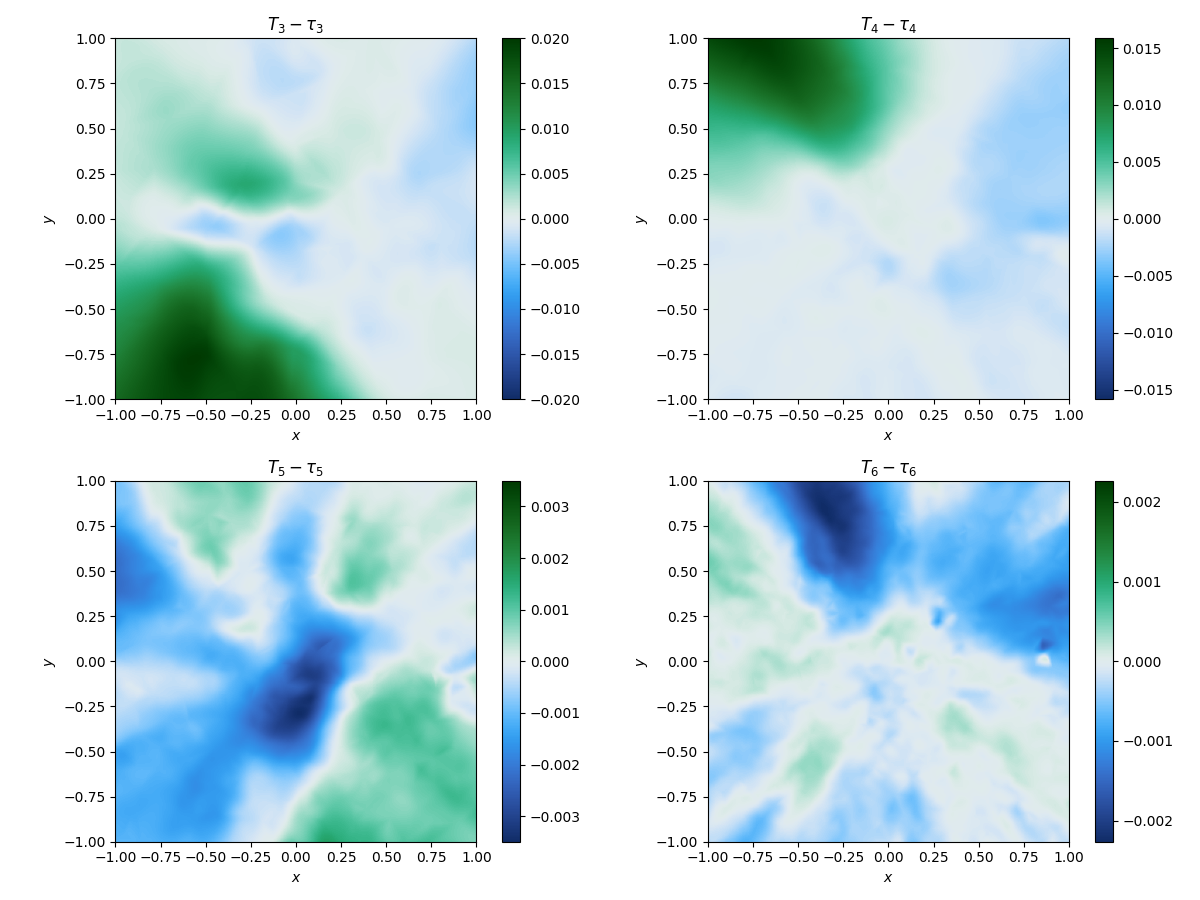}
  \caption{A slice of the pointwise error of the gridded eikonal data for the nonconstant speed of sound problem on the cube $\Omega = [-1, 1]^3$. The slice shown a uniform $256 \times 256$ grid discretizing the horizontal slice $[-1, 1] \times [-1, 1] \times \{1/2\}$ (i.e., the intersection of $\Omega$ with the horizontal plane at $z = 0.5$). The eikonal for the $k$th mesh is $T_k$, and the signed error $T_k - T_{k-1}$ is shown. In order to show the decrease in the error, we keep the color limit for each heat map the same. We can also observe that the error is smooth but somewhat random and unpredictable. This stems from the fact that tetrahedron meshes are unstructured.}\label{fig:nonconstant-slowness-pointwise-errors}
\end{figure}

\paragraph{A numerical test} To check that this method produces plausible results, we run the jet marching method with these modified updates on the building test problem from \Cref{ssec:building}. We use the speed of sound:
\begin{equation}
  c(\mx) = v_0 + \m{v}^\top \m{x}, \qquad v_0 = 1, \qquad \m{v} = (0.025, -0.025, 0.05). \label{eq:nonconstant-c}
\end{equation}
This is a rather extreme speed function for acoustics, but is a reasonable test for our method. We compute the eikonal for a single point source and show a quiver plot in \Cref{fig:nonconstant-slowness}. To get a sense of the numerical accuracy of our modified updates, we discretize the unit cube $\Omega = [-1, 1]^3$ into a sequence of tetrahedron meshes with maximum volume constraints equal to $0.1 \sqrt{8}^k$, where $k = 0, 1, 2, 3, 4, 5, 6$. The resulting meshes have 48, 138, 285, 754, 1917, 5022, and 13513 vertices, respectively, while the average edge lengths are approximately 0.875, 0.572, 0.423, 0.293, 0.207, 0.147, and 0.104. Each point source problem is solved with a factoring radius equal to 0.2. We denote the average edge length for the $k$th mesh by $H_k$. We expect our numerical method to converge at a rate roughly equal to $O(H_K^2)$ in the $\ell_1$ and $\ell_2$ norms.

To evaluate the accuracy of our solution, we first march the numerical solution on the tetrahedron mesh to get the eikonal and its gradient at each vertex of the mesh. We then build a degree 3 tetrahedral cubic \Bezier{} spline interpolating this data (as in \Cref{ssec:extracting-level-sets}). Using this spline, we evaluate the solution on a $256 \times 256$ uniform grid discretizing the horizontal plane $[-1, 1] \times [-1, 1] \times \{0.5\}$ (i.e., the $xy$-plane intersected with the cube at height $z = 0.5$). Let $T_k$ be the grid of values obtained this way for the $k$th mesh, and let $\tau_k$ denote the corresponding analytical eikonal values, which are given by:
\begin{equation}
  \tau(\mx) = \frac{1}{\norm{\m{v}}} \cosh^{-1}\left(1 + \frac{1}{2} s_0 s(\mx) \norm{\m{v}}^2 \norm{\m{x}}^2 \right).
\end{equation}
We plot the pointwise errors $T_k - \tau_k$ for $k = 3, 4, 5, 6$ in \Cref{fig:nonconstant-slowness-pointwise-errors}. We also plot histograms of the pointwise errors (where our samples correspond to the $256^2$ grid nodes) and the relative $\ell_1$, $\ell_2$, and $\ell_\infty$ errors of this gridded data with respect to $H_k$ in \Cref{fig:nonconstant-slowness-hist-and-loglog}. We emphasize that the staircasing behavior in the relative $\ell_1$ error stems from the fact that the tetrahedron meshes are of varying quality, which has a significant effect on the interpolation error. This can be seen in the pointwise error plots, which are akin to random fields---they do not exhibit any sort of pattern as the mesh refines. The included histograms show the distribution of the absolute pointwise errors, which do not decrease in a predictable fashion---again, a function of the lack of structure in the tetrahedron meshes. From our preliminary numerical results, we speculate that the error of our method in the relative $\ell_1$ and $\ell_2$ norms is $O(H_k^2)$, while for the relative $\ell_\infty$ norm, it is somewhere between $O(H_k)$ and $O(H_k^2)$.

\section{Conclusion}

In this work, we present:
\begin{enumerate}
\item A jet marching method which can be used to compute the eikonal and its gradient with second order accuracy in a complicated domain for a constant speed of sound on a tetrahedron mesh.
\item We show how to average over neighboring cells  in this context to compute the Hessian of the eikonal. Results suggest it is possible to do so with first order accuracy.
\item We introduce a new approach to handling causality for a marching method for solving the eikonal equation on an unstructured mesh by explicitly minimizing over the $\mathtt{valid}$ front. To minimize the cost of doing so and to avoid complicating the local search for updates, we carefully characterize what makes an update physical and introduce a novel method of caching previous updates.
\item We introduce algorithms for accurately reinitializing the eikonal in the shadow zone to accurately compute the viscosity solution in the presence of diffraction. This provides an alternative to local factoring.
\item Auxiliary algorithms for marching derived quantities, such as a function $\mathtt{org}$ whose $1/2$-level set can approximate the shadow and reflection GO boundaries, the amplitude prefactor $A$, and the $\mtin$ and $\mtout$ fields needed to compute the UTD coefficient needed for the ray-dependent BCs for the edge-diffracted amplitude.
\item Numerical experiments demonstrating the accuracy of our method in computing $T$, $\grad T$, and $\hess T$.
\item Numerical experiments showing the feasibility of using our method to compute acoustic fields of rays in a simplified architectural model with some challenging features.
\item A straightforward demonstration of how to extend the method to cope with nonconstant speeds of sound.
\end{enumerate}
The primary advantage of our approach is that it allows one to simultaneously trace fields of rays, avoiding issues with ray spreading and interpolation in raytracing. For this reason, it can be thought of a semi-Lagrangian form of \emph{beam tracing}~\cite{heckbert1984beam}. An explicit goal is to avoid global queries, such as global raytracing, visibility tests, or using the method of images to handle reflections.

The results presented here can be viewed as part of a long term project investigating the feasibility of computing a high-frequency geometric optic (WKB) approximation of the Helmholtz Green's function in 3D in a domain with a complicated boundary. We have considered the simplified case of a constant speed of sound and a piecewise linear boundary. The purpose of this paper is to present a new algorithmic toolbox for approaching this problem which can be expanded on in the future to handle more complicated problems. Our belief is that the algorithms presented here have the essential features to accurately compute the early arrivals of an acoustic wave propagating in a 3D environment. This is useful, since the early arrivals contain a significant amount of the ``perceptual content'' of sound~\cite{kuttruff2016room}. There are several principle extensions remaining: 1) handling nonconstant speeds of sound, 2) handling freespace caustics, and 3) handling curved surfaces. These extensions are sufficiently complicated to warrant a paper treatment each, at least. The amplitude equation deserves special attention, since new types of distinct diffraction phenomena arise in each case.

A fundamental limitation of our approach is that we do not attempt to refine the mesh around the shadow boundary, the reflection boundary, or the silhouette lines. Not only is doing so a delicate proposition, it is straightforward to construct examples where allowing surfaces with highly faceted structures leads to GO boundaries which are extremely complicated---clearly, as these boundaries undergo repeated reflection and diffraction, they become too complex to compute. At the same time, we make a modeling assumption when we employ UTD to model edge-diffraction---that is, UTD models diffraction from a \emph{semi-infinite} wedge, but we violate this when we apply this to diffraction from a wedge of finite extent. It is unclear what the consequences of doing so are. These two limitations are related. If we employ a high-frequency diffraction theory such as UTD, we essentially require that the boundary of our surface can be decomposed into a small number of acoustically large smooth facets and singular features. Whether these features are acoustically large depends on the wavenumber. In the future, a more careful study of the discrepancies that exist between the various approaches to modeling wave propagation would be useful in terms of being able to decide the appropriateness of using a particular approach to solving a problem.

We also briefly remark on whether ``numerical geometric acoustics'' should be considered a time or frequency domain algorithm. The situation is a bit circuitous. We initially start by thinking of computing the Green's function for the wave equation. Making the time-harmonic assumption leads us to the Helmholtz equation, the PDE which the Fourier coefficients of the solution of the Helmholtz equation must satisfy. Once we make the WKB ansatz to approximate the solutions of the Helmholtz equation, we obtain the eikonal and amplitude PDEs, describing the behavior of a propagating GO wave. Considering that we can reinitializing these PDEs to propagate secondary waves, it becomes clear that we are dealing with time-domain quantities. Indeed, the multiple branches of the ``multipath eikonal'' produced this way essentially give the Green's function for the wave equation, leading us back to where we started. From these considerations, it is perhaps best to think of our approach as a hybrid.

The code used to carry out the numerical results will be released as a part of the \texttt{jmm} library, which is available on GitHub under the open source Apache 2.0 license:
\begin{center}
  \texttt{https://github.com/sampotter/jmm}
\end{center}

\section{Acknowledgments}

This work was supported in part by the NSF CAREER grant DMS1554907 and MTECH grant No.\ 6205.

\bibliographystyle{plain}
\bibliography{jmm3d}

\appendix

\section{Evaluating the UTD edge-diffraction coefficient}\label{sec:UTD-coef}

We now explain the edge-diffraction coefficient and how to evaluate it, summarizing standard material on UTD~\cite{mcnamara1990introduction}.

The edge diffraction coefficient is parametrized by the local wedge and ray geometry at the point of diffraction. In the UTD literature, we label the two faces of the wedge the \emph{o-face} and the \emph{n-face}. See \Cref{fig:edge-centered-coords} for this setup. For each face, we have an outward facing unit surface normal, denoted $\mn_o$ and $\mn_n$, respectively. We let $\mt_e$ denote the unit tangent vector of the edge. By convention, the $o$-face is the face of the wedge for which the face
normal crossed with the edge tangent vector points into the face.

The angle the incident ray makes with the edge ($\beta_{\In}$) and the angle the cone makes with the edge ($\beta_{\Out}$) are the same. Denote this angle $\beta$. If we let $\mtin$ and $\mtout$ denote the unit tangent vectors of an ingoing and outgoing ray, respectively, then we have:
\begin{equation}
  \beta = \cos^{-1}(-\mt_e\cdot\mtin) = \cos^{-1}(\mt_e\cdot\mtout).
\end{equation}
If we let $\mt_e$ denote the $z$-axis of a spherical coordinate system whose origin is the point of diffraction, then $\beta$ gives the elevation angles of the ingoing and outgoing rays. For the azimuthal angle, we let $\mt_o$ be the $x$-axis of the spherical coordinate system, and let $\phi$ be the azimuth so that $\phi = 0$ corresponds to the $o$-face and $n\pi$ corresponds to the $n$-face. Let $\mtin^\perp$ and $\mtout^\perp$ be the projection of $\mtin$ and $\mtout$ into the $xy$-plane of this \emph{edge-centered coordinate system}:
\begin{equation}
  \mtin^\perp = \frac{\mtin - (\mt_e\cdot\mtin)\mt_e}{\norm{\mtin - (\mt_e\cdot\mtin)\mt_e}}, \qquad \mtout^\perp = \frac{\mtout - (\mt_e\cdot\mtout)\mt_e}{\norm{\mtout - (\mt_e\cdot\mtout)\mt_e}}.
\end{equation}
Then, the azimuth of the ingoing and outgoing rays can be computed from:
\begin{equation}
  \phiin = \tan^{-1}\parens{\frac{\mn_o\cdot\mtin^\perp}{\mt_o\cdot\mtin^\perp}}, \qquad \phiout = \tan^{-1}\parens{\frac{\mn_o\cdot\mtout^\perp}{\mt_o\cdot\mtout^\perp}}, \qquad \phiin, \phiout \in [0, n\pi].
\end{equation}
These angles allow us to parametrize the ingoing and outgoing rays in the edge-centered coordinate system.

We now introduce some auxiliary quantities used to define $D$. We start with the length parameter $L$, given by:
\begin{equation}
  L(\mx) = \frac{\rhodiff(\mx)\big(\rho_e(\mx) + \rhodiff(\mx)\big)\rho_1(\mx) \rho_2(\mx)}{\rho_e(\mx)\big(\rho_1(\mx) + \rhodiff(\mx)\big)\big(\rho_2(\mx) + \rhodiff(\mx)\big)} \sin(\beta(\mx))^2.
\end{equation}
As before, $\rhodiff(\mx)$ is the distance along the ray from the point of the diffraction to the observation point. The parameter $\rho_e(\mx) = 1/\kappa_{\mq_e}(\mx_e)$ is the sectional radius of curvature determined by the unit vector:
\begin{equation}
  \mq_e = \frac{\mt_e - \big(\mtin(\mx_e)\cdot\mt_e\big)\mtin(\mx_e)}{\norm{\mt_e - \big(\mtin(\mx_e)\cdot\mt_e\big)\mtin(\mx_e)}},
\end{equation}
where $\mn$ is the surface normal of the directly visible face. The parameters $\rho_1(\mx) = 1/\kappa_1(\mx_e)$ and $\rho_2(\mx) = 1/\kappa_2(\mx_e)$ are the radii of curvature associated with the principal curvatures of the GO wave at the point of diffraction.

Note that some of these parameters are nominally functions of $\mx$ (the field point), but are actually defined in terms of $\mx_e$ (the point of diffraction). Writing $\mx$ as a function of $\mx_e$ is simple if we adopt a purely Lagrangian perspective, but it is less straightforward for a semi-Lagrangian (or even Eulerian) solver. We address this issue in \Cref{sec:amplitude}.
In order to define the edge diffraction coefficient, we also need to define the weight function:
\begin{equation}
  a^{\pm}(\beta) = 2\cos\parens{\frac{2\pi{}nN^{\pm}(\beta) - \beta}{2}}^2, \quad \mbox{where} \quad N^{\pm}(\beta) = \Arg\min_{N\in\mathbb{Z}} \abs{\beta \pm \pi - 2\pi{}nN}.,
\end{equation}
and the Kouyoumjian transition function:
\begin{equation}
  F(x) = 2i\sqrt{x}e^{ix}\int_{\sqrt{x}}^\infty e^{-iy^2}dy = 2i\sqrt{x}e^{ix} F_{-}(\sqrt{x}),
\end{equation}
where $F_- = \int_x^\infty e^{-iy^2}dy$ is the modified negative Fresnel integral.

\begin{figure}
  \centering
  \includegraphics[width=0.9\linewidth]{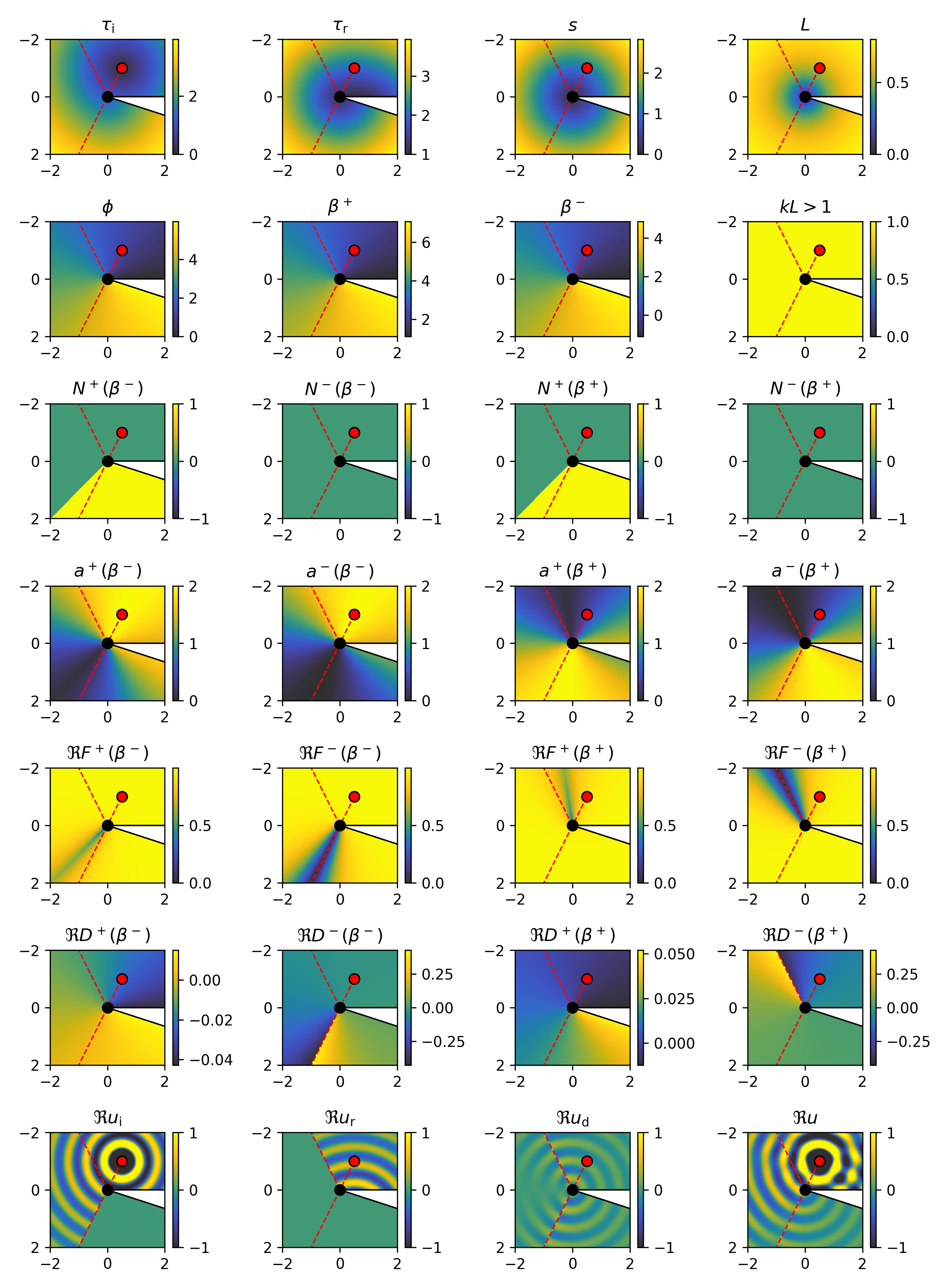}
  \caption{Plots of the various parameters involved in evaluating the edge diffraction coefficient. Note that the terms $D_1, D_2, D_3$, and $D_4$ are each associated with a different part of the edge-diffracted field.}\label{fig:utd}
\end{figure}

Finally, the edge diffraction coefficient is given by:
\begin{equation}
  D = D_1 + D_2 + R \cdot (D_3 + D_4),
\end{equation}
where $R$ is the reflection coefficient. Letting $\beta^\pm = \phi_{\Out} \pm \phi_{\In}$ and using the foregoing definitions, the terms $D_1, D_2, D_3$, and $D_4$ making up the edge-diffraction coefficient $D$ are given by (omitting some arguments):
\begin{equation}
  \begin{split}
    D_1 &= \frac{-e^{i\pi/4}}{2n\sqrt{2\pi{}k}\sin(\beta)} \cot\parens{\frac{\pi + \beta^-}{2n}} F\big(kLa^+(\beta^-)\big), \\
    D_2 &= \frac{-e^{i\pi/4}}{2n\sqrt{2\pi{}k}\sin(\beta)} \cot\parens{\frac{\pi - \beta^-}{2n}} F\big(kLa^-(\beta^-)\big), \\
    D_3 &= \frac{-e^{i\pi/4}}{2n\sqrt{2\pi{}k}\sin(\beta)} \cot\parens{\frac{\pi + \beta^+}{2n}} F\big(kLa^+(\beta^+)\big), \\
    D_4 &= \frac{-e^{i\pi/4}}{2n\sqrt{2\pi{}k}\sin(\beta)} \cot\parens{\frac{\pi - \beta^+}{2n}} F\big(kLa^-(\beta^+)\big). \\
  \end{split}
\end{equation}
Recall that $k = \omega/c$ is the wavenumber. These coefficients are attached to different combinations of faces and shadow boundaries. See \Cref{fig:utd} for a depiction of each of the parameters participating in $D$'s definition.

\section{Minima common to multiple adjacent updates}\label{sec:common-minima}

\begin{figure}
  \centering
  \includegraphics[width=\linewidth]{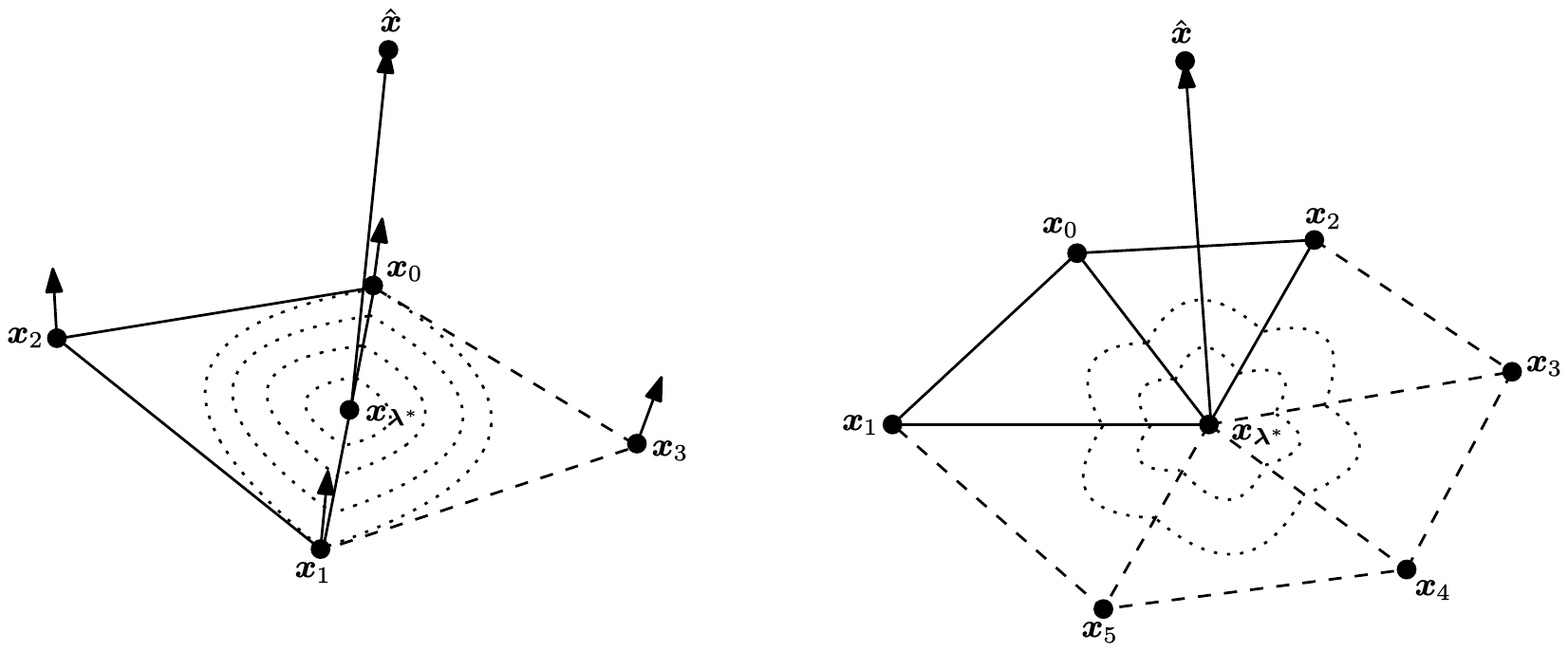}
  \caption{Special cases which need to be handled using an update cache. Note: in keeping with our notation, in both figures, $\mx_0$ is the newly $\mathtt{valid}$ node. Triangles with a solid outline are updates for the current update fan, while dashed triangles are not. \emph{Left}: an \emph{edge minimzer}. The point $\mx_{\mlam^*}$ when viewed as the minimizer for each incident tetrahedron update is a boundary minimizer in the interior of the edge $[\mx_0, \mx_1]$ with nonzero Lagrange multipliers. To determine that it is a physical update, we need to take both updates together---however, these two updates may not occur in the same phase of updates, making it necessary to cache updates for later comparison. \emph{Right}: a \emph{vertex minimizer}, for which it is necessary to track a ring of tetrahedron updates surrounding $\mx_{\mlam^*}$.}\label{fig:minimizers}
\end{figure}

Let $f : \RR^3 \to \RR$ be a smooth function, and let $(\mx_0, \mx_1, \mx_2)$ and $(\mx_0, \mx_1, \mx_3)$ define a pair of triangles which share the edge $[\mx_0, \mx_1]$ but are \emph{not} coplanar---e.g., see \Cref{fig:minimizers}. This gives a slightly more general setting from which to consider adjacent tetrahedron updates. Let $\mv_i = \mx_i - \mx_0$ ($1 \leq i \leq 3$) and define the restriction of $f$ to each triangle:
\begin{equation}
  f_1(\lambda_1, \lambda_2) = f(\mx_0 + \lambda_1\mv_1 + \lambda_2\mv_2), \qquad f_2(\lambda_1, \lambda_3) = f(\mx_0 + \lambda_1\mv_1 + \lambda_3\mv_3).
\end{equation}
The gradients of $f_1$ and $f_2$ are given by:
\begin{align}
    \nabla f_1(\lambda_1, \lambda_2) &= \begin{bmatrix} \mv_1 & \mv_2 \end{bmatrix}^\top \nabla f(\mx_0 + \lambda_1 \mv_1 + \lambda_2 \mv_2), \\
    \nabla f_2(\lambda_1, \lambda_3) &= \begin{bmatrix} \mv_1 & \mv_3 \end{bmatrix}^\top \nabla f(\mx_0 + \lambda_1 \mv_1 + \lambda_3 \mv_3).
\end{align}
To minimize $f$ over the first triangle, we solve:
\begin{equation}
  \begin{split}
    \mbox{minimize} \quad& f_1(\lambda_1, \lambda_2) \\
    \mbox{subject to} \quad& \lambda_1 \geq 0, \quad \lambda_2 \geq 0, \quad 1 - \lambda_1 - \lambda_2 \geq 0.
  \end{split}
\end{equation}
This is a nonlinear optimization problem with linear inequality constraints, analogous to the cost function for the tetrahedron update. Its Lagrangian is given by:
\begin{equation}
  L_1(\alpha_1, \beta_1, \gamma_1, \lambda_1, \lambda_2) = f_1(\lambda_1, \lambda_2) - \alpha_1 \lambda_1 - \beta_1 \lambda_2 - \gamma_1 (1 - \lambda_1 - \lambda_2).
\end{equation}
We find that the first-order necessary conditions for optimality are:
\begin{align}
  \pd{f_1}{\lambda_1} = \alpha_1 - \gamma_1, \qquad \pd{f_1}{\lambda_2} = \beta_1 - \gamma_1,
\end{align}
and the complementary slackness conditions are:
\begin{equation}
  \alpha_1 \lambda_1 = 0, \qquad \beta_1 \lambda_2 = 0, \qquad \gamma_1 (1 - \lambda_1 - \lambda_2) = 0.
\end{equation}
Now, assume that a minimizer for $f_1$ and $f_2$ exists on the interior of the edge
$[\mx_0, \mx_1]$. That is, there exists an optimum
$(\lambda_1^*, \lambda_2^*) = (\lambda_1^*, \lambda_3^*) =
(\lambda_1^*, 0)$ such that $0 < \lambda_1^* < 1$, where:
\begin{align}
  (\lambda_1^*, 0) = \Arg\min_{(\lambda_1, \lambda_2) \in \Delta^2} f_1(\lambda_1, \lambda_2) = \Arg\min_{(\lambda_1, \lambda_3) \in \Delta^2} f_2(\lambda_1, \lambda_3).
\end{align}
From the complementary slackness conditions, we know that
$\alpha_1^* = \gamma_1^* = 0$ since $\lambda_1^* \neq 0$. Likewise,
since $\lambda_2^* = 0$, we may have $\beta_1^* \neq 0$ in
general. This gives:
\begin{equation}
  \begin{bmatrix} 0 \\ \beta_1^* \end{bmatrix} = \nabla f_1(\lambda_1^*, 0) = \begin{bmatrix} \mv_1 & \mv_2 \end{bmatrix}^\top \nabla f(\mx_0 + \lambda_1^* \mv_1) = \begin{bmatrix} \mv_1^\top \nabla f(\mx_0 + \lambda_1^* \mv_1) \\ \mv_2^\top \nabla f(\mx_0 + \lambda_1^* \mv_1) \end{bmatrix}.
\end{equation}
Hence, the optimal Lagrange multiplier $\beta_1^*$ satisfies:
\begin{equation}
  \beta_1^* = \mv_2^\top \nabla f(\mx_0 + \lambda_1^* \mv_1).
\end{equation}
Along very similar lines, the corresponding Lagrange multiplier for $f_2$ satisfies:
\begin{equation}
  \beta_2^* = \mv_3^\top \nabla f(\mx_0 + \lambda_1^* \mv_1).
\end{equation}
Hence, we have found a global optimum of $f$ over $\conv(\mx_0, \mx_1, \mx_2) \cup \conv(\mx_0, \mx_1, \mx_3)$ which yields two nontrivial boundary minimizers over each subtriangle.

\end{document}